\documentclass[12pt,english]{article}

\usepackage{babel}
\makeatletter

\usepackage{amsmath,latexsym}
\usepackage{amsfonts}
\usepackage{amssymb}
\usepackage{fontenc}
\usepackage{textcomp}
\usepackage{graphicx}

\graphicspath{{fig/}}

\newtheorem{thrm}{Theorem}[section]
\newtheorem{prop}[thrm]{Proposition}

\newtheorem{cllry}[thrm]{Corollary}
\newtheorem{lmma}[thrm]{Lemma}
\newtheorem{remk}[thrm]{Remark}
\newtheorem{defn}[thrm]{Definition}

\numberwithin{equation}{section}

	\def\bbl{{\mathbb L}}
	\def\bbk{{\mathbb K}}
	
	\def\calB{{\mathcal B}}	
	
	\def\calM{{\mathcal M}}

	\def\calL{{\mathcal L}}
	
	\def\calJ{{\mathcal J}}
	\def\calL{{\mathcal L}}
	\def\calN{{\mathcal N}}

	\def\calV{{\mathcal V}}
	
	\def\calS{{\mathcal S}}
	
	\def\bbr{{\mathbb R}}
	\def\bbe{{\mathbb E}}
	\def\bbp{{\mathbb P}}

	\def\bfX{{\mathbf X}}
	
	\def\bfx{{\mathbf x}}

	\def\bfXt{{\mathbf{X}(t)}}

	\def\ld{{\stackrel{\frak {D}}{\rightarrow}}}

\makeatother
\begin{document}

\title{Nonparametric estimation for L\'evy processes 
with a view towards mathematical finance}
\author{Enrique Figueroa-L\'opez{}
\\
Department of Mathematics\\
Purdue University\\
150 N. University Street,\\ West Lafayette, IN 47906\\
jfiguero@math.purdue.edu \and 
Christian Houdr\'e\thanks{Research supported in part
by NSF and NSA grants. 
}\\
Laboratoire d' Analyse et de\\
Math\'ematiques Appliqu\'ees\\
CNRS UMR 8050\\
Universit\'e Paris XII\\
94010 Cr\'eteil Cedex, France\\
and\\
School of Mathematics\\
Georgia Institute of Technology\\ Atlanta, GA 30332-0160, USA \\
houdre@math.gatech.edu}
\maketitle

\begin{abstract} { 
   Nonparametric methods for the estimation of
    the L\'evy density of a L\'evy process $X$ are developed.
    Estimators that can be  written in terms of the
    ``jumps'' of $X$ are introduced, and so are discrete-data based
    approximations.  A model selection approach made up of two steps
    is investigated.  The first step consists in the selection of a
    good estimator from a linear
    model of proposed L\'evy densities, while the second
    is a data-driven selection of a linear model among a
    given collection of linear models.
    By providing lower bounds for the minimax risk
    of estimation over Besov L\'evy densities,
    our estimators are shown to achieve the ``best'' rate of
    convergence.  A numerical study for the case of histogram
    estimators and for 
    variance Gamma processes, models of key importance in risky
    asset price modeling driven by L\'evy processes, is presented.
   }
\end{abstract}
\section{Introduction}
The class of L\'evy processes is central to the theory of 
stochastic processes 
(see \cite{Sato} and \cite{Bertoin} for excellent monographs 
on the topic). 
Recently, new subclasses of L\'evy processes have been 
introduced and actively investigated 
mostly because of their relevance to mathematical finance.
Among the better known models are the \emph{variance Gamma model} 
of \cite{Madan1}, the \emph{CGMY model} of \cite{Madan}, 
and the \emph{generalized hyperbolic motion} of \cite{Barndorff:1998} 
and \cite{Eberlein:1995} (see also  \cite{Barndorff2} and 
\cite{Eberlein}). 
This phenomenon is not surprising if one brings to mind the traditional 
model for risky assets, namely the \emph{Black-Scholes model}. 
In this model the price $S(t)$ of an asset at time $t$ is assumed 
to be governed by 
\[
	S(t)=S(0)e^{\sigma B(t)+ \mu t},
\]
where $B(t)$ is a standard Brownian motion. 
However, a well-documented empirical evidence against the 
Black-Scholes model, specially in describing high-frequency data 
and option prices, have led researchers to consider non-Gaussian  
based models (see for instance \cite{Madan}, \cite{Madan2}, 
\cite{Barndorff}, \cite{Barndorff2}, \cite{Eberlein}, 
and references therein). 
The transition to L\'evy processes is natural since these 
preserve the statistical qualities of 
Brownian Motion's increments, but relax the path continuity
by allowing jump-alike discontinuities (a specification that is more 
consistent with the real evolution of stock prices through time).
Another point that naturally led to L\'evy processes was the use of  
economic-relevant random clocks (like volume or number of trades) 
instead of the historical time. 
Specifically, a more robust and sensible model is to take 
\begin{equation}\label{GBM}
	S(t) = S(0) e^{\sigma B(T(t)) +\mu T(t)},
\end{equation}
where $T(t)$ is a general increasing stochastic process with $T(0)=0$. 
In that case, if $T(t)$ has independent and stationary increments, 
the log return process is necessarily a L\'evy process 
(see 30.1 in \cite{Sato}). 
Such considerations led to the study  of 
\emph{exponential L\'evy processes} of the form 
\begin{equation}\label{GLP}
	S(t)=S(0) e^{X(t)}, 
\end{equation}
where $X(t)$ is a L\'evy process. 
This paradigm has proved to be successful to account for many of the 
empirical features of financial data. 
However, among other drawbacks, the high computational intensity and 
numerical issues involved in calibrating such models have prevented 
them from being more widely used in practice.  
In particular, these difficulties become very serious when
dealing with  ``high-frequency'' data.
 
L\'evy processes are determined by three ``parameters'': 
a non-negative real $\sigma^{2}$, a real $\mu$, and 
a measure $\nu$ on $\mathbb{R}\backslash \{0\}$. 
These three parameters characterize  the dynamic of a 
L\'evy process $\{ X(t) \}_{t\geq{0}}$
 as the superposition of a Brownian motion with drift, 
$\sigma B(t)+\mu t$, and a pure-jump
 L\'evy process whose jump behavior is specified 
by the measure $\nu$ as follows:
\[
	\nu(A)=\frac{1}{t}\bbe\left[\sum_{s\leq{t}} 
	\chi_{A}\left(\Delta{X}(s)\right)\right],
\]
where $\Delta{X}(t)\equiv{X}(t)-X(t^{-})$ is the jump of $X$ 
at time $t$ and $A$ is such that the indicator $\chi_{A}(\cdot)$ 
vanishes in a neighborhood of the origin
(this is a consequence of the so called L\'evy-It\^o 
decomposition for the sample paths of processes with 
independent increments; 
see Theorem 13.4. of \cite{Kallenberg} or Section 19 of \cite{Sato}).
We assume throughout that $\nu$ is determined by a function 
$p:\bbr\backslash\{0\}\rightarrow [0,\infty)$, called 
the \emph{L\'evy density}, in the following sense:
\[
	\nu(A)=\int_{A} p(x) dx, \;\; \forall 
	A\in\calB(\bbr\backslash\{0\}).
\]
In that case, the value of $p$ at $x_{0}$ provides, roughly speaking, 
information on the frequency of jumps with sizes ``close'' to $x_{0}$.

Estimating the L\'evy density poses a nontrivial problem, 
even when $p$ takes simple parametric forms. 
Parsimonious L\'evy densities usually produces not only 
intractable but sometimes not even expressible densities 
for the marginals $X(t)$, $t>0$.  
The current practice of estimation relies on approximations 
of the density function using \emph{inversion formulas} 
combined with likelihood methods (see for instance \cite{Madan}). 
Such approximations make the estimation particularly susceptible to
numerical errors and mis-specification; 
that is, slight changes in the model can produce quite different 
results. It is important to notice that these problems become quite 
critical for ``high-frequency'' data.
Other common calibration methods include 
simulation based methods and multinomial log likelihoods 
(see for instance \cite{Jiang} and \cite{Bibby}). 


In the present paper, we introduce new estimation methods for 
the L\'evy density.
We concentrate on model-free estimation schemes that allow to 
efficiently retrieve a fairly general L\'evy density. 
Being nonparametric, we relax the dependency on the model 
and expect that data itself validates the best model. 
Three theories serve as foundations for our methodology: 
i) the characterization of the jumps associated with a 
L\'evy process as a spatial Poisson process, 
ii) some recent methods for the nonparametric estimation of 
spatial Poisson processes introduced in \cite{Reynaud}, and 
iii) the short-term properties of L\'evy processes to 
approximate jump-dependent quantities.  
To the best of our knowledge, such connection between the 
L\'evy density and the statistical properties of the process 
in small time spans has not been used for calibration purposes 
before the present work. 
It is relevant to point out that our procedures are suitable 
for high-frequency data, which is widely available nowadays. 
Furthermore, it is precisely for such data that standard 
statistical estimation methods are not viable, 
the traditional \emph{geometric Brownian motion} model is 
totally inaccurate, 
and general exponential L\'evy models may be more relevant.


Let us describe the outline of the paper. 
In Section \ref{BscMethod}, we construct functional estimators 
which can be written in terms of integrals of deterministic functions  with respect to the random measure associated with the jumps of $X$. 
The proposed method follows the reasoning of the works on minimum 
contrast estimation on sieves and model selection developed 
in the context of density estimation and nonlinear regression 
in \cite{Birge:1994} (see \cite{Birge:1995} 
and \cite{Birge}) and recently extended to the estimation of 
intensity functions for Poisson processes 
in \cite{Reynaud}.
Concretely, the procedure addresses two problems: 
1) the selection of a good estimator, called the 
\emph{projection estimator}, from a linear model ${\cal S}$ of 
possible estimators, and 2) the selection of a linear model among 
a given collection of linear models using a penalization technique 
that led to a \emph{penalized projection estimator} (p.p.e.). 
A bound for the \emph{risk} of the p.p.e. is found in Section 
\ref{SectOracleIneq}.
As a consequence, \emph{Oracle inequalities}, 
that ensure to approximately reach the best expected error 
(using projection estimators) up to a constant,
are obtained. 
We also assess the rate of convergence of the p.p.e. on 
\emph{regular splines}, when the L\'evy density belongs 
to some \emph{Besov spaces}.  
By analyzing the \emph{minimax risk} of estimation on these Besov spaces, 
it is actually proved in Section \ref{SectionMinimaxRisk}
that the p.p.e. attains the best possible rate in the minimax sense,
when the estimation is based on jumps bounded away from the origin. 
In Sections \ref{SecAproxPoissnInt} and \ref{EstimationmMethod},
 we examine the problem that the Poisson jump measure 
cannot be retrieved  from discrete observation, and devise an approximation procedure for Poisson integrals based on equally space sampling observations of the process. 
Finally, in the last part our methods are applied to the estimation of 
a classical model used in mathematical finance:
the \emph{Variance Gamma model} of \cite{Madan1}.
The L\'evy processes are simulated using time series representations 
and ``discrete skeletons'',
whereas the considered estimators are mainly regular histograms.


\section{A model-free estimation method}\label{BscMethod}
Consider  a real L\'evy  process 
$X=\left\{ X(t) \right\}_{t   \geq  0}$ with
 L\'evy density
\index{L\'evy processes!L\'evy density}\index{L\'evy density}  
$p$. That  is,  $X$   is a  c\`adl\`ag   process  with
independent and  stationary   increments such  that 
the characteristic function of its marginals is  given by  
  \begin{equation}\label{CharcFuntChap3}  
  	\bbe\left[e^{iuX(t)}\right]=
  	\exp\left\{t\left(iub-\frac{u^{2} \sigma^{2}}{2}+
  	\int_{\bbr_{0}}\left\{e^{iux}-1-iux1_{[\vert  x \vert
  \leq1]}\right\}p(x)dx\right)\right\}, 
  \end{equation}
  where $\mathbb{R}_{0}=\mathbb{R}\backslash   \{0\}$
  and $p:\mathbb{R}_{0}\rightarrow\bbr_{+}$ satisfies
  \begin{equation}
    \label{LevyDensity}
    \int_{\bbr_{0}} (1 \wedge x^{2})p(x) dx < \infty.
  \end{equation}
  Being a c\`adl\`ag process, the set of jump times
 \(
	\left\{t>0: X (t)-X(t^{-})>0\right\}
 \)
 is countable and, for Borel subsets $B$ of $[0,\infty) \times \mathbb{R}_{0}$,
\begin{equation}\label{JumpMeasure}
	{\cal J}(B) \equiv \# 
	\left\{t>0: (t,X (t)-X(t^{-})) \in B \right\},
\end{equation}
is a well-defined random measure on 
$[0,\infty) \times \mathbb{R}_{0}$, 
with $ \# $ denoting cardinality. 
The L\'evy-It\^o decomposition of the sample paths
\index{L\'evy-It\^o decomposition} (see Theorem 19.2 of
\cite{Sato}) implies that ${\cal J}$ is a Poisson process 
on the Borel sets of ${\cal B}([0,\infty) \times \mathbb{R}_{0})$ 
with mean measure given by
\begin{equation}\label{MeanMeasure}
	\mu (B) = \iint\limits_{B} p(x)\, d t\, d  x.
\end{equation}

We study the problem of estimating the L\'evy density $p$ 
on a Borel set $D\in {\cal B}\left(\bbr_{0}\right)$ using 
a \emph{projection estimation approach}. 
According to this paradigm, $p$ is estimated by estimating
the best approximating function in a finite-dimensional 
linear space $\calS$.
The linear space $\calS$ is taken so that it has good 
approximation qualities in general classes of functions. 
Typical choices are piecewise polynomials or wavelets. 
In order for this approach to be general enough but still 
feassible, it is usually assumed that the function to be estimated 
is bounded and belongs to an $\bbl^{2}$ space on $D$, simplying
the task of specifying the best approximating function. 
The simplest case is when $p$ is taken bounded and 
\(	
	\int_{D} p^{2}(x) dx < \infty.	
\) 
This condition is quite general if $D$ is away from the origin, 
since (\ref{LevyDensity}) entails 
\begin{equation}
\label{Prop1LevyDnsty}
\int_{\vert x \vert > \varepsilon} p^{2}(x) dx < \infty,
\end{equation}
for any $\varepsilon >0$, when $p$ is bounded on $\{ x:
|x|>\varepsilon\}$. However, around the origin the 
L\'evy density is not bounded in most applications.
This motivates the use of measures different from the Lebesgue
measure. Concretely, it is assumed that the L\'evy measure 
$\nu(dx)\equiv p(x)dx$ is absolutely continuous with respect 
to a known measure $\eta$ on ${\cal B}\left(D\right)$ and that 
the Radon-Nikodym derivative
\begin{equation}\label{RegDnsty}
	\frac{d\nu}{d\eta}\left(x\right) = s(x), \;\;x \in D, 
\end{equation}
is positive, bounded, and satisfies
\begin{equation}\label{CondtRegDnsty}
	\int_{D} s^{2}(x) \eta(dx) < \infty.
\end{equation}
\noindent \begin{defn}\label{Regularization} If (\ref{RegDnsty}) and (\ref{CondtRegDnsty}) are verified, 
we say that $\eta$ is a \emph{regularizing measure}
\index{regularizing measure} for the L\'evy density $p$. 
In that case, $s$ is referred to as the 
\emph{regularized (under $\eta$) L\'evy density}
\index{regularized L\'evy density} of $p$ on $D$.
\end{defn}
Notice that under the previous regularization assumption, 
the measure ${\cal J}$ of  (\ref{JumpMeasure}), 
when restricted to ${\cal B}([0,\infty) \times D)$, 
is a Poisson process with mean measure
\begin{equation}
\label{RegMeanMeasure}
\mu (B) = \iint\limits_{B} s(x) \,d t\,\eta(dx), \;\; B\in {\cal B}([0,\infty) \times D).
\end{equation} 
Our goal will be to estimate the regularized L\'evy density $s$, 
and  using (\ref{RegDnsty}) to retrieve $p$ on $D$ from $s$. 
To illustrate this strategy consider a continuous L\'evy density 
$p$ such that 
\[
p(x)=O\left(x^{-1}\right),\;\;{\rm as}\;\; x\rightarrow 0. 
\] 
This type of densities admit the regularizing measure 
$\eta(dx)=x^{-2} dx$ on domains of the form 
$D=\{ x : 0<|x|<b\}$.
Indeed, $s(x)= x^{2} p(x)$ will be bounded and fulfills 
(\ref{CondtRegDnsty}).  
Clearly,  each estimator $\hat{s}$ for $s$ will induce 
the natural estimator $x^{-2}\hat{s}(x)$ for $p$.

The previous methodology is motivated by recent 
results on the estimation of intensity functions of 
non-homogeneous Poisson processes (see \cite{Reynaud}). 
In that  paper, a type of  \emph{projection
estimator} is proposed,
whereas \emph{penalized projection estimation}\index{penalized 
projection estimation}
is used as a data-driven criterion for selecting the
best space among a family of linear spaces. 
However, these procedures focus on finite Poisson point processes and 
on classes of intensity functions that are defined with 
respect to a finite reference measure 
(see Section \ref{SectOracleIneq} for a more detailed description of this hypothesis).
Actually, the value of the reference measure plays a key role 
in the definitions of projection estimators and penalization. 
Our job in this section is to implement and justify a 
projection estimation approach that does not rely on the 
finiteness of the Poisson process.

Let us describe the main ingredients of our procedure. 
Consider the random functional 
\begin{equation}\label{Contrast}
	\gamma_{D}(f)\equiv -\frac{2}{T} 
	\iint\limits_{[0,T] \times D} f(x)\, 
	{\cal J}(dt,dx) + \int\limits_{D} f^{2}(x)\, \eta(dx),
\end{equation}
which is well defined for any function 
$f \in L^{2}\left((D,\eta)\right)$, where 
$D \in {\cal B}\left(\bbr_{0}\right)$ and $\eta$ is as in equations (\ref{RegDnsty})-(\ref{RegMeanMeasure}). 
Following the terminology of \cite{Birge:1994} and \cite{Reynaud}, 
we call $\gamma_{D}$ the \emph{contrast function}.
\index{contrast function} 
Throughout this section,
\[
	\| f \|^{2} \equiv  \int_{D} f^{2}(x)\, \eta(dx),
\]
for any $f \in L^{2}\left((D,\eta)\right)$.  
Let ${\cal S}$ be  a finite dimensional subspace of  
$L^{2}\equiv L^{2}\left((D,\eta)\right)$. The \emph{projection estimator}\index{projection estimator}\index{orthogonal projection!estimator of} of $s$ on $S$ is defined by 
\begin{equation}
	\label{DefPE}
	\hat s(x) \equiv \sum_{i=1}^{d} \hat\beta_{i} \varphi_{i} (x),
\end{equation}
where $\left\{ \varphi_{1}, \dots, \varphi_{d} \right\}$ is an arbitrary orthonormal
basis of ${\cal S}$ and 
\begin{equation}
\label{DefCoefPE}
\hat\beta_{i} \equiv \frac{1}{T}\iint\limits_{[0,T] \times D} \varphi_{i}(x) {\cal
J}(dt,dx).
\end{equation}
Let us give another characterization of the projection estimator.
	\noindent \begin{remk}
	\label{Remark_1}\index{contrast function!relation with projection estimation}
	The projection estimator is the unique minimizer of the contrast function $\gamma_{D}$ over
	$S$. Indeed, plugging $f=\sum_{i=1}^{d} \beta_{i} \varphi_{i}$ in (\ref{Contrast}) gives 
	$\gamma_{D} (f)=\sum_{i=1}^{d} \left( -2\beta_{i}\hat\beta_{i} + \beta^{2}_{i}\right)$, 
	and thus, $\gamma_{D} (f) \geq -\sum_{i=1}^{d} \hat\beta_{i}^{2}$, for all $f \in S$.
	In particular,  this characterization implies that $\hat s$ does not depend on the 
	choice of the orthonormal basis, and suggests a mechanism to numerically approximate
	$\hat{s}$ when we do not have an explicit orthonormal basis for ${\cal S}$. 
\end{remk}
The remark above helps to make sense of $\hat s$ as an estimator
 of the regularized L\'evy density $s$ because  the minimizer of $\bbe \left[ \gamma_{D}(f) \right]$ over all $f\in {\cal S}$ is precisely the closest function in ${\cal S}$ to $s$. Concretely, the \emph{orthogonal projection}\index{orthogonal projection} of $s$ on the subspace ${\cal S}$, namely 
 \begin{equation} \label{OrthProj}
 	s^{\bot} \equiv \sum_{i=1}^{d} \left(\int_{D} \varphi_{i}(y) s(y) \eta(dy)\right)\varphi_{i} (x),
\end{equation}
is such that  
\begin{equation}
-\| s^{\bot} \|^{2} = \bbe\left[ \gamma_{D}(s^{\bot}) \right]  \leq \bbe \left[
\gamma_{D} (f) \right],\;\;\;\forall f\in{\cal S}.
\end{equation}
Moreover, we can readily corroborate that $\hat s$ is an unbiased estimator\index{projection estimator!mean} of the  orthogonal projection $s^{\bot}$. In order to assess the quality of estimation, we compute the ``square error'' of
$\hat s$:
\begin{equation}
\label{Chi}
\chi^{2} \equiv \| s^{\bot} - \hat s \|^{2} = 
\sum_{i=1}^{d} \left[\;\iint\limits_{[0,T] \times D} \varphi_{i}(x) 
 \frac{{\cal J}(dt,dx)-s(x)\,dt\,\eta(dx)}{T} \;\right]^{2}.
\end{equation}
Then, by the standard formula for the variance of Poisson integrals, 
the mean square\index{projection estimator!mean-square error} error takes the form
\begin{equation}
\label{RiskPEonOP}
\bbe\left[ \chi^{2} \right] =
\frac{1}{T} \sum\limits_{i=1}^{d} \int\limits_{D} \varphi^{2}_{i}(x)s(x)\,\eta(dx).
\end{equation}
The quantity $ \bbe\left[ \chi^{2} \right]$ is called the \emph{variance
term}\index{variance term} and the equation above shows that this term will shrink to $0$ when the time horizon $T$ goes to infinity.  
Moreover, the \emph{risk}\index{risk of an estimator}\index{projection estimator!risk of}  of 
$\hat s$, $\bbe \left[ \| s - \hat s \|^{2}\right]$, 
can be decomposed into a nonrandom term plus the previous variance term: 
\begin{equation}
\label{RiskPEonR}
\bbe \left[ \| s - \hat s \|^{2}\right]= \|s-s^{\bot}  \|^{2} 
+ \bbe \left[\chi ^{2}\right].
\end{equation}
The first term, called the \emph{bias term},\index{bias term} accounts for the distance of
the unknown function $s$ to the model ${\cal S}$ and does not depend on the
estimation criteria we use within the model.

The next natural problem to tackle is to design a data-driven scheme for selecting a ``good'' model\index{model selection problem} from a collection of linear models $\left\{ {\cal S}_{m},
m\in {\cal M} \right\}$. Namely, we wish to select a model that
approximately realizes the best trade-off between the risk
of estimation within the model and the distance of the unknown L\'evy
density to the model. Let $\hat s_{m}$ and $s^{\bot}_{m}$ be respectively the
projection estimator and the orthogonal projection of $s$ on ${\cal
S}_{m}$. For each $m\in {\cal M}$, let $\chi^{2}_{m}$ be as in (\ref{Chi}). 
The following simplifications of (\ref{RiskPEonR})  give insight on a possible solution:
\begin{align}
\bbe \left[ \| s - \hat s_{m} \|^{2}\right] &=\|s-s_{m}^{\bot}  \|^{2} 
+ \bbe \left[ \chi_{m}^{2} \right] \nonumber \\
\label{SimplifyRisk}
&= \|s \|^{2} - \|s_{m}^{\bot} \|^{2} 
+ \bbe \left[ \chi_{m}^{2} \right] \\
&= \|s \|^{2} - \bbe \left[ \|\hat s_{m} \|^{2}\right] 
+ 2\bbe \left[ \chi_{m}^{2} \right] \nonumber\\
&= \|s \|^{2} + \bbe\left[ \gamma_{D}\left( \hat s_{m} \right) 
+ {\rm pen}(m) \right], \nonumber
\end{align}
where ${\rm pen}(m)$ is defined in terms of an orthonormal basis 
$\left\{ \varphi_{1,m}, \dots ,\varphi_{d_{m},m} \right\}$ for ${\cal S}_{m}$
by the equation: 
\begin{equation}
\label{Penalty1}
{\rm pen}(m)=\frac{2}{T^{2}} \iint\limits_{[0,T] \times D} 
\;\left(\sum_{i=1}^{d_{m}} \varphi^{2}_{i,m}(x) \right) {\cal J}(dt,dx).
\end{equation}
Equation (\ref{SimplifyRisk}) shows that the risk of $\hat
s_{m}$ moves ``parallel'' to the expectation of the \emph{observable statistics}
$\gamma_{D}\left( \hat s_{m} \right) + {\rm pen}(m)$. This
fact heuristically justifies to choose the model that minimizes such a 
penalized contrast value. 
In general, 
 it makes sense to consider \textbf{penalized projection estimators} (p.p.e.)
 \index{penalized projection estimator} of the form 
\begin{equation}
\label{PPE}
\tilde s \equiv \hat s_{\hat m},
\end{equation}
where ${\rm pen}:{\cal M}\rightarrow [0,\infty)$, $\hat s_{m}$ is the projection estimator on ${\cal S}_{m}$ (see  (\ref{DefPE})), and 
\(
\hat m \equiv {\rm argmin}_{m \in {\cal M}} 
\left\{ \gamma_{D}\left( \hat s_{m} \right) + {\rm pen}(m) \right\}.
\)

Methods of estimation based on the minimization of penalty functions have a long history in the literature on regression and density estimation (for instance,  \cite{Akaike}, \cite{Mallows}, and \cite{Schwarz}). The general idea is to choose among a given collection of parametric models the model that minimizes a loss function plus a penalty term that controls the complexity of the model. Such penalized estimation was promoted for nonparametric density estimation  in \cite{Birge}, and in the context of non-homogeneous Poisson processes in \cite{Reynaud}. There are two main accomplishments obtained in these works both in the context of density estimation and intensity estimation of nonhomogeneous Poisson processes: Oracles inequalities and competitive performance against minimax estimators. 
The following section shows that the method outlined above preserves Oracle inequalities. 


\section{Risk bounds, oracle inequalities, and 
rates of convergence}\label{SectOracleIneq}
Consider the problem of model selection among a collection of linear models, 
$\left\{{\cal S}_{m}, m \in {\cal M} \right\}$, for the regularized L\'evy density $s$ 
on $D$ as outlined in the previous section. We showed through
(\ref{SimplifyRisk}) that a sensible criterion to decide for a projection
estimator is to penalize its contrast value with a properly chosen penalty function 
${\rm pen}: {\cal M}\rightarrow [0,\infty)$. Of course, the
``best'' model, namely 
	\begin{equation}\label{OracleModel}
		\bar m \equiv {\rm argmin}_{m\in {\cal M}}\; \bbe \left[ \| s - \hat
		s_{m} \|^{2}\right], 
	\end{equation}
is not accessible, but we can aspire to achieve the smallest possible
risk up to a constant. In other words, it is desirable that our estimator
$\tilde s$ comply with an inequality of the form  
	\begin{equation}
	\label{OracleIneqFrst}
		\bbe \left[ \| s - \tilde s \|^{2}\right] \leq 
		C  \inf_{m\in {\cal M}} \bbe \left[ \| s - \hat s_{m} \|^{2}\right],
	\end{equation}
for a constant $C$ independent of the linear models. The model that achieves the minimal risk 
of projection estimation is called the \emph{Oracle model}\index{Oracle model} and 
inequalities of the type (\ref{OracleIneqFrst}) are called 
\emph{Oracle inequalities}\index{Oracle inequality}. Approximate Oracle inequalities were proved
in \cite{Reynaud} for the intensity function of a nonhomogeneous Poisson
process $\left\{ N_{A} \right\}_{A\in{\cal V}}$ on a measurable space $({\rm
V},{\cal V})$. Concretely, \cite{Reynaud} defines projection estimators $\hat{s}_{m}$ 
and penalized projection estimators $\tilde{s}$ satisfying
\begin{equation}\label{OracleIneqReynaud}
	\bbe \left[ \int_{V} \left| s({\bf v}) - \tilde{s}({\bf v}) \right|^{2}\frac{\zeta(d{\bf v})}{\zeta(V)} \right]
	 \leq C  \inf_{m\in {\cal M}} 
	 \bbe \left[  \int_{V} \left| s({\bf v}) - \hat{s}_{m}({\bf v}) \right|^{2}\frac{\zeta(d{\bf v})}{\zeta(V)} \right] 
	 +\frac{C'}{\zeta(V)},
\end{equation}
where $s$ and $\zeta$ are respectively a bounded measurable function and 
a finite measure on $V$ such that 
	\[
		\bbe[N_{A}]=\int_{A} s({\bf v}) d\zeta({\bf v}),\;\; A\in\calV.
	\]
The finiteness of $\zeta$ plays an important role in the definition of the 
estimators, and in obtaining the Oracle inequality  
(\ref{OracleIneqReynaud}). 
However,  such a property is not necessarily satisfied by 
the mean measure of the Poisson process ${\cal J}(\cdot)$ of (\ref{JumpMeasure}) on 
${\cal B}([0,T]\times D)$ (for instance, if  
$D=\left\{ \vert x \vert>\varepsilon\right\}$ under $\zeta(d{\bf v})=dx\,dt$ as in (\ref{MeanMeasure}), 
or if $D=\left\{ 0<\vert x \vert<b\right\}$ and $\zeta(d{\bf v})=x^{-2}dx\,dt$ as in the example described after Definition \ref{Regularization}). 
In this section we show that, based on one sample of the L\'evy process $X$ on $[0,T]$,  
the projection estimators $\{\hat{s}_{m}\}_{m\in {\cal M}}$ introduced in Section \ref{BscMethod},
and certain penalized projection estimators \(\tilde{s}\) 
satisfy the approximate Oracle inequality
\[
	\bbe \left[ \| s - \tilde s \|^{2}\right] \leq 
	C  \inf_{m\in {\cal M}} \bbe \left[ \| s - \hat s_{m} \|^{2}\right] +\frac{C'}{T},
\]
where $s$ is a regularized L\'evy density, and the constants $C, C'$ depend only on the 
``complexity'' of the family of linear models. 
Actually, we will be able to estimate the order of the constants $C$ and $C'$ appearing in the Oracle inequality. 

The main tool in obtaining Oracle inequalities is an upper 
bound for the risk of the penalized projection estimator 
$\tilde{s}$ of  (\ref{PPE}). 
The proof of this bound is a simple variation of the 
argument of \cite{Reynaud}; 
however, to overcome the possible lack of finiteness on 
$\zeta$ and to avoid unnecessary use of upper bounds, 
the dimension of the linear model is explicitly included
in the penalization. 
Finally, the obtained risk bound is used to assess the 
rate of convergence of $\tilde{s}$ to $s$ in the long run 
(as $T\rightarrow\infty$) when $s$ is ``smooth'' and 
the considered linear spaces are piecewise polynomials.

The following regularity condition was introduced in \cite{Reynaud} to make a distinction between not too ``large'' families of linear models and certain wavelet-type linear models. We will focus here on the simplest case: 
\noindent \begin{defn}\label{Polynomial}\index{polynomial collection} 
	A collection of models $\left\{{\cal S}_{m}, m \in {\cal M} \right\}$ is said to be polynomial if there 
	exist constants $\Gamma>0$ and $R \geq 0$ such that for every positive integer $n$
	\[
		\#\left\{ m \in {\cal M}: d_{m}=n \right\} \leq \Gamma n^{R},
	\]
	where $d_{m}$ stands for the dimension of the model $\calS_{m}$, while $\#$
	denotes cardinality.
\end{defn}
Below, we return to the setting of Section \ref{BscMethod}; 
that is to say, $X=\left\{ X(t) \right\}_{0\leq{t}\leq{T}}$ is a L\'evy process with L\'evy density $p$ and  \emph{regularized} L\'evy density $s$ on a domain $D\in{\cal B}\left(\bbr_{0}\right)$ under a 
\emph{regularizing measure} $\eta$ (see Definition \ref{Regularization}). Define also
	\begin{equation}\label{DefOfDm}
		D_{m}= \sup\left\{ \| f\|_{\infty}^{2}: f\in S_{m}, \| f\|^{2} \equiv  \int_{D}
		f^{2}(x)\eta(dx)=1 \right\}.
	\end{equation}
\noindent\begin{remk}\label{DmInOB} 
	If $\left\{ \varphi_{1,m}, \dots,\varphi_{d_{m},m} \right\}$ is an arbitrary orthonormal basis of 
	${\cal S}_{m}$, then  $D_{m}=\|\sum_{i=1}^{d_{m}}\varphi^{2}_{i,m}\|_{\infty}$ 
	(see Section \ref{proofsChap3} for a verification).
\end{remk}
We now present the main result of this section 
(see Section \ref{SectionProofOracle} for the proof):
\noindent\begin{thrm}\label{PreOracleIneq}\index{penalized projection estimation!bound for the risk}
	Let $\left\{{\cal S}_{m}, m \in {\cal M} \right\}$ 
	be a polynomial family of finite dimensional linear subspaces 
	of $L^{2}((D,\eta))$ and let 
	${\cal M}_{T}\equiv\{m\in\calM: D_{m}\leq T\}$.
	If $\hat s_{m}$ and $s^{\bot}_{m}$ are respectively the projection estimator 
	and the orthogonal projection of the regularized L\'evy density $s$ 
	on ${\cal S}_{m}$ then, the penalized projection estimator $\tilde s_{_{T}}$ 
	on $\left\{\calS_{m}\right\}_{m\in\calM_{T}}$ defined by (\ref{PPE}) is such that 
	\begin{equation}
	\label{PseudoOracleIneq}
		\bbe \left[ \| s - \tilde s_{_{T}} \|^{2}\right] \leq 
		C  \inf_{m\in {\cal M}_{T}} \left\{ \| s - s^{\bot}_{m}\|^{2} +
		\bbe \left[ {\rm pen}(m)\right]\right\}+\frac{C'}{T},
	\end{equation}
	whenever ${\rm pen}: {\cal M} \rightarrow [0,\infty)$ takes one of the following 
	forms for some fixed (but arbitrary) constants $c>1$, $c'>0$, and $c''>0$:
	\item{\textbf{(a)}} 
	${\rm pen}(m) \geq c \frac{D_{m} {\cal N}}{T^{2}}
	+c' \frac{d_{m}}{T}$, where 
	${\cal N}\equiv{\cal J}([0,T] \times D)$ is the number of 
	jumps prior to $T$ with sizes falling in $D$ and where 
	it is assumed that $\rho \equiv \int_{D} s(x) \eta(dx)<\infty$;
	\item{\textbf{(b)}} 
	${\rm pen}(m) \geq c \frac{\hat V_{m}}{T}$, where 
	$\hat V_{m}$ is defined in terms of an orthonormal basis 
	$\left\{ \varphi_{i,m}\right\}_{i=1}^{d_{m}}$ of ${\cal S}_{m}$ by
	\begin{equation}\label{DefV}
		\hat V_{m}\equiv \frac{1}{T}\iint\limits_{[0,T] \times D} 
		\;\left(\sum_{i=1}^{d_{m}} \varphi^{2}_{i,m}(x) \right) 
		{\cal J}(dt,dx),
	\end{equation}
	and where it is assumed that 
	$\beta \equiv \inf_{m\in {\cal M}} 
	\frac{\bbe\left[ \hat V_{m}\right]}{D_{m}}>0$ and that 
	$\phi \equiv \inf_{m\in {\cal M}}\frac{D_{m}}{d_{m}}>0$;
	\item{\textbf{(c)}} ${\rm pen}(m) \geq c \frac{\hat V_{m}}{T}+
		c'\frac{D_{m}}{T}+c''\frac{d_{m}}{T}$.\\
	Moreover, the constant $C$ depends only on $c$, $c'$ and $c''$, while 
	$C'$ varies with $c$, $c'$, $c''$, $\Gamma$, $R$, $\| s \|$, $\| s \|_{\infty}$,
	$\rho$, $\beta$, and $\phi$.
\end{thrm}

\noindent\begin{remk} \label{OnTheConstants}
In the Remark \ref{Constants}, the order of the 
constants $C$ and $C'$ is analyzed. We will show that 
for $c\geq{2}$ and for arbitrary $\varepsilon>0$,
there is a constant $C'(\varepsilon)$ (increasing) so that 
\begin{equation}
	\bbe \| s - \tilde s \|^{2} \leq 
	\left( 1+\varepsilon\right) \inf_{m\in {\cal M}} 
	\left\{ \| s - s^{\bot}_{m}\|^{2} +
	\bbe \left[ {\rm pen}(m)\right]\right\}
	+\frac{C'(\varepsilon)}{T}.
\end{equation}
\end{remk} 

As a first use of the previous risk bound, we obtain Oracle inequalities for our p.p.e.
The next corollary immediately follows from the first equality in 
(\ref{SimplifyRisk}), equation (\ref{RiskPEonOP}), and part (b) above:
\noindent\begin{cllry} \label{OracleIneq}\index{Oracle inequality}
	\index{penalized projection estimation!Oracle inequality}
	In the setting of Theorem \ref{PreOracleIneq}, 
	if the penalty function is of
	the form ${\rm pen}(m) \equiv c \frac{\hat V_{m}}{T}$,
	for every $m\in {\cal M}_{T}$, $\beta>0$, and $\phi >0$, then
	\begin{equation}\label{OracleIneq}
		\bbe \left[ \| s - \tilde s_{_{T}} \|^{2}\right] \leq 
		\tilde{C} \inf_{m\in {\cal M}_{T}} 
		\left\{ \bbe\left[ \| s - \hat{s}_{m}\|^{2} 
		\right]\right\}+\frac{\tilde{C}'}{T},
	\end{equation}
	for a constant $C_{1}$ depending only on $c$, and a constant $C_{2}$ depending on
	$c$, $\Gamma$, $R$, $\| s \|$, $\| s \|_{\infty}$, $\beta$, and $\phi$.
\end{cllry} 

As a second application of (\ref{PseudoOracleIneq}), we analyze the 
``long run'' ($T\rightarrow\infty$) rate of convergence of penalized 
projection estimators on regular piecewise polynomials, 
when the L\'evy density is ``smooth''.
More precisely,  restricted to the window of estimation 
$D\equiv[a,b]\subset\bbr_{0}$,  the L\'evy density $s$ is
assumed to belong to the
\emph{Besov space}\footnote{
These Besov spaces are also called \emph{Lipschitz} or \emph{H\"older} spaces.}   
$\mathcal{B}^{\alpha}_{\infty}\left(\bbl^{p}([a,b])\right)$
with  some $p\in[2,\infty]$ and $\alpha>0$
(see Section 2.9-10 of \cite{Devore} for the definition). 
An important reason for the choice of this class of functions is
the availability of estimates for the error of approximation by 
\emph{splines}\footnote{
Piecewise polynomial functions $f$ such that on each compact 
interval, 
$f$ is made up of only finitely many polynomial pieces.}, 
trigonometric polynomials, and wavelet expansions 
(see for instance Chapter 12 of \cite{Devore}, 
and Lemma 13 of \cite{Birge:1995}). 
In particular, if $\calS^{k}_{m}$ denotes the space of 
piecewise polynomials of degree bounded by $k$, based on 
the regular partition of $[a,b]$ with $m$ pieces ($m\geq{1}$), 
Theorem 12.2.4 in \cite{Devore} implies that 
for any 
$s\in\mathcal{B}^{\alpha}_{\infty}\left(\bbl^{p}([a,b])\right)$ 
with $k>\alpha-1$, there exists a constant $C(s)$ such that
\begin{equation}\label{ApproxBesovIneq}
	{\rm d}_{p}\left(s,\calS^{k}_{m}\right)\leq 
	C(s) m^{-\alpha},
\end{equation}
where ${\rm d}_{p}$ is the distance induced by the 
$\bbl^{p}$-norm on $([a,b],dx)$. 
Actually, $C(s)$ can be taken to be increasing on   
$|s |_{\mathcal{B}^{\alpha}_{\infty}\left(\bbl^{p}\right)}$, 
the standard seminorm on 
$\mathcal{B}^{\alpha}_{\infty}\left(\bbl^{p}([a,b])\right)$
(see (10.1) Chapter 2 in \cite{Devore}). Combining 
(\ref{ApproxBesovIneq}) with (\ref{PseudoOracleIneq}), 
we obtain the following result (see Section \ref{proofsChap3}
for a proof).
\begin{cllry} \label{ConvergenceBesov}
Let $D\equiv [a,b]\subset\bbr_{0}$ and let $\calS_{m}^{k}$
be the space of piecewise 
polynomials of degree at most $k$ based on the regular 
partition of $[a,b]$ with 
$m$ pieces ($m\geq{1}$). Following the notation of Theorem  \ref{PreOracleIneq},
let $\tilde{s}_{_{T}}$ be the penalized projection estimator on 
$\{ \calS_{m}^{k}\}_{m\in\calM_{T}}$ with penalization
${\rm pen}(m)\equiv c\frac{\hat{V}_{m}}{T}+c'\frac{D_{m}}{T}+c''\frac{d_{m}}{T}$ 
(for some fixed $c>1$ and $c',c''>0$). Then, if the restriction of 
the L\'evy density $s$ to $[a,b]$ is a member of 
$\mathcal{B}^{\alpha}_{\infty}\left(\bbl^{p}([a,b])\right)$ with 
$2\leq{p}\leq\infty$ and $0<\alpha<k+1$, then
\[
	\limsup_{T\rightarrow\infty}T^{2\alpha/(2\alpha+1)}
	\bbe \left[ \| s - \tilde s_{_{T}} \|^{2}\right] <\infty.
\]
Moreover, for any $R>0$ and $L>0$,
\begin{equation}\label{LongRunRiskppe}
	\limsup_{T\rightarrow\infty}T^{2\alpha/(2\alpha+1)}
	\sup_{s\in \Theta(R,L)}\bbe 
	\left[ \| s - \tilde s_{_{T}} \|^{2}\right] <\infty,
\end{equation}
where $\Theta(R,L)$ consists of all L\'evy densities $s$ such that 
$\|s\|_{\bbl^{\infty}([a,b])}<R$, 
and $s$ restricted to $[a,b]$ is a member of  
$\mathcal{B}^{\alpha}_{\infty}\left(\bbl^{p}([a,b])\right)$
with seminorm 
$| s |_{\mathcal{B}^{\alpha}_{\infty}\left(\bbl^{p}\right)}<L$.

\end{cllry}
The previous result implies that the p.p.e. on regular splines 
has a rate of convergence of order $T^{-2\alpha/(2\alpha+1)}$
for the class of Besov L\'evy densities $\Theta(R,L)$. 
We will see in the next section  that the rate cannot be improved 
(see Corollary \ref{MinimaxResult3} and Remark \ref{MinimaxMore}).


\section{On the minimax risk for the estimation of 
smooth L\'evy densities}\label{SectionMinimaxRisk}
This section presents some results on the \emph{minimax risk} 
of estimation for certain families of smooth L\'evy densities. 
Roughly speaking, a \emph{minimax risk} on a given family 
$\Theta$ of  ``parameters'' has the following general form:
\[
	\inf_{\hat{s}}\sup_{{s}\in{\Theta}}\bbe_{s}
	\left[d\left({s},\hat{s}\right)\right],
\]
where the $\inf$ is taken over all the estimators $\hat{s}$ 
(based on the available random data, whose law distribution 
is itself determined by the parameter ${s}$), 
 and $d(s,\hat{s})$ is a function that measures how distant 
 ${s}$ and $\hat{s}$ are from each other. 
In some sense, 
$\sup_{{s}\in{\Theta}}\bbe_{s}
\left[d\left({s},\hat{s}\right)\right]$ 
measures the maximum error that can arise when using the estimator 
$\hat{s}$. Therefore, an estimator that approximately accomplishes 
a minimax risk is desirable. 
Comparisons to the minimax risks is one of the most solicited 
measures of performance in statistical estimation. 
In fact, minimax type results have been obtained in very general 
contexts (see for instance \cite{Ibragimov} and \cite{Birge:1995} 
in the case of density estimation based on i.i.d. random variables, 
and \cite{Kutoyants} and \cite{Reynaud} in the case of intensity 
estimation based on finite Poisson point processes).

Since the jumps of a L\'evy process can be associated with a 
Poisson point process on $\bbr_{+}\times\mathbb{R}\backslash\{0\}$, 
many results and techniques for the statistical inference 
of Poisson processes can be translated into  the context of 
L\'evy processes. Following this approach, 
we adapt below a result of Kutoyants \cite{Kutoyants} 
(Theorem 6.5) on the asymptotic minimax risk for the 
estimation of ``smooth'' intensity functions of a Poisson point 
processes on $[0,1]$, based on $n$ independent copies. 
The idea of the proof is due to 
Ibragimov and Has'minskii and is based on the statistical 
tools for distributions satisfying 
the \emph{Local Asymptotic Normality} (LAN) property 
(see Chapters II and Section IV.5 of \cite{Ibragimov}).
Some generalizations and consequences are also deduced.

Let us introduce a \emph{loss function} 
$\ell:\bbr\rightarrow\bbr$ with the following properties:
\begin{itemize}
	\item[(i)] $\ell(\cdot)$ is nonnegative, 
	$\ell(0)=0$ but not identically $0$, and continuous at $0$;
	\item[(ii)] it is symmetric: $\ell(u)=\ell(-u)$ for all $u$;
	\item[(iii)] for any $c>0$, $\{u:\ell(u)<c\}$ is a convex set;
	\item[(iv)] $\ell(u)\exp\{\varepsilon|u|^{2}\}\rightarrow{0}$ 
	as $|u|\rightarrow\infty$, for any $\varepsilon>0$.
\end{itemize}
Consider the problem of estimating the L\'evy density $s$ 
of a L\'evy process $\left\{ X(t) \right\}_{0\leq{t}\leq{T}}$. 
We are interested in the error of estimation at a fixed point $x_{0}\in\bbr_{0}$ and in minimax results of the form:
\begin{equation}\label{MinMaxIneq1}
	\liminf_{T \rightarrow\infty} 
	\left\{\inf_{\hat{s}_{_{T}}}\sup_{{s}\in\Theta}
	\bbe_{s} \left[\ell\left(T^{\gamma}
	\left(\hat{s}_{_{T}}(x_{0})-s(x_{0})\right)\right)\right]
	\right\}>0,
\end{equation}
where the infimum is over all the ``estimators'' 
$\hat{s}_{_{T}}$ based on the jumps of the L\'evy process 
$\left\{ X(t) \right\}_{0\leq{t}\leq{T}}$, 
$\Theta$ is a collection of L\'evy densities, 
and $\gamma>0$ is a constant depending on the family $\Theta$. 
In other words, (\ref{MinMaxIneq1}) implies the existence of 
a lower bound $B>0$ and a time $T_{0}$ such that from 
that time on, all
\emph{non-anticipative}\footnote{Here, non-anticipative means 
that the estimator is based on the jumps that 
occurred up to the present.} 
estimators $\hat{s}_{_{T}}$ will not do better than 
$T^{-\gamma}$ uniformly on $\Theta$, in the sense that 
there would exist an $s\in\Theta$ for which
	\[
		\bbe_{s} \left[l\left(T^{\gamma}
		\left(\hat{s}_{_{T}}(x_{0})-{s}(x_{0})\right)\right)\right]>B.
	\]
Therefore, the inequality (\ref{MinMaxIneq1}) impose a constraint  
on the rate of convergence at $x_{0}$ that the estimators can 
attain. By estimators, we mean a ``process''
$\hat{s}:\bbr_{0}\times\Omega\rightarrow \bbr$ such that
for each $x\in\bbr_{0}$, the random variable $\hat{s}(x;\cdot)$
is measurable with respect to the $\sigma$-field generated by the 
point process $\calJ$, while for each $\omega\in\Omega$, 
$\hat{s}(\cdot;\omega)$ is measurable with respect to the product 
$\sigma$-field.



The considered L\'evy densities satisfy
 a H\"older condition of order $\beta$ on a given window of 
estimation. 
Concretely, fix an interval $[a,b]\subset\bbr\backslash\{0\}$, and let $k\in\{0,1,\dots\}$ and
$\beta\in(0,1]$. Define the family $\Theta_{k+\beta}\left(L;[a,b]\right)$ of functions 
$f:\bbr\backslash\{0\}\rightarrow \bbr$ such that $f$ is $k$ times differentiable on $[a,b]$ and 
	\begin{equation}\label{HolderCond}
		| f^{(k)}(x_{1})-f^{(k)}(x_{2})|\leq L| x_{1}-x_{2} |^{\beta},\;\;\forall\;\;x_{1},x_{2}\in [a,b].
	\end{equation}
Below, $\calL$ stands for the class of all L\'evy densities; 
that is, all functions ${s}:\bbr_{0}\rightarrow\bbr_{+}$ such that
\[
	\int_{\bbr_{0}} \left( x^{2}\wedge{1}\right) s(x) dx <\infty.
\]
The following result is a minor variation of Theorem 6.5 of \cite{Kutoyants}. For completeness, we present its proof in 
Section \ref{SectionProofMinimax}.
\begin{thrm}\label{MinimaxResult1}
If $x_{0}$ is an interior point of the interval 
$[a,b]\subset\bbr\backslash\{0\}$, then
\begin{equation}\label{MinMaxIneq1b}
	\liminf_{T \rightarrow\infty} 
	\left\{\inf_{\hat{s}_{_{T}}}\sup_{{s}\in\Theta}
	\bbe_{p} \left[\ell\left(T^{\alpha/(2\alpha+1)}
	\left(\hat{s}_{_{T}}(x_{0})-{s}(x_{0})\right)\right)\right]
	\right\}>0,
\end{equation}
where $\alpha:=k+\beta$, 
$\Theta:=\calL\cap\Theta_{\alpha}\left(L;[a,b]\right)$, and 
 the infimum is over all the estimators $\hat{s}_{_{T}}$ 
based on those jumps of the 
L\'evy process $\left\{ X(t) \right\}_{0\leq{t}\leq{T}}$ 
whose sizes lie in $[a,b]$.
\end{thrm}
 As already noticed in \cite{Ibragimov}, the previous result can 
be strengthen to be in a certain sense uniform in $x_{0}\in(a,b)$
(see  Section \ref{SectionProofMinimax} for a proof)
\begin{cllry}\label{MinimaxResult2}
With the notation and hypothesis of Theorem \ref{MinimaxResult1},
\begin{equation}\label{MinMaxIneq2a}
	\liminf_{T \rightarrow\infty} 
	\left\{\inf_{\hat{s}_{_{T}}}\inf_{x\in(a,b)}
	\sup_{{s}\in\Theta}
	\bbe_{s} \left[\ell\left(T^{\alpha/(2\alpha+1)}
	\left(\hat{s}_{_{T}}(x)-{s}(x)\right)\right)\right]
	\right\}>0.
\end{equation}
\end{cllry}

Let us now apply the above assertion to obtain the long run 
minimax risk of \emph{measurable} estimators, 
under the $\bbl^{2}$-norm.
Here, measurable means that for each $\omega\in\Omega$,
$\hat{s}(\cdot;\omega)$ is a measurable function on 
$\left([a,b],\calB\left([a,b]\right)\right)$.
In Section \ref{SectionProofMinimax} a proof is given.
\begin{cllry}\label{MinimaxResult3}
Let $[a,b]$ be a closed interval of $\bbr\backslash\{0\}$, then
\begin{equation}\label{MinMaxIneq3a}
	\liminf_{T \rightarrow\infty} T^{2\alpha/(2\alpha+1)}
	\left\{\inf_{\hat{s}_{_{T}}}\sup_{{s}\in\Theta}
	\bbe_{s} \left[\int_{a}^{b}\left(\hat{s}_{_{T}}(x)-s(x)
	\right)^{2}dx\right]\right\}>0,
\end{equation}
where $\alpha:=k+\beta$, 
$\Theta:=\calL\cap\Theta_{\alpha}\left(L;[a,b]\right)$, and 
 the infimum is over all the measurable estimators 
$\hat{s}_{_{T}}$ based on the jumps of the 
L\'evy process $\left\{ X(t) \right\}_{0\leq{t}\leq{T}}$ 
whose sizes lie on $[a,b]$.
\end{cllry}
\begin{remk}\label{MinimaxMore}
The proofs of the previous results can be readily 
modified to cover even smaller classes 
of L\'evy densities $\Theta$. 
For instance, $\Theta = \calL \cap
\Theta_{\alpha}(L;[a,b]) \cap 
\{ s : \|s\|_{\bbl^{\infty}([a,b])}<R\}$. 
This class has a very close relationship with the family of 
Besov densities $\Theta(R,L)$ introduced in (\ref{LongRunRiskppe}). 
Indeed, the class $\Theta_{\alpha}(L;[a,b])$ 
is contained in 
$\mathcal{B}^{\alpha}_{\infty}\left(\bbl^{\infty}([a,b])\right)$
(see Section 2.9 of \cite{Devore}).
Since 
$\mathcal{B}^{\alpha}_{\infty}\left(\bbl^{\infty}\right)
\subset 
\mathcal{B}^{\alpha}_{\infty}\left(\bbl^{p}\right)$,
(\ref{MinMaxIneq3a}) holds true on $\Theta=\Theta(R,L)$.
Therefore, the p.p.e. on regular splines, described in the 
previous section, has the best possible 
 rate of convergence and moreover, 
achieves the minimax rate of convergence on 
$\Theta(R,L)$. This type of 
property is called \emph{adaptivity} in that, 
without knowing the smoothness of $s$
(controlled by $\alpha$), the p.p.e. reaches asymptotically 
the minimax risk up to a constant. 
See for instance Section 4 of \cite{Birge:1995} for a 
discussion on adaptivity.
 
\end{remk}


\section{Calibration based on discrete time data: approximation of Poisson integrals}
\label{SecAproxPoissnInt}\index{Poisson integrals!approximation of}
One drawback to the method outlined in Section \ref{BscMethod} is that in general we do not observe the jumps of a L\'evy process $X=\left\{X(t)\right\}_{t\geq 0}$. In practice, we can aspire to
sample the process $X(t)$ at discrete times, but we are neither able
to measure the size of the jumps $\Delta X(t) \equiv X(t)-X(t^{-})$ nor the 
times of jumps $\left\{t: \Delta X(t)>0\right\}$. Poisson integrals\index{Poisson integrals} of the
type
	\begin{equation}
	\label{GenericIntegral}
		I\left(f\right)\equiv \iint\limits_{[0,T] \times \bbr_{0}} f(x) {\cal J}(dt,dx)=
		\sum\limits_{t\leq T} f(\Delta X(t)),
	\end{equation}
are simply not accessible. In this section, we discuss the approximation
of the integral (\ref{GenericIntegral}) based on time series of the
form  $\left\{X(t_{k}^{n})\right\}_{k=0}^{n}$, where $t_{k}^{n} = \frac{kT}{n}$. 

Let us motivate our approximation scheme. The natural way of interpolating the sample 
path of a L\'evy process from the sampling observations $\left\{X(t_{k}^{n})\right\}_{k=0}^{n}$
is to take a c\`adl\`ag piecewise constant approximation of the form
	\begin{equation}
	\label{DiscreteSkeleton}
		X^{n}(t) \equiv \sum\limits_{k=1}^{n} X\left(t_{k-1}^{n}\right) \mathbf{1}\left(t\in [t^{n}_{k-1},t^{n}_{k})\right),\;\;\;t\in [0,T),
	\end{equation}
where as usual $\mathbf{1}$ is the indicator function of the corresponding set. It
is quite simple to prove that $X^{n}$ converges to
$X$ at finitely many points with probability one (a quality shared by any
right-continuous process $X$).  Furthermore,  the approximated 
process $X^{n}$, having independent increments, converges to ${\bf X}$ in $D[0,\infty)$, under the Skorohod metric (see VI of \cite{Pollard} and for concrete Example VI.18). 
Hence, we might expect that
	\begin{equation}
	\label{GenericApprox}
		I_{n}\left(f\right) \equiv \sum\limits_{t\leq T} f(\Delta X^{n}(t))
		=\sum\limits_{k=1}^{n} f\left(X\left(t_{k}^{n}\right)-X\left(t_{k-1}^{n}\right)\right),
	\end{equation}
converges to (\ref{GenericIntegral}) as $ n \rightarrow \infty$. Indeed, we prove 
the weak convergence of  (\ref{GenericApprox}) to (\ref{GenericIntegral}) using 
well-know facts on the transition distributions of $X$ in small time (see for instance pp. 39 of \cite{Bertoin}, Corollary 8.9 of \cite{Sato}, or Corollary 3 of \cite{Ruschendorf}).
More precisely, 
	\noindent \begin{lmma}
	 Let   $X=\left\{ X(t) \right\}_{t \geq 0}$ be a L\'evy process with L\'evy measure $\nu$.
	Then:\\
		\textbf{\emph 1)}  For each $a>0$,
		\begin{equation}\label{ConvDistr}
			\lim_{t\rightarrow 0} \frac{1}{t} \bbp \left( X(t) >
			a\right)=\nu([a,\infty)),\;\;{\rm and}\;\;
			\lim_{t\rightarrow 0} \frac{1}{t} \bbp \left( X(t) \leq
			-a\right)=\nu((-\infty,-a]).
		\end{equation}
		\textbf{\emph 2)}\index{L\'evy processes!small time distribution!moments}
		 For any continuous bounded function $h$ vanishing on a neighborhood of the origin,
		\begin{equation}
		\label{WeakLimit1}
			\lim_{t\rightarrow 0} \frac{1}{t} \bbe \left[ h\left(X(t)\right) \right]
			= \int_{\bbr_{0}} h(x) \nu(dx).
		\end{equation}
	\end{lmma}
	\begin{remk} \label{EasyRmk}
		In particular, the two parts in the previous Lemma imply (\ref{WeakLimit1}) when  
		$h(x) =  \mathbf{1}_{(a,b]}(x) f(x)$,  where $[a,b]$ is an interval of $ \bbr_{0}$ 
		and $f$ is a continuous function.
	\end{remk}
It is worth mentioning that \cite{Ruschendorf} provides stronger results for the small-time distributional properties of $X(t)$. The following theorem summarizes some of their results. 
	\noindent \begin{thrm}
	\label{TranstionDistrResults}\index{L\'evy processes!small time distribution}
	\index{Spectral function}\index{L\'evy processes!small time density}
	 Let  $X=\left\{ X(t) \right\}_{t \geq 0}$ be a L\'evy process with L\'evy measure $\nu$.
	 Let $F_{t}$ be the distribution function of $X(t)$ and $G$ the \emph{spectral function} of $\nu$; 
	 i.e. $G(x)=\nu([x,\infty))$ for $x>0$ and $G(x)=\nu((-\infty,x])$ for $x<0$. 
	 The following properties  hold:\\
		\emph{(i)}\index{L\'evy processes!small time density}
		 	If $F_{t}$ and $G$ have densities\footnote{
			The function $g \geq 0$ is said to be the density of the spectral function $G$
			 if $G'(x)=g(x)$ for $x<0$ and $G'(x)=-g(x)$ for $x>0$.}
			 $f_{t}$ and $g$, then for $x\neq 0$
			 \begin{equation}\label{Dnsties}
			 	\lim_{t\rightarrow 0} \frac{1}{t} f_{t}(x) =
				\frac{\partial}{\partial t} f_{t}(x) \Big|_{t=0} = g(x),
			\end{equation}
			where we additionally assume that $F_{t}(x)$ is continuous in a neighborhood of 
			$(t=0,x)$ and that moreover 
			 $(\partial/\partial t) F_{t}(x)$, $(\partial/\partial x) F_{t}(x)$, 
			and 	$(\partial/\partial t)(\partial/\partial x) F_{t}(x)$ 
			exist and are continuous in $(t=0,x)$.\\
If $h$ is continuous and bounded and if $\lim_{|x|\rightarrow 0} h(x)| x |^{-2} =0$,
then
\[
\lim_{t\rightarrow 0} \frac{1}{t}\bbe\left[h\left(X(t)\right)\right] = 
\int_{\bbr_{0}} h(x) \nu(dx).
\]
Moreover, if  $\int_{\bbr_{0}} (|x| \wedge 1) \nu(dx)<\infty$, it is 
enough to postulate that $h(x) (| x | \wedge 1)^{-1}$ is continuous and bounded.
\end{thrm}
	
\noindent Limiting results like (\ref{WeakLimit1}) are useful to establish the convergence in distribution  of ${\rm I}_{n}\left(f\right)$ since
	\[
		\bbe\left[ e^{{\rm i}u I_{n}(f) }\right]=
		\left(\bbe \left[e^{{\rm i}uf\left(X\left(\frac{T}{n}\right)\right)}\right]
		\right)^{n} 
		 = \left(1+\frac{a_{n}}{n} \right)^{n},
	\]
where $a_{n}=n\bbe\left[h\left(X\left(\frac{T}{n}\right)\right)\right]$ with
$h(x)=e^{{\rm i} u f(x)}-1$. So, if $f$ is such that 
\begin{equation}\label{WeakLimit1b}
	\lim_{t\rightarrow 0} \frac{1}{t} 
	\bbe \left[ e^{{\rm i} u f(X(t))}-1 \right]
			= \int_{\bbr_{0}} \left(e^{{\rm i} u f(x)}-1\right) \nu(dx),
\end{equation}
then $a_{n}$ converges to $a\equiv T\int_{\bbr_{0}} h(x) \nu(dx)$, and thus
\[
		\lim_{n \rightarrow \infty} \left(1+\frac{a_{n}}{n} \right)^{n}=
		\lim_{n \rightarrow \infty} e^{n \log\left(1+\frac{a_{n}}{n} \right)}=
		e^{a}.
	\]
We thus have the following result (see Section \ref{proofsChap3} for verification):
\noindent \begin{prop}\label{ConvergenceApprox} 
\index{Poisson integrals!approximation of} 
Let $X=\left\{ X(t) \right\}_{t \geq 0}$ be a L\'evy process with 
L\'evy measure $\nu$. Then, 
\[
	\lim_{n\rightarrow \infty} \bbe\left[ e^{{\rm i}u I_{n}(f) }\right]
	=\exp\left\{ T \int_{\bbr_{0}} \left(e^{{\rm i} u f(x)}-1\right)\nu(dx)\right\},
\] 
if $f$ satisfies either one of the following:
\item{1)} $f(x)= \mathbf{1}_{(a,b]}(x) h(x)$ for an interval 
$[a,b] \subset \bbr_{0}$ and a continuous function $h$;
\item{2)} $f(x)$ is continuous on 
$\bbr_{0}$ and $\lim_{|x|\rightarrow 0} f(x)| x |^{-2} =0$.\\
In particular, $I_{n}(f)$ converges in distribution to $I(f)$ under any of the two 
previous conditions. 
\end{prop}

	\begin{remk}
	\label{AsympMeanVarAppxInt}\index{Poisson integrals!approximation of}
		Clearly, if $f$ and $f^{2}$ satisfy (\ref{WeakLimit1}), then the mean and variance of 
		$I_{n}(f)$ obey the asymptotics:
		\[
			\lim_{n\rightarrow \infty} \bbe \left[ I_{n}(f) \right]=
			T \int_{\bbr_{0}} f(x) \nu(dx);
		\]
		\[
			\lim_{n\rightarrow \infty} {\rm Var} \left[ I_{n}(f) \right]=
			T \int_{\bbr_{0}} f^{2}(x) \nu(dx).
		\]
	\end{remk}


\section{Estimation Method}
\label{EstimationmMethod}
Let us summarize the previous sections and outline the proposed
algorithm of estimation:
	\begin{description}
	\item[Statistician's parameters:] The procedure is fed with a Borel
	\emph{window of estimation} $D\subset\bbr_{0}$, a collection 
	$\left\{{\cal S}_{m}\right\}_{m \in {\cal M}}$ of finite dimensional \emph{linear models}
	 of $L^{2}\left((D,\eta)\right)$, and a \emph{level of penalization} $c>1$.

	\item[Model and data:] It is assumed that a L\'evy process
	$\left\{X\left(t\right)\right\}_{t \in [0,T]}$ is monitored at 
	equally spaced times $t_{k}^{n}= k\frac{T}{n}$, $k=1,\dots,n$, during the time period $[0,T]$.
	The data consists of the time series $\left\{X\left(t^{n}_{k}\right)\right\}_{k=1}^{n}$.
	The L\'evy process admits a regularized L\'evy density $s$ 
under the measure $\eta$ on $D$ (see Definition \ref{Regularization}).

	\item[Estimators:] Inside the linear model ${\cal S}_{m}$, the estimator of $s$ is
	the \emph{approximated projection estimator}: 
	\begin{equation}\label{ApprxProjEstm}
	\index{projection estimator!approximation of}
	\index{approximated projection estimator}
		\hat s^{n}_{m}(x) \equiv \sum_{i=1}^{d_{m}} 
		\hat\beta^{n}_{i,m} \varphi_{i,m} (x),
	\end{equation}
	where $\left\{ \varphi_{1,m},\dots,\varphi_{d_{m},m}\right\}$ is an orthonormal 
	basis for ${\cal S}_{m}$, and 
	\begin{equation}\label{DefApprBetas}
		\hat\beta^{n}_{i,m} \equiv \frac{1}{T} 
		\sum\limits_{k=1}^{n}
		\varphi_{i,m}\left(X\left(t_{k}^{n}\right)-X\left(t_{k-1}^{n}\right)\right),
	\end{equation}
	is the estimator of the inner product 
	$\beta_{i,m}\equiv\int_{D} \varphi_{i,m}(x) s(x) \eta(dx)$, for $i=1,\dots,d_{m}$.
	Across the collection of linear models 
	$\left\{{\cal S}_{m}: m \in {\cal M}\right\}$, the estimator $\hat s^{n}_{m}$ which minimizes 
	\(
		-\| \hat s ^{n}_{m} \|^{2}  + c\; {\rm pen}^{n}(m),
	\)
	is selected, where
	\[
		{\rm pen}^{n}(m)=\frac{1}{T^{2}} \sum\limits_{k=1}^{n} 
		\;\left(\sum_{i=1}^{d_{m}} 
		\varphi^{2}_{i,m}\left(X\left(t_{k}^{n}\right)-X\left(t_{k-1}^{n}\right)\right) \right).
	\]
	\end{description}

\begin{remk} \label{ComprsnDenstyEst}
It is worthwhile to point out the great similarity of the 
scheme above to some methods of density estimation introduced 
in  \cite{Birge}. 
In this paper, the authors estimate the probability density 
function $f$ of a random sample $X_{1},\cdots,X_{n}$ 
by projection estimators of the form: 
	\begin{equation}\index{estimation of density functions}
		\hat f (x) = \sum_{i=1}^{d} \left\{\frac{1}{n} \sum_{k=1}^{n} 
		\varphi_{i}\left(X_{k}\right)\right\}\varphi_{i} (x),
	\end{equation}
where $\left\{\varphi_{i}\right\}_{i=1}^{d}$ is an orthonormal 
basis of a linear space ${\cal S}$ of $L^{2}(\bbr,dx)$. 
More generally, $f$ can be the density function with respect 
to a measure $\mu$ in the sense that 
$\bbp\left[X_{i}\in{\cdot}\right]=\int_{\cdot} f(x)\mu(dx)$, 
and the projection estimator above will be well defined 
provided that $f\in L^{2}(\bbr,\mu(dx))$. 
To solve the problem of model selection, they introduced 
penalized projection estimators. 
One considered penalty function there is
	\[
		{\rm pen}({\cal S}) =\frac{2}{n(n+1)} 
		\sum\limits_{k=1}^{n} 
		\;\sum_{i=1}^{d} \varphi^{2}_{i}\left(X_{k}\right).
	\]
In some sense, the method outlined at the beginning of this section 
``works'' as a byproduct of the small time qualities of L\'evy processes 
and of standard methods of nonparametric estimation for probability 
densities. Indeed, consider the statistics
	\[
		\hat\beta^{n,j}_{i,m} \equiv \frac{n}{Tj} 
		\sum\limits_{k=1}^{j}
		\varphi_{i,m}\left(X\left(t_{k}^{n}\right)-
		X\left(t_{k-1}^{n}\right)\right),
	\]
where $T/n$ is the time span of the increments 
and $j$ is the number of increments in the sample.
From \cite{Birge:1994}, 
as $j$ progresses,
	\begin{equation}\label{ProjEstpdf}
		\hat s^{n,j}_{m}(x) \equiv 
		\sum_{i=1}^{d_{m}} \hat\beta^{n,j}_{i,m} 
		\;\varphi_{i,m} (x),
	\end{equation}	 
estimates the orthogonal projection of $\frac{n}{T} f_{T/n}(x)$ 
on $\calS_{m}$, where $f_{t}$ stands for the probability density 
function of $X(t)$ (if it exists).  
On the other hand, \cite{Ruschendorf} proves that 
$\frac{n}{T} f_{T/n}(x)$ converges to the L\'evy density $p$, 
as $n\rightarrow\infty$ (under some regularity conditions). 
Therefore, for large $n$ and $j$, (\ref{ProjEstpdf}) will 
approximate the projection of $p$ on $\calS_{m}$. 
Notice that in general, a.s.
	 \[
	 	\lim_{n \rightarrow \infty} \lim_{j \rightarrow \infty} 
		\frac{n}{Tj} \sum\limits_{k=1}^{j}
		\varphi\left(X\left(t_{k}^{n}\right)-
		X\left(t_{k-1}^{n}\right)\right) = 
		\int_{\bbr_{0}} \varphi(x) \nu(dx),
	 \]
whenever $\varphi$ is such that  the limit (\ref{WeakLimit1}) holds. 
Our penalized projection estimators (\ref{ApprxProjEstm}) 
are obtained from (\ref{ProjEstpdf}) by taking $n=j$. 
It is not clear from the references just mentioned whether taking 
$n=j\rightarrow\infty$ will produce good results or not. 
We will see below that this is the case.
\end{remk}

Let ${\cal R}(X)$ be the linear space of measurable functions 
$h$ such that (\ref{WeakLimit1}) is satisfied. For instance, 
${\cal R}(X)$ contains the functions $f$ satisfying 
conditions (1) or (2) in Proposition \ref{ConvergenceApprox}.
The following result holds true 
(see Section \ref{proofsChap3} for a proof).
\begin{prop} \label{ApprProjAsymt}
Let $s^{\bot}_{m}$ be the orthogonal projection of $s$ on 
${\cal S}_{m}$. If $\varphi_{i,m}$ and $\varphi^{2}_{i,m}$ 
belong to ${\cal R}(X)$ for every $m \in {\cal M}$ 
and $i=1,\dots,d_{m}$, then the approximated projection 
estimator $\hat{s}^{n}_{m}$ of $s$ on $\mathcal{S}_{m}$ 
(based on $n$ equally spaced observations) satisfies:
	\begin{equation}\label{AsymRisk}
		\lim_{n\rightarrow \infty} 
		\bbe\left[\|\hat{s}^{n}_{m} - s^{\bot}_{m}\|^{2}\right]= 
		\bbe\left[\|\hat{s}_{m} - s^{\bot}_{m}\|^{2}\right].
	 \end{equation}
Moreover, 
	\[
		\lim_{n\rightarrow \infty} 
		\bbe\left[\|\hat{s}^{n}_{m} - s\|^{2}\right]= 
		\bbe\left[\|\hat{s}_{m} - s\|^{2}\right].
	\]
\end{prop}


\section{Numerical tests of projection estimators}
In this section, we try to assess the performance 
of some penalized projection estimators based on 
simulated L\'evy processes. 
Piecewise constant functions are considered, and 
for their intrinsic relevance in mathematical finance, 
two classes of L\'evy processes are 
studied: Gamma  and variance Gamma processes. 
A method of least-squares errors is also applied to generate 
parametric L\'evy densities that closely fit the nonparametric 
outputs.
\subsection{Specifications of the statistical methods}\label{specifModel}
Let us describe in greater details the considered projection 
estimators. 
To simplify notation, ${\cal J}\left( B \right)$ is used 
instead of ${\cal J}\left( [0,T] \times B \right)$ of 
(\ref{JumpMeasure}) when referring to the number of jumps 
of sizes in $B \in  {\cal B}(\bbr_{0})$ occurring prior to $T$.
Let  ${\cal C}:a = x_{0} < x_{1} < \dots < x_{m} = b $ be a partition of
the interval $D \equiv [a,b]$ ($0<a$ or $b<0$), and let ${\cal S}_{{\cal C}}$ be the span of the indicator 
functions $\chi_{[x_{0},x_{1})} , \dots, \chi_{[x_{m-1},x_{m})}$. 
In other words, the linear model ${\cal S}_{{\cal C}}$ consists of ``histogram functions''
\index{histogram projection estimators}\index{projection estimator!on histogram type functions}
 on the window $D$ with cutoff points in ${\cal C}$. 
We assume that the L\'evy process has a L\'evy density $s$ bounded outside of any
neighborhood of the origin. This assumption is very mild, and yet  good enough for the 
integral $\int_{D} s^{2}(x) dx$ to be finite. In that case, the orthogonal projection of $s$ onto 
${\cal S}_{{\cal C}}$ exists (under the standard inner product  of $L^{2}\left(D,dx\right)$), 
and thus the projection estimation on ${\cal S}_{{\cal C}}$ is meaningful. In the terminology 
 of Section \ref{BscMethod}, the regularizingn measure is simply $dx$, the regularized 
L\'evy density coincides with the L\'evy density, and the orthonormal basis 
$\{ \varphi_{1},\dots,\varphi_{m}\}$ for ${\cal S}_{{\cal C}}$ is 
	\[
		\varphi_{i}(x)= \frac{1}{\sqrt{x_{i}-x_{i-1}}} \chi_{[x_{i-1},x_{i})}(x),\;\; i=1,\dots,m.
	\]	
According to the \emph{basic estimation method} outlined in Section \ref{BscMethod},
the projection estimator on the linear model ${\cal S}_{{\cal C}}$ is given by 
	\begin{equation}
	\label{PrjEstHist}
		\hat s_{{\cal C}}(x) = \frac{1}{T}
		\sum_{i=1}^{m} \; \frac{{\cal J}([x_{i-1},x_{i}))}{x_{i}-x_{i-1}}
		\; \chi_{[x_{i-1},x_{i})}(x).
	\end{equation}
Following the heuristics of Section \ref{BscMethod} and Theorem \ref{PreOracleIneq} part (b), an appealing procedure to select a projection estimator of the form (\ref{PrjEstHist}) is to 
look for the minimization of the following penalized contrast value 
	\begin{equation}
	\label{PenConstrValue} 
		\frac{1}{T^{2}} \sum_{i=1}^{m} \frac{1}{x_{i}-x_{i-1}} 
		 \left\{ c \; {\cal J}([x_{i-1},x_{i})) 
		-\left[{\cal J}([x_{i-1},x_{i}))\right]^{2}\right\}.
	\end{equation}
Here, $c>1$ is a constant that controls the level of penalization. 
In fact, Theorem \ref{PreOracleIneq} and Corollary \ref{OracleIneq} 
ensure that, 
 for large enough $T$, the previous procedure will yield competitive results against the best 
projection estimator. For that to happen it is necessary to restrict ourselves to models ${\cal C}$ 
satisfying $D_{{\cal C}} \leq T$, where $D_{{\cal C}}$ is defined as in (\ref{DefOfDm}). In this case, 
the constant $D_{{\cal C}}$ is $1/\min_{1\leq i \leq m} \{x_{i} - x_{i-1}\}$ as seen from 
Remark \ref{DmInOB}. 

\index{histogram projection estimators!regular partitions}
\index{projection estimator!on histogram type functions}
The simplest case is to take regular partitions $\left\{ x_{i} = a+i\Delta x
\right\}_{i=0}^{m}$, where $\Delta x = (b-a)/m$ is the mesh of the
partition. Then, the projection estimators of (\ref{PrjEstHist}) becomes 
	\begin{equation}
	\label{PrjEstHistRegPart}
		\hat s_{m}(x) \equiv \frac{m}{T(b-a)}
		\sum_{i=1}^{m} \; {\cal J}([x_{i-1},x_{i}))
		\; \chi_{[x_{i-1},x_{i})}(x),
	\end{equation}
and penalized projection estimation will look to minimize
	\begin{equation}
	\label{PenProjHist}
		 \frac{m}{T^{2} (b-a)} \left( c {\cal J}([a,b)) 
		- \sum_{i=1}^{m}\left({\cal J}([x_{i-1},x_{i}))\right)^{2}\right),
	\end{equation}
over all $m$ such that $D_{m} = m/(b-a)$ is smaller than $T$. 

\index{projection estimator!application to fit parametric models}
For comparisons against  other procedures and to assess the 
goodness of fit to 
specific parametric models, it is useful to determine the 
parametric model of a given type 
that ``best fits'' our non-parametric estimators; 
for instance, suppose the we want 
to assess whether or not the nonparametric results 
supports the parametric 
Gamma model for the L\'evy density.
 The method of least square errors provides 
an easy solution to this problem.\index{least-squares methods}
\index{projection estimator!with least-squares methods}
For instance, if $s_{\theta}(x)$ is the
parametric form of the L\'evy density, a plausible estimator 
of $\theta$ is
	\[
		\hat\theta = {\rm argmin}_{\theta} \; 
		d(s_{\theta}, \hat s),
	\]
where $\hat s$ is the (penalized) projection estimator on a given
family of linear models, and $d$ is a function that accounts for
the difference between $s_{\theta}$ and $\hat s$. For instance, 
for a fixed set of points 
$\left\{ x_{i} \right\}^{k}_{i=1}\subset D$, $d(\cdot,\cdot)$ can 
simply be defined for functions $f$ and $g$ as
	\[
		d(f,g) \equiv 
		\sum_{i=1}^{k} \left[f(x_{i})-g(x_{i})\right]^{2}.
	\]
It is sometimes preferable to use a least-square method that is 
linear in the parameters, and hence, is robust against 
numerical errors.  In that case,  we can look for a functional $L$ so
that $L(s_{\theta})$ is linear in $\theta$ and define
	\[
		d(f,g) \equiv \sum_{i=1}^{k} \left[L(f)(x_{i})-L(g)(x_{i})\right]^{2}.
	\]
As an example, consider the L\'evy density of a Gamma L\'evy 
process with parameters $\alpha$ and $\beta$:
	\[
		s(x) = \frac{\alpha}{x} e^{-x/\beta},\;\; x>0.
	\] 
Given a projection estimator $\hat{s}$ of $s$, 
least-square estimates of $\alpha$ and $\beta$ can be
constructed from
\begin{equation}
\label{RegresGamma1}
{\rm argmin}_{\alpha,\beta} \sum_{i=1}^{m} \left( \frac{\alpha}{x_{i}}
\exp\left(-\frac{x_{i}}{\beta}\right) - \hat s(x_{i})\right)^{2},
\end{equation}
where $\left\{ x_{i} \right\}^{k}_{i=1}\subset D$.
Notice that the estimation would be very susceptible 
to the points close to the origin.
Instead, a regression method that is linear in the
parameters can be devised using a logarithmic transformation
as follows
\begin{equation}
\label{RegresGamma2}
	{\rm argmin}_{\alpha,\beta} \sum_{i=1}^{m} 
	\left( -\frac{1}{\beta}x_{i}+\log(\alpha)
	- \log(x_{i}\hat s(x_{i}))\right)^{2}.
\end{equation}

\subsection{Estimation of Gamma L\'evy densities.}\label{GammaLevyProc}
\subsubsection{The model} 
L\'evy Gamma processes are fundamental building blocks in the 
construction of other L\'evy processes like  
the variance Gamma model \cite{Madan1} and the generalized Gamma 
convolutions \cite{Bondesson}. 
Moreover, by Bernstein's theorem, any L\'evy density of the form 
$u(x)/\vert x \vert$, where $u$ is a completely monotone function, 
is the limit of superpositions of Gamma L\'evy densities. 

The Gamma L\'evy process 
$X=\left\{ X(t)\right\}_{t \geq 0}$ is determined by two
positive parameters $\alpha$ and $\beta$ so that the probability 
density function of $X(t)$ is
\begin{equation}\label{DnstyFunctGamma}
	f_{t}(x) = \frac{x^{\alpha t -1} e^{-x/\beta}}{\Gamma(\alpha t)\beta^{\alpha t}},
\end{equation}
for $x>0$. The characteristic function of $X(t)$ is  
\[
E\left[e^{iuX(t)}\right]= \left(1 - i \beta t \right)^{\alpha t} =
\exp\left[t\left(\alpha \int_{0}^{\infty} \left(e^{iux}-1\right)\nu(dx)\right) \right],
\]
where the L\'evy measure $\nu$ is
\begin{equation}
\label{LevyDstyGamma}
\nu(dx) =\frac{\alpha}{x} \exp\left(-\frac{x}{\beta}\right) dx,\;\;{\rm for}\;\;x>0;
\end{equation}
see \cite{Kolmogorov} pp. 87 or Example 8.10. of \cite{Sato}. From the point of view of the marginal densities, 
$\beta$ is a \emph{scale parameter} and $\alpha$ is a \emph{shape parameter}. 
In terms of the jump activity, $\alpha$ controls
the overall activity of the jumps, while $\beta$ takes charge of the heaviness of
the L\'evy density tail, and hence, of the frequency of big jumps. Notice that changes
in the time units is statistically equivalent to changes of the parameter $\alpha$, 
while changes in the units  at which the values of $X$ are measured 
are statistically reflected on changes of the parameter $\beta$. That is to say, 
the scaled process $\left\{ c X(ht)\right\}_{t \geq 0}$ is also a  
Gamma L\'evy process with shape parameter 
$\alpha h$ and scale parameter $\beta c$. This property is consistent with the previous remark on
$\alpha$ taking charge of the jump activity and on 
$\beta$ taking charge of the frequency of large jumps. 

\subsubsection{The simulation procedure}\label{SectSimProcGammaProc}
Simulation schemes based on series representation are used 
to generate Gamma L\'evy processes. 
Such procedures allow to retrieve a sample of the jumps of 
the process. Concretely, following \cite{Rosinski:2001}, the process  
\begin{equation}
	X(t)\equiv \beta \sum_{i=1}^{\infty} V_{i}
	\exp\left(-\frac{\Gamma_{i}}{\alpha}\right)
	{\bf 1}[U_{i} \leq t],
\end{equation}
is a Gamma L\'evy process on $[0,T]$ with shape parameter 
$\alpha$ and scale parameter $\beta$ provided that 
$\{\Gamma_{i} \}_{i\geq{1}}$ is a homogenous Poisson process 
with intensity  $1$, $\{V_{i} \}_{i\geq{1}}$ are independent 
exponential r.v. with mean 1, $\{U_{i} \}_{i\geq{1}}$ are 
i.i.d.  uniformly distributed on $[0,T]$, and these three series 
are mutually independent. 
Below, we shall truncate the series to $n$ terms in order to 
generate a sample path, and in particular,  to approximate the  
jump process ${\cal J}$ of $X$ by 
 \begin{equation}\label{ApproxPoissProc}
	{\cal J}_{n}(\cdot) \equiv 
	\sum_{i=1}^{n} \delta_{(U_{i},J_{i})} (\cdot),
\end{equation}
where 
$J_{i}\equiv\beta V_{i}\exp\left(-\frac{\Gamma_{i}}{\alpha}\right)$.

\subsubsection{The numerical results}\label{NumResults}	
\begin{figure}[htp]
	{\par\centering
	\includegraphics[width=8cm,height=8cm]{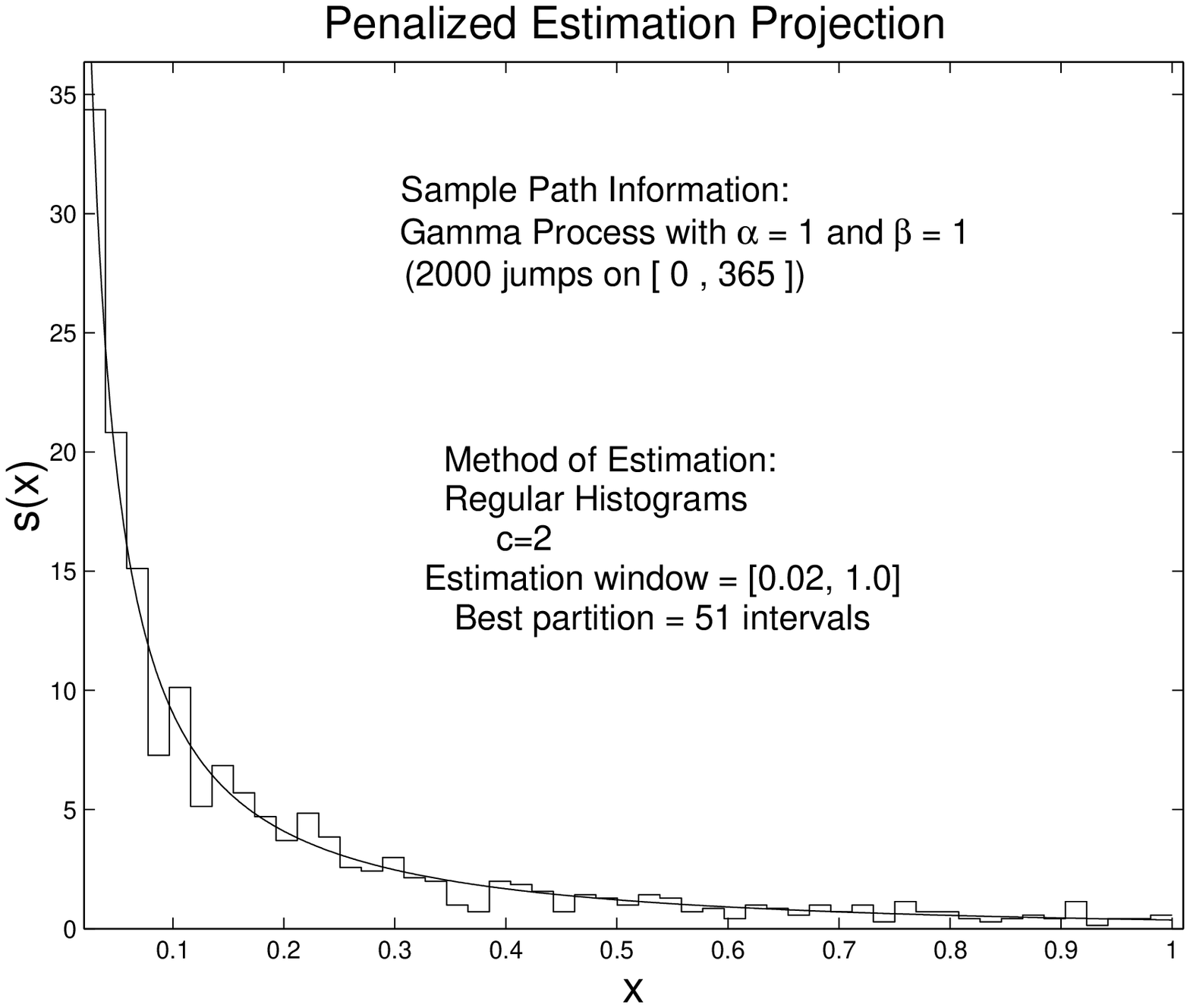} 
	\par}
	\caption{\label{EstGamma1} 
	Penalized projection estimation of $\frac{e^{-x}}{x}$.}

{\par \centering 
	\includegraphics[width=8cm,height=8cm]{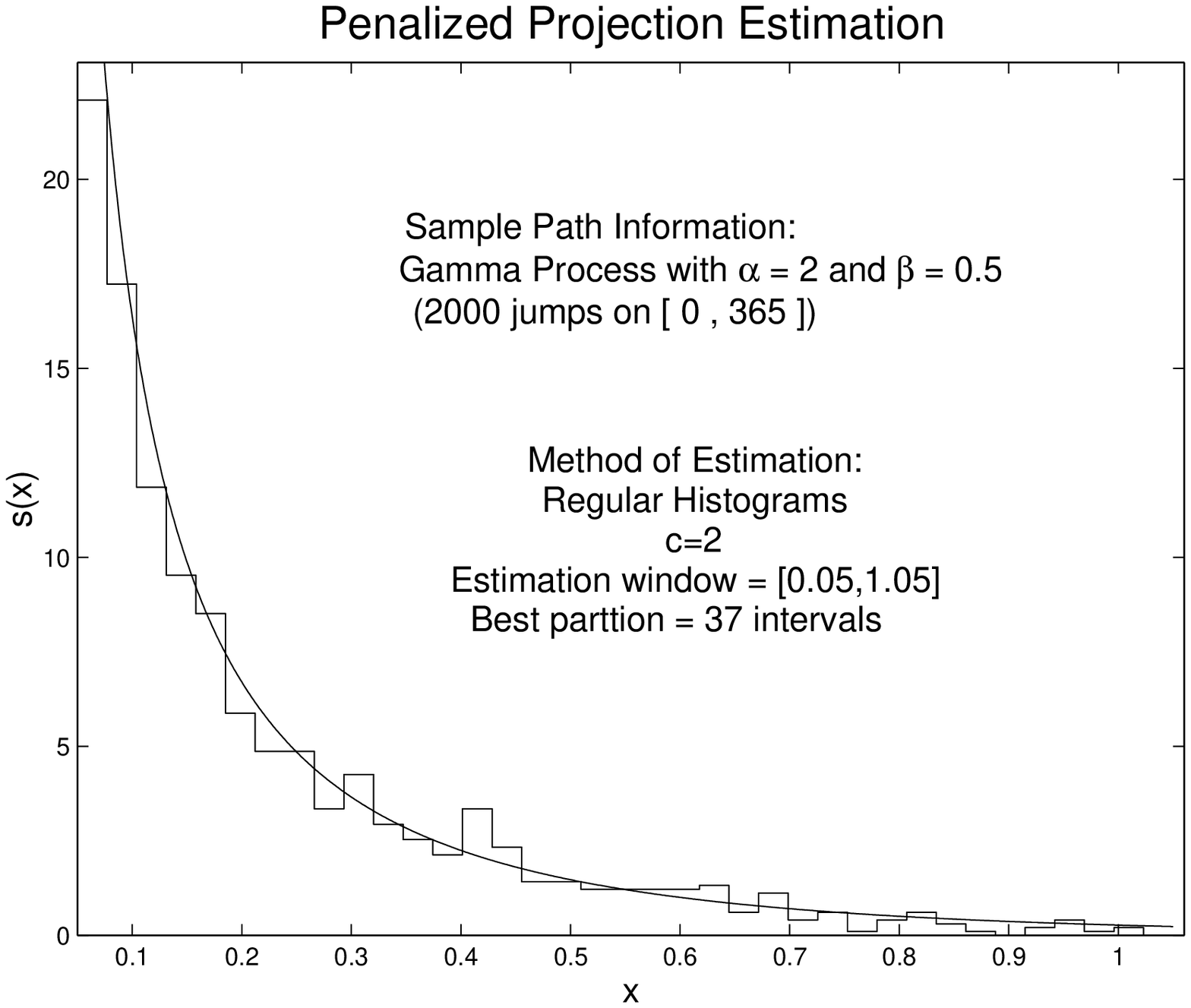} \par}
	\caption{\label{EstGamma2a} Penalized projection estimation of
	$\frac{2}{x}\exp\left(-2x\right)$.}
\end{figure}

We now present a few examples to illustrate the technique of projection estimation on histogram 
functions based on regular partitions (see Section \ref{specifModel} for the specifications of the 
estimation method). Figure \ref{EstGamma1} shows the Gamma L\'evy
density  with $\alpha = 1$ and $\beta =1$, and the penalized projection histogram of the form 
(\ref{PrjEstHistRegPart}). The estimation is based on 2000 jumps of the Gamma 
L\'evy process on $[0,365]$, and their resulting Poisson integrals obtained by using 
(\ref{ApproxPoissProc}) instead of $\calJ$.   
The least-square method (\ref{RegresGamma2}), taking the $x_{i}$'s as the mid
points of the partition intervals, yields the estimators $\hat 
\alpha = 0.932$ and $\hat \beta = 1.055$. The maximum likelihood estimators
based on the increments of the sample path of time length 1 are 1.015 for
$\alpha$ and 0.949 for $\beta$ (we do not observe real improvement if the time
length of the increments is reduced).    

\begin{figure}[htp]
	{\par \centering 
	\includegraphics[width=8cm,height=8cm]{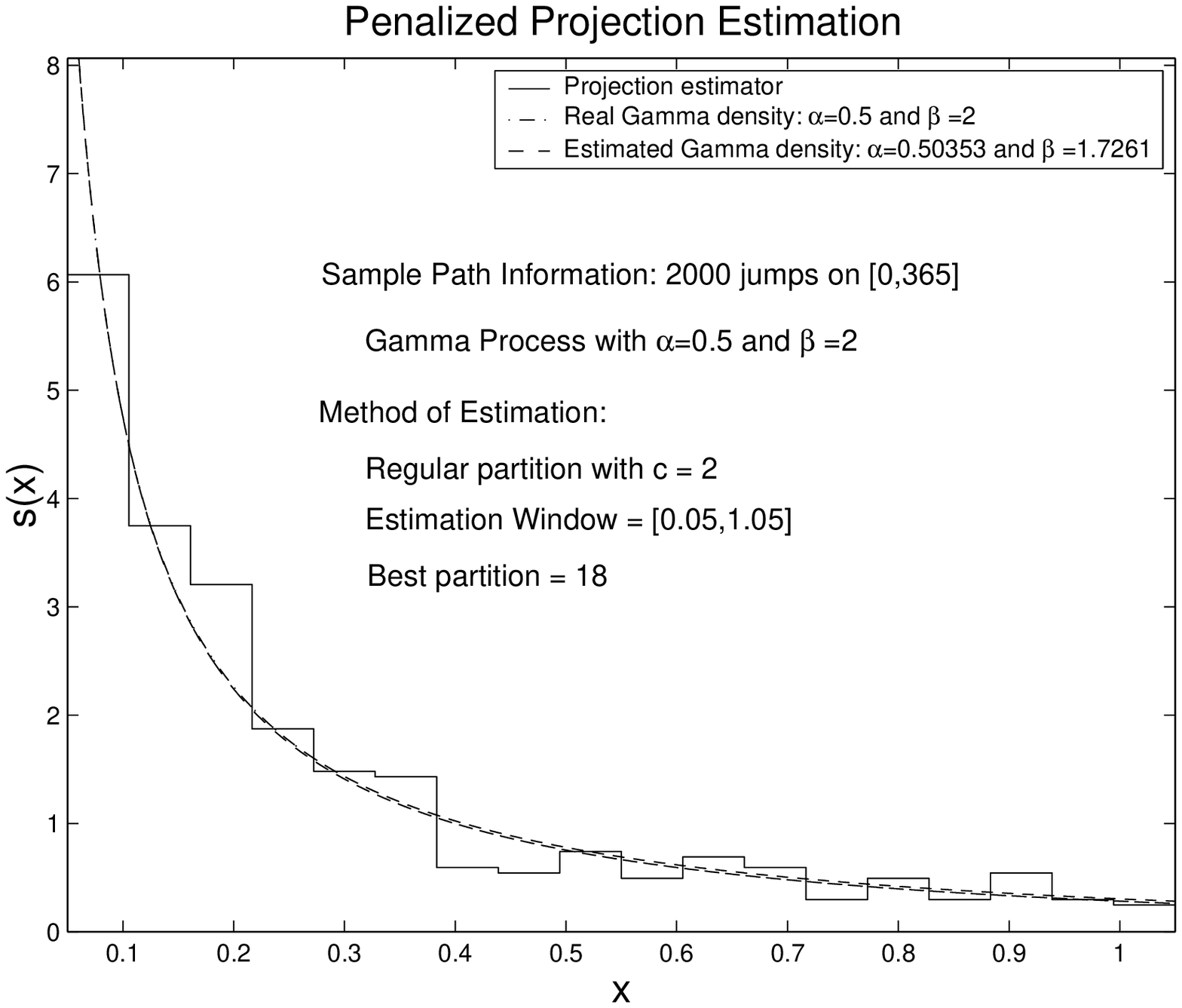} \par}
	\caption{\label{EstGamma2b} Penalized projection estimation of
	$\frac{1}{2x}\exp\left(-\frac{x}{2}\right)$.}
\end{figure}

In the next simulation, we consider a Gamma density with a lighter tail
($\beta = 0.5$) and more jump activity ($\alpha=2$). The opposite
setting was also studied: a heavier tail determined by a $\beta = 2$ and
a lower jump activity given by an $\alpha = 0.5$ (see Figures \ref{EstGamma2a} and   \ref{EstGamma2b}).
In the first scenario, the least-square method estimators are $\hat \alpha =
1.907$ and $\hat \beta = 0.472$, while the maximum likelihood estimators are
$1.924$ and $0.527$, respectively. For this second Gamma density, the
least-square method (\ref{RegresGamma1}), taking the $x_{i}$'s as the
midpoints of the partition intervals, produce estimators $\hat \alpha = 0.5$
and $\hat \beta= 1.72$, while the maximum likelihood estimators are $0.55$
and $1.99$, respectively.

\begin{figure}[htp]
	{\par \centering
	\includegraphics[width=8cm,height=8cm]{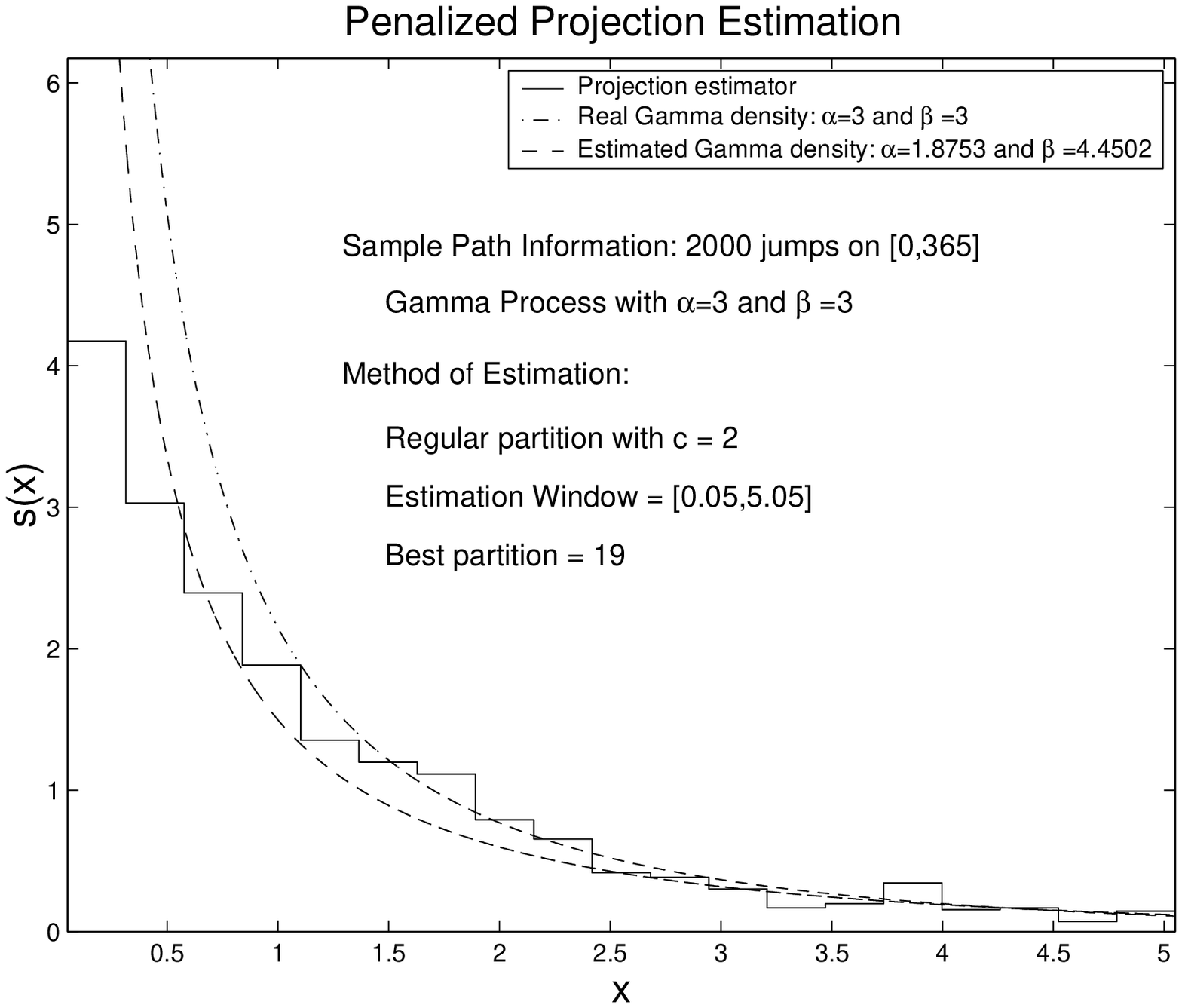} \par}
	\caption{\label{EstGamma3a} Penalized projection estimation of
	$\frac{3}{x}\exp\left(-\frac{x}{3}\right)$.}
\end{figure}

\begin{figure}[htp]
	{\par\centering
	\includegraphics[width=8cm,height=8cm]{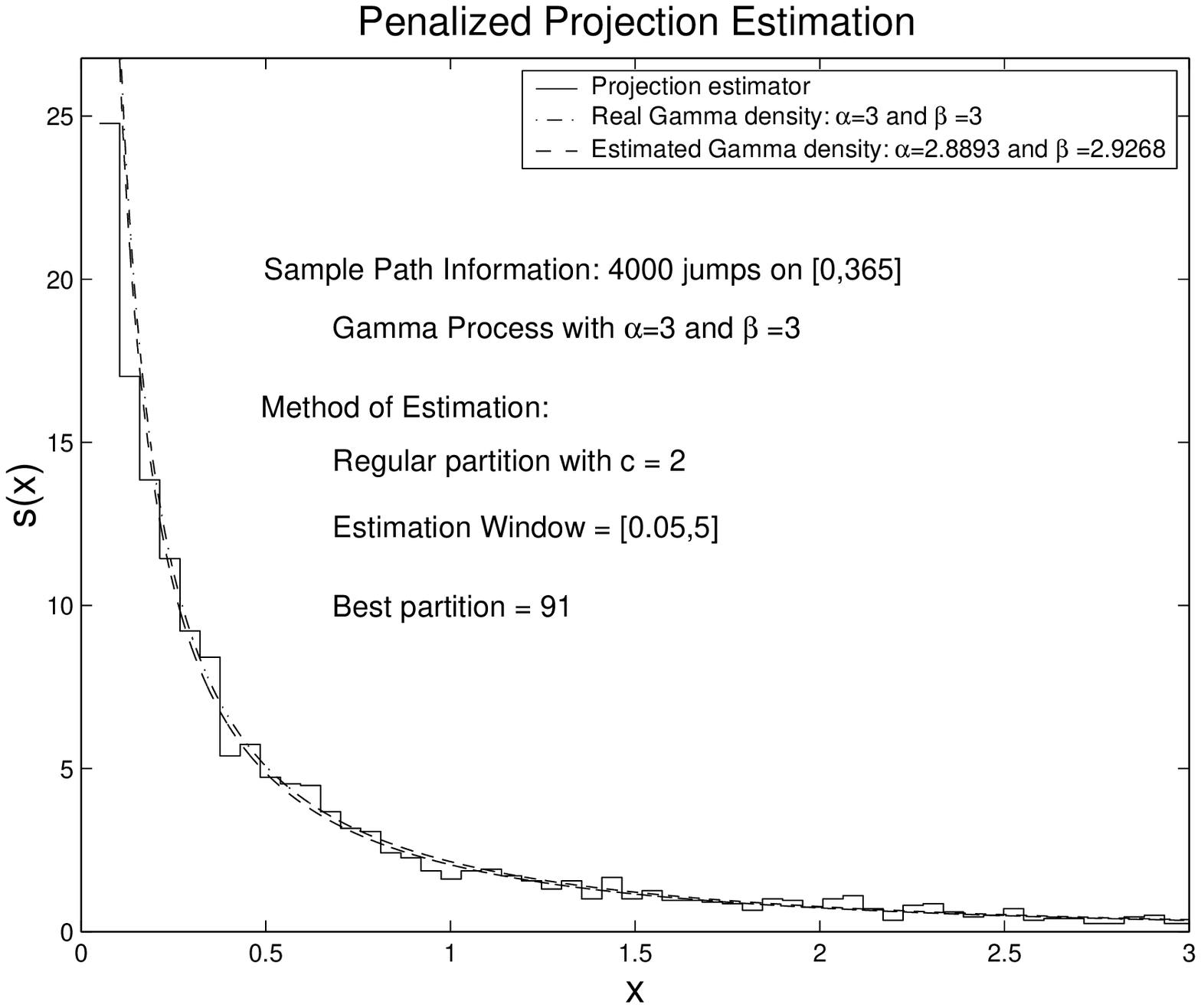} 
	\par}
	\caption{\label{EstGamma3b} Penalized projection estimation of
	$\frac{3}{x}\exp\left(-\frac{x}{3}\right)$. }
\end{figure}

\begin{figure}[htp]
	{\par \centering
	\includegraphics[width=8cm,height=8cm]{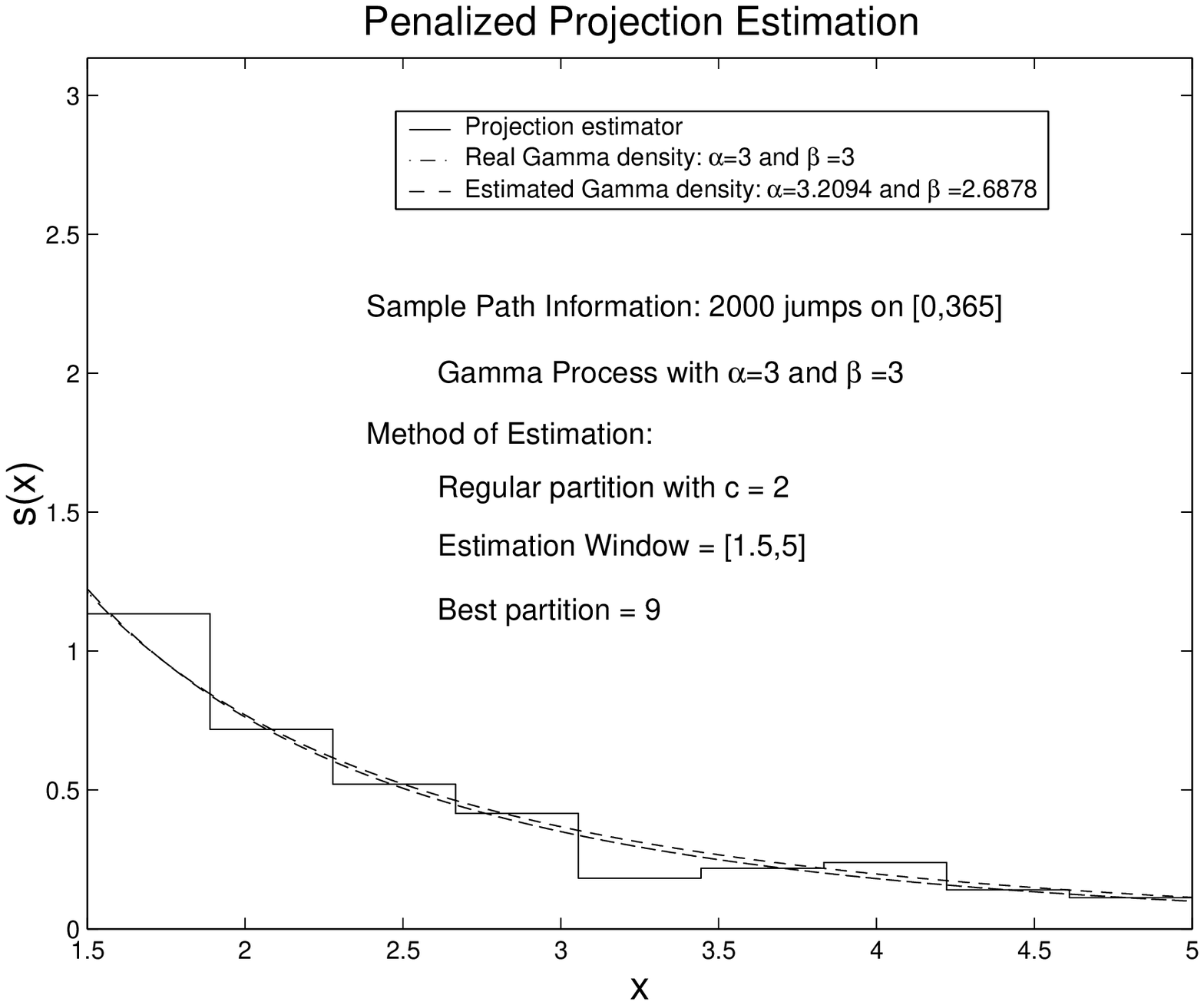} 
	\par}
	\caption{\label{EstGamma3c} Penalized projection estimation of
	$\frac{3}{x}\exp\left(-\frac{x}{3}\right)$.}
\end{figure}

Approximate histogram estimation on regular partitions is less successful in case of
high activity levels. This problem is particularly evident when we have in
addition heavy tails in the L\'evy density. For instance, if $\alpha = 3$
and $\beta =3$, the method requires a large sample size to 
satisfactorily retrieve the behavior around the origin (see Figures \ref{EstGamma3a}
and \ref{EstGamma3b}).
For $2000$ jumps, the least square estimates are $\hat{\alpha}=1.87$ and $\hat{\beta}=4.45$, 
while the estimates are $\hat \alpha = 2.8893$
and $\hat\beta = 2.9268$ for twice as many jumps. 
The maximum likelihood estimators
based on the increments of time length 0.5 are 2.4134054 for
$\alpha$ and 3.30971 for $\beta$ when the approximating 
process is made out of 2000 jumps,
while when the process is approximated using $4000$ jumps, 
these estimates are 
$2.8281$ and $3.1007$ for $\alpha$ and $\beta$, respectively.
We also notice in our experiments that the estimates for the first simulation improve 
considerably if the window of estimation is taken ``far away'' from the origin 
(for example, $\hat \alpha = 3.20944$ and $\hat \beta = 2.68775$ on $[a,b]=[1.5,5]$; see Figure
\ref{EstGamma3c} ).

\subsubsection{Regularized projection estimation around the origin}
We present another way to estimate the Gamma L\'evy density even around the
origin based on the regularization technique described in Section \ref{BscMethod}. 
The key observation is the following: the Gamma L\'evy measure (\ref{LevyDstyGamma})
can be written as 
\begin{equation}
\label{LevyDstyGamma2}
\nu(dx) =\alpha x \exp\left(-\frac{x}{\beta}\right) \eta(dx),
\end{equation}
where $\eta(dx)=\frac{1}{x^{2}}dx$. Then, $s(x)\equiv \alpha x
\exp\left(-\frac{x}{\beta} \right)$ is square integrable with respect to
$\eta$, opening the possibility to use the projection estimation of $s$ on
a linear space ${\cal S}$ of $L^{2}\left((0,\infty),\eta\right)$. Once an
estimator $\hat s$ for $s$ has been obtained, $\hat p$ defined by
$\hat{p}(x)=\hat{s}(x)/x^{2}$ can work as an estimator for the L\'evy
density $p(x)\equiv\alpha\exp{(-x/\beta)}/x$. In the terminology introduced in 
Section \ref{BscMethod}, $\eta$ is a regularizing measure for the Gamma 
L\'evy density $p$, and $s$ is the corresponding regularized L\'evy density
(see Definition \ref{Regularization}).

Let us specify this method for the linear model
	\[
		{\cal S}_{{\cal C}}=
		\left\{  
		f(x) \equiv c_{1} x \chi_{[x_{0},x_{1})}(x) + \sum_{i=2}^{m}c_{i}\; \chi_{[x_{i},x_{i+1})}(x):
		c_{1},\dots,c_{m} \in \bbr\right\},
	\]
where  ${\cal C}:0 = x_{0} < x_{1} < \dots < x_{m} = b $ is a partition of
a chosen interval $D=[0,b]$. The projection estimator, say $\hat s_{{\cal C}}$, onto ${\cal S}_{{\cal
C}}$, under the standard inner product of $L^{2}\left((0,\infty),\eta\right)$, takes on the value
\[
\hat s _{{\cal C}} (x) = x  \frac{1}{T x_{1}} \sum_{t\leq T} \Delta X (t) {\rm
I} \left[ \Delta X(t) <
x_{1} \right],
\]
if $x < x_{1}$, while if $x_{i-1}\leq x < x_{i}$, for some $i \in \{2,\dots,m\}$, then 
\[
\hat s_{{\cal C}}(x) = \frac{x_{i-1}x_{i}}{T(x_{i} - x_{i-1})} {\cal J}([x_{i-1},x_{i})).
\]
We shall use the penalty function of Theorem \ref{PreOracleIneq} part (b) 
to perform the model selection procedure. 
That is, among different partitions ${\cal C}$
that satisfy
\[
D_{{\cal C}} = \max \left\{ \frac{1}{x_{1}},\frac{x_{2}x_{1}}{x_{2}-x_{1}}, \dots,
\frac{x_{m}x_{m-1}}{x_{m}-x_{m-1}}\right\} \leq T,
\]
we choose the projection estimator $\hat s_{{\cal C}}$ that minimize 
\begin{align}
\label{PenConstrValue2} 
\gamma(\hat s_{{\cal C}}) + \hat V_{{\cal C}} &=
\frac{1}{T^{2}}\sum_{i=2}^{m} \frac{x_{i}x_{i-1}}{x_{i}-x_{i-1}} 
 \left[ c \; {\cal J}([x_{i-1},x_{i})) 
-\left({\cal J}([x_{i-1},x_{i}))\right)^{2}\right]  \nonumber \\
 &\;\; + \frac{c}{T^{2}x_{1}} 
   \sum_{{t\leq T :}\atop{\Delta X(t) < x_{1}}} (\Delta X (t))^{2}   -
 \frac{1}{x_{1}}\left(\sum_{{t\leq T :}\atop{\Delta X(t) < x_{1}}} \Delta X
(t)\right)^{2}.\nonumber
\end{align}
The previous formulas are found directly from the definitions and results given in
Section  \ref{BscMethod} (see for instance formulas (\ref{Contrast}), 
(\ref{DefPE}), (\ref{DefOfDm}), and (\ref{DefV})).


\begin{remk} Observe that the previous procedure is appropriate to estimate the density
function $s(x)=\frac{\alpha}{x}\exp(-\frac{x}{\beta})$ around the origin as far as 
\[
\hat{\alpha} \equiv \frac{1}{T x_{1}} \sum_{t\leq T} \Delta X (t) {\rm
I} \left[ \Delta X(t) <
x_{1} \right],
\]
is a good estimator of $\alpha$. It is not hard to check that
the bias of $\hat{\alpha}$ tend to zero as $x_{1}\searrow 0$. However, the
variance of $\hat{\alpha}$ converges to $\frac{\alpha}{2T}$, suggesting that
the method works better when $T$ is ``large'' and $\alpha$ is ``small''.
\end{remk}

\begin{figure}[hbp]
	{\par\centering
	\includegraphics[width=8cm,height=8cm]{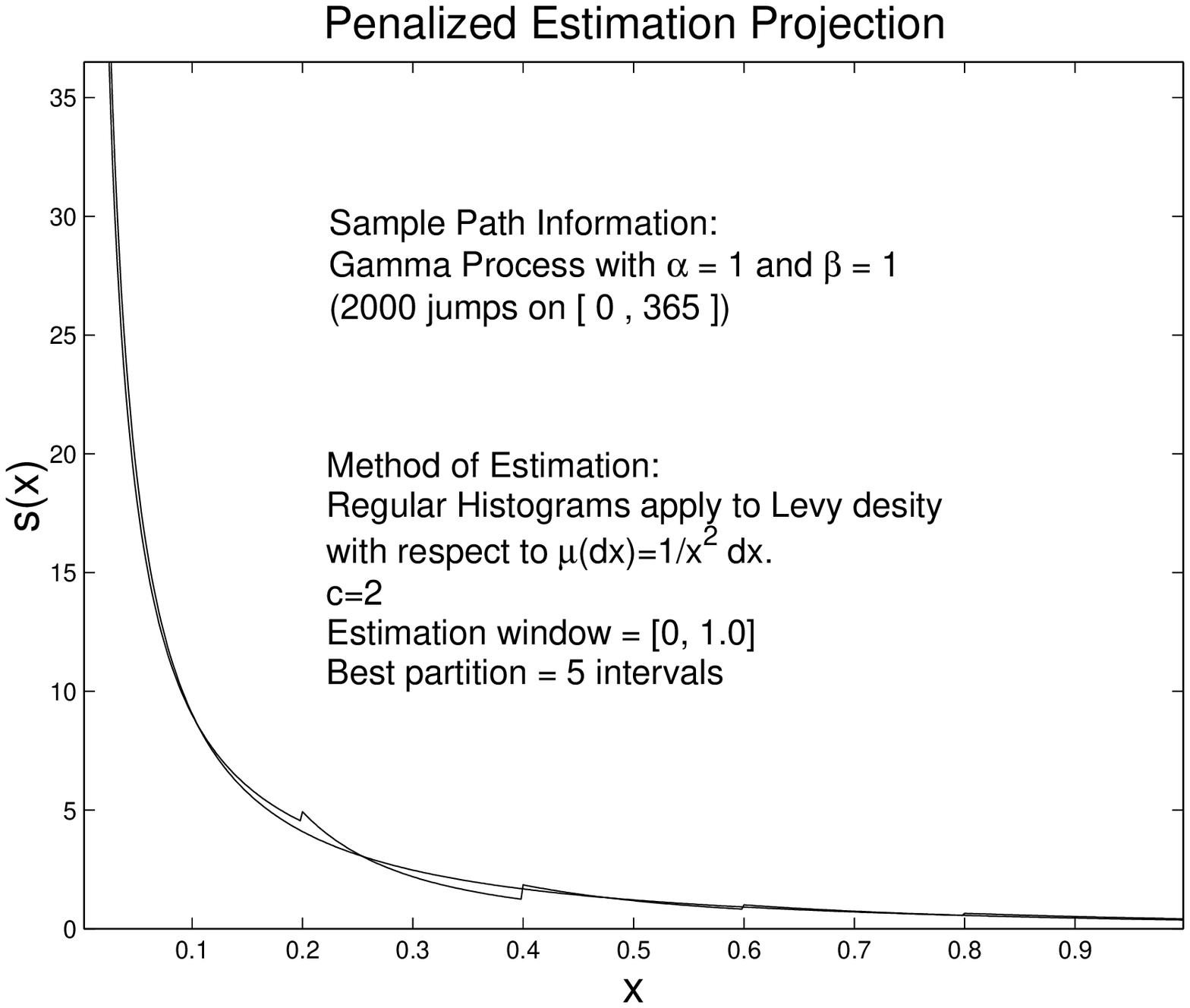}
	 \par}
	\caption{\label{EstGamma1_Mth2} Regularized penalized projection 
	estimation of
	$\frac{e^{-x}}{x}$.}
\end{figure}

We apply the above method to the simulated L\'evy process used in Figure \ref{EstGamma1}; i.e. a Gamma process with $\alpha =1$ and $\beta =1$. Figure \ref{EstGamma1_Mth2} shows the
estimator $\hat p _{2}(x) = \hat s(x)/x^{2}$ and the actual L\'evy density
$p(x)=\exp(-x)/x$ for $x \in [0.02,1]$ (we used regular partitions on
$[0,1]$). From Figure \ref{EstGamma1}, the improvement is notorious,
 and moreover, we accomplish a good estimation around the
origin of $\hat p_{2}(x)=0.9/x$, for $x \in [0,0.2)$.

\begin{figure}[htp]
	{\par \centering
	\includegraphics[width=8cm,height=8cm]{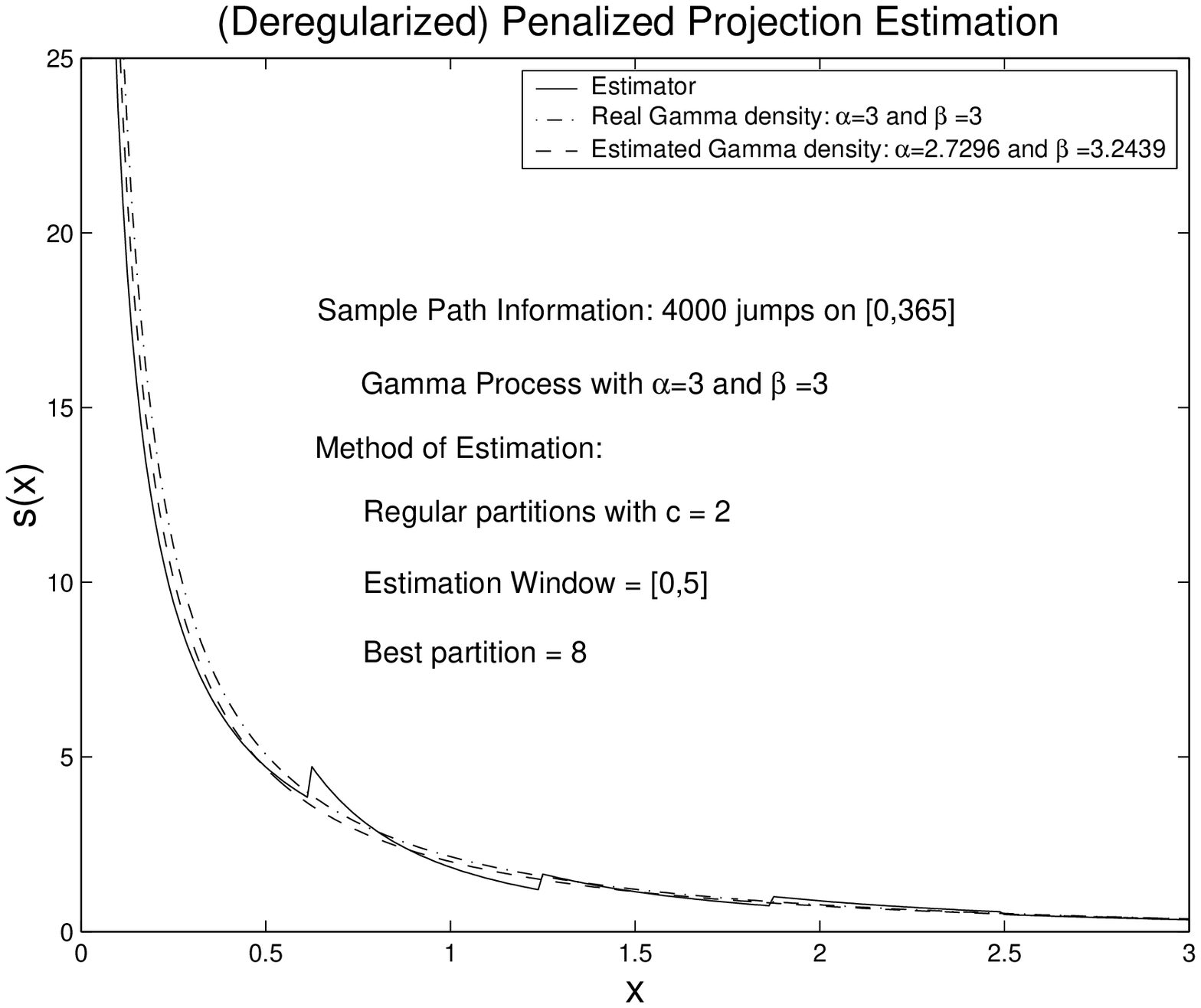} 
	\par}
	\caption{\label{EstGamma6_Mth2a} 
	Regularized penalized projection estimation of
	$\frac{3}{x}\exp\left(-\frac{x}{3}\right)$.}
\end{figure}

\begin{figure}[htp]
	{\par \centering
	\includegraphics[width=8cm,height=8cm]{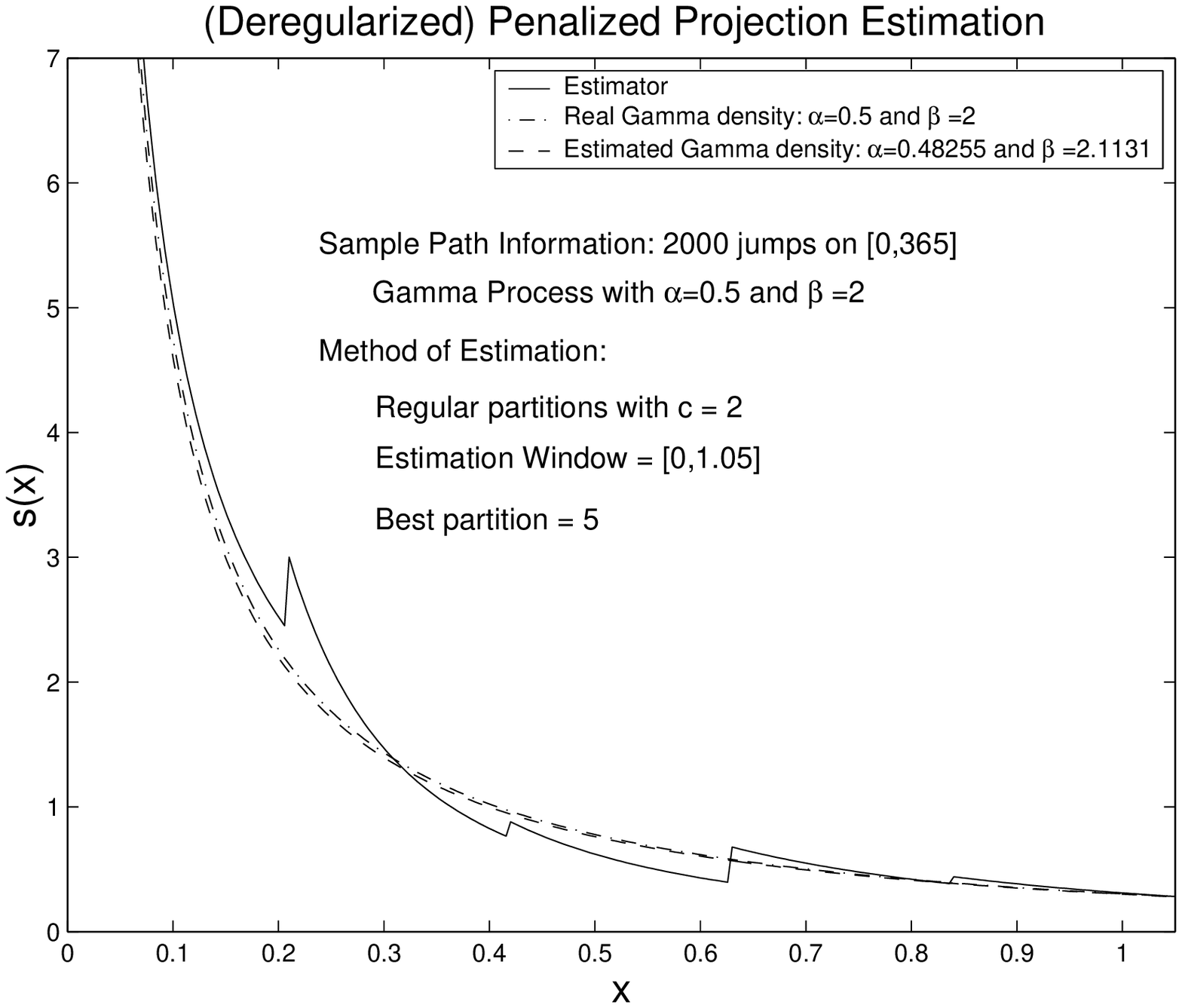} 
	\par}
	\caption{\label{EstGamma6_Mth2b} 
	Regularized penalized projection estimation of
	$\frac{1}{2x}\exp\left(-\frac{x}{2}\right)$.}
\end{figure}

This regularization procedure was also applied to the 
simulations of the Gamma L\'evy processes with 
$(\alpha =3, \beta=3)$ and with $(\alpha =1/2, \beta=2)$.
The results are plotted in Figures \ref{EstGamma6_Mth2a} and  
\ref{EstGamma6_Mth2b} below 
(compare with Figures \ref{EstGamma2b}  and \ref{EstGamma3b}). 
We observe an improvement on both sample data. 
For instance, for $\alpha=\beta =3$, 
the nonparametric estimator $\hat s(x)/x^{2}$ combined with a 
method of least-squares errors estimate  
$\alpha$ by $2.7296$ and $\beta$ by $3.2439$. 
Similarly, when $\alpha=0.5$ and $\beta=2$, 
least-square errors estimates $\hat\alpha=0.4825$ 
and $\hat\beta=2.1131$.

\subsubsection{Performance of projection estimation based on 
finitely many observation}	
In this part, we study the performance of the (approximate) 
projection estimators introduced in Section \ref{SecAproxPoissnInt}, 
and formally stated in Section \ref{EstimationmMethod}. 
Namely, the method obtained by approximating the 
Poisson process of jumps ${\cal J}$ by 
\[
	{\cal J}^{n} (\cdot)=
	\sum_{i=1}^{n} \delta_{(t_{i}^{n}J_{i})}(\cdot),
\]
where $J_{i}$ is the $i^{th}$ increment of $X$ from 
$t_{i-1}^{n}$ to $t_{i}^{n}$ and $t_{i}^{n} = i T/n$. 
The time span between increments is denoted by $\Delta t = T/n$.
Again, the considered estimators are histograms as defined in 
Section \ref{specifModel} and applied in Section \ref{NumResults}.
	
\begin{figure}[htp]
	{\par \centering
	\includegraphics[width=8cm,height=8cm]{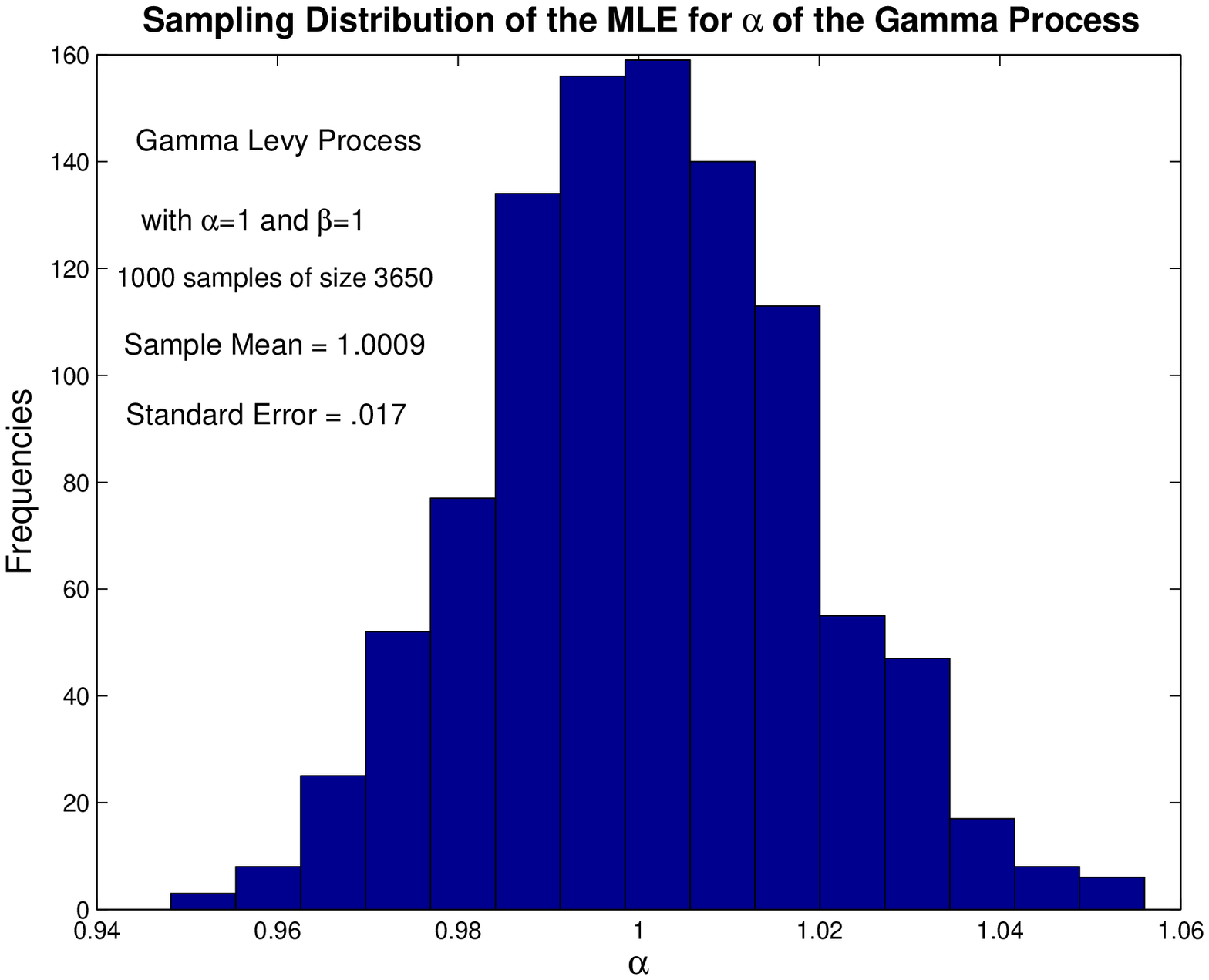} 
	\par}
	\caption{\label{SmplDistr_Gamma_f1} 
	Sampling Distribution for the MLE of the $\alpha$ of 
	the Gamma L\'evy Process.}
\end{figure}

\begin{figure}[htp]
	{\par \centering
	\includegraphics[width=8cm,height=8cm]{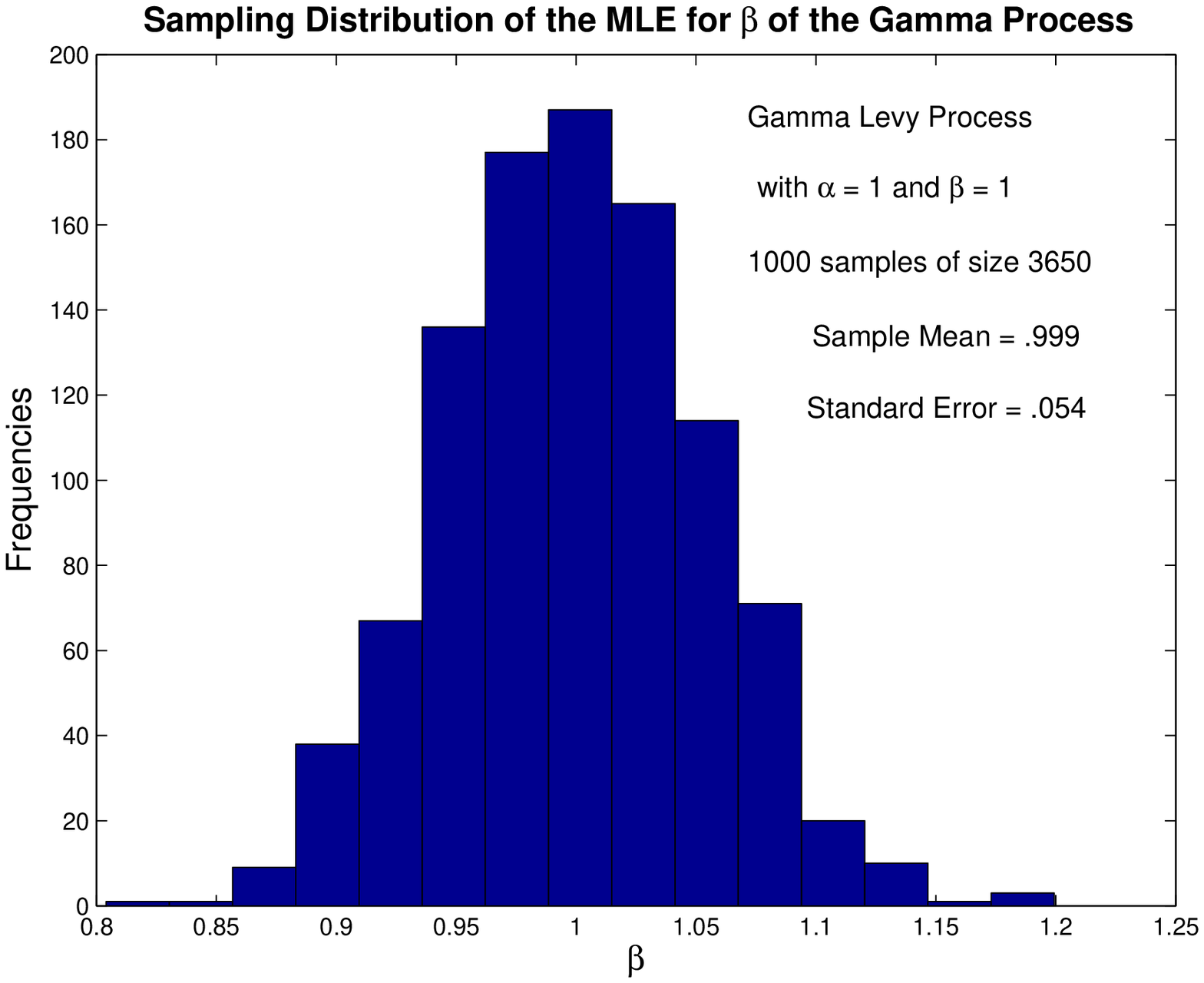} 
	\par}
	\caption{\label{SmplDistr_Gamma_f2} 
	Sampling Distribution for the MLE of the $\beta$ of 
	the Gamma L\'evy Process.}
\end{figure}

\begin{figure}[htp]
	{\par \centering
	\includegraphics[width=8cm,height=8cm]{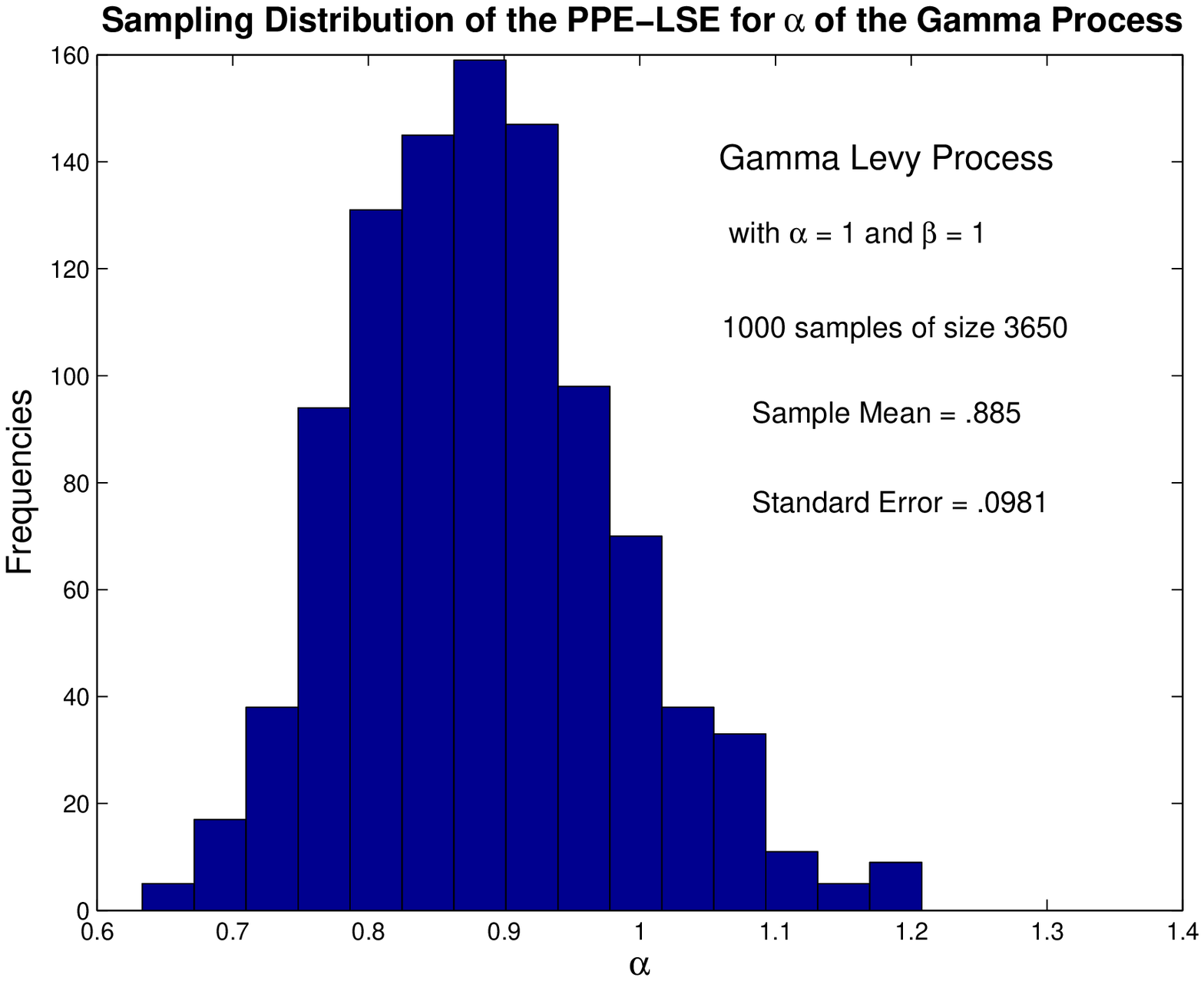} 
	\par}
	\caption{\label{SmplDistr_Gamma_f3} 
	Sampling Distribution for the Estimates of the $\alpha$
	of a Gamma L\'evy process 
	obtained from the PPE and the LSE method.}
\end{figure}

\begin{figure}[htp]
	{\par \centering
	\includegraphics[width=8cm,height=8cm]{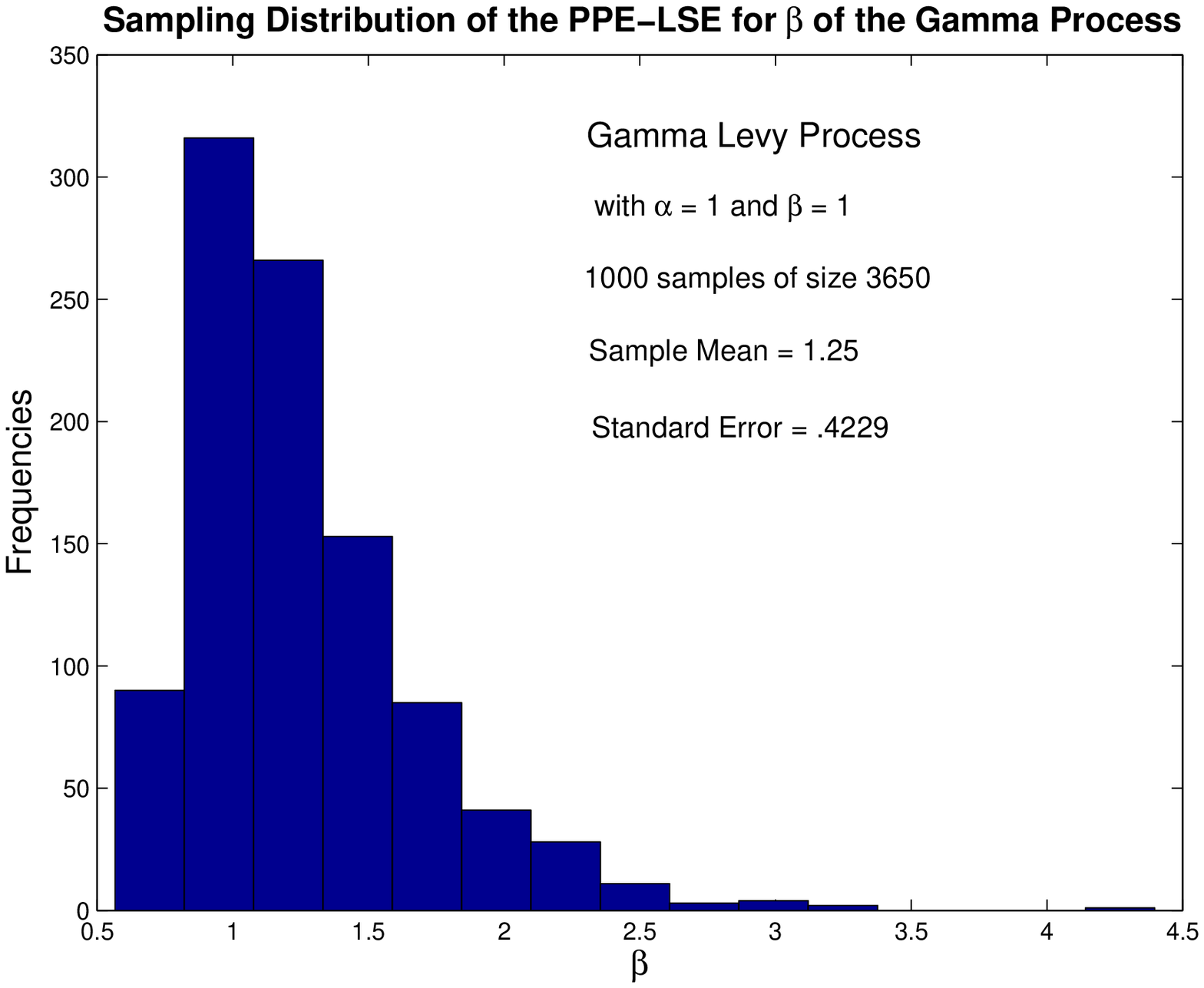} 
	\par}
	\caption{\label{SmplDistr_Gamma_f4} 
	Sampling Distribution for the Estimates of the $\beta$
	of a Gamma L\'evy process 
	obtained from the PPE and the LSE method.}
\end{figure}

Table \ref{TableEst1} compares the (approximate) penalized projection 
estimators with least-square errors (PPE-LSE) 
to the maximum likelihood estimators (MLE) for the Gamma 
L\'evy process with $\alpha=\beta=1$ using different 
time spans $\Delta t$.
We also consider two types of simulations: jump-based and 
increment-based.
The method based on jumps uses series representation with 
$n=36500$ jumps occurring during the time period $[0,365]$
(notice that if we think of 365 as days, the number of jumps 
corresponds to a rate of about 1 jump every 5 minute). 
The increment-based method is a \emph{discrete skeleton} with 
mesh of $0.001$.
Notice that maximum likelihood estimation does not do well for small 
time spans when the approximate sample path is based on jumps. 
On the other hand, penalized projection estimation does not 
provide good results for long time spans when the approximate 
sample path is based on increments. The sampling distributions of 
the MLE for $\alpha$ and $\beta$ are shown in Figures 
\ref{SmplDistr_Gamma_f1} and \ref{SmplDistr_Gamma_f2} 
in the case of $\Delta t=0.1$.
On the other hand, the sampling distributions of the estimates 
for $\alpha$ and $\beta$ obtained from fitting the PPE are given 
in Figures \ref{SmplDistr_Gamma_f3} and \ref{SmplDistr_Gamma_f4}. 
Even though, the MLE are much more superior, 
the estimates based PPE have good performance considering that they
are model-free. 
\begin{table}[!htbp]
{\par\centering		
\begin{tabular}{| c | c | c | c | c |  c | c | c | c |} 
\hline 
\hspace{1 cm} &
 \multicolumn{4}{| c |}{Jump-based Simulation} &
 \multicolumn{4}{| c |}{Increment-based Simulation} \\
 \hline
$\Delta t$ &
 \multicolumn{2}{| c |}{PPE-LSE} &
  \multicolumn{2}{| c |}{MLE} &
 \multicolumn{2}{| c |}{PPE-LSE} &
  \multicolumn{2}{| c |}{MLE} \\
  \hline
 1  & 1.01& 1.46 & .997 & .995 & .73& 1.78 & 1.09 & .99 \\
 \hline
 0.5 & 1.03 & 1.09 & .972 & .978 & .9 & 1.49 & 1.01 & 1.06 \\
 \hline
 0.1 & .944 & .995 & 1.179 & .837 & .923 & 1.03 & .989 & 1.09 \\
 \hline
 0.01 & .969 & .924 & 6.92 & .5  & .955 & 1.019 & .9974 & 1.083 \\
 \hline
\end{tabular}
\par}
\caption{\label{TableEst1}Estimation of a L\'evy Gamma process with $\alpha =\beta = 1$. Two types of simulation are considered: series-representation based and increments-based. The estimations are based on equally
spaced sampling observation at the time span $\Delta t$. Results for the approximate penalized projection estimators with least-squares errors, and for the maximum likelihood estimators are given. }
\end{table}

\subsection{Estimation of variance Gamma processes.}\label{VGM}
\subsubsection{The model}
Variance Gamma processes were proposed in \cite{Madan1} as substitutes 
to the Brownian Motion in the Black-Scholes model. There are two
useful representations for this type of processes. In short, a variance Gamma
process $X=\{ X(t)\}_{t \geq 0}$ is a Brownian motion with drift, time
changed by a Gamma L\'evy process. Concretely, 
\begin{equation}
\label{TimeChangeDef}
X(t)= \theta U(t) + \sigma W(U(t)),
\end{equation}
where $\{ W(t)\}_{t\geq 0}$ is a standard Brownian motion, $\theta \in
\mathbb{R}$, $\sigma>0$, and $U=\{U(t)\}_{t \geq 0}$ is an independent Gamma L\'evy process with 
density at time $t$ given by  
\begin{equation}
\label{DstyRandomTime}
f_{t}(x) = \frac{x^{t/\nu-1}\exp\left(-\frac{x}{\nu}\right)}
       {\nu^{t/\nu}\Gamma\left(\frac{t}{\nu}\right)}.
\end{equation}
Notice that ${\rm E}\left[U(t)\right]=t$ and ${\rm
Var}\left[U(t)\right]=\nu t$; therefore, the random clock $U$ has a ``mean rate'' of 
one and a ``variance rate'' of $\nu$. There is no loss of generality in
restricting the mean rate of the Gamma process $U$  to one since, 
as a matter of fact, any process of the form 
\[
	\theta_{1} V(t)+\sigma_{1} W(V(t)),
\]
where $V(t)$ is an arbitrary Gamma L\'evy process, $\theta_{1}\in\mathbb{R}$, and $\sigma_{1}>0$, has the same law as a process of the form (\ref{TimeChangeDef}) with suitably chosen $\theta$, $\sigma$, and $\nu$. This a consequence of the \emph{self-similarity}\footnote{namely, $\{ W(ct)\}_{t\geq 0}\stackrel{\frak{D}}{=}\{ c^{1/2}W(t)\}_{t\geq 0}$, for any $c>0$.\index{selfsimilarity}} property of Brownian motion and the fact that $\nu$ in (\ref{DstyRandomTime}) is a scale parameter. 

The process $X$ is itself a L\'evy process since Gamma processes are \emph{subordinators} (see Theorem 30.1 of \cite{Sato}). Moreover, it is not hard to check that ``statistically'' $X$ is the difference of two Gamma L\'evy processes (see 2.1 of \cite{Madan}): 
\begin{equation}
\label{DiffGammaDef}
 \{X(t)\}_{t\geq 0} \stackrel{\frak{D}}{=} \{X_{+}(t)-X_{-}(t)\}_{t\geq 0},
\end{equation}
where $\{X_{+}(t)\}_{t\geq 0}$ and  $\{X_{-}(t)\}_{t\geq 0}$ are Gamma
L\'evy processes with respective L\'evy measures
\[
\nu_{\pm}(dx)= \alpha
\exp\left(-\frac{x}{\beta_{\pm}}\right)dx,\;\;{\rm for}\; x>0.
\]
Here, $\alpha=1/\nu$ and 
\[
\beta_{\pm} = \sqrt{\frac{\theta^{2}\nu^{2}}{4}+\frac{\sigma^{2}\nu}{2}}
\pm \frac{\theta \nu}{2}.
\]
As a consequence of this decomposition, the L\'evy density of $X$ takes
the form
\[
s(x)=\left\{\begin{array}{ll}
	\frac{\alpha}{|x|}\exp\left(-\frac{|x|}{\beta_{-}}\right) &\:{\rm if}\; x<0, \\
	\frac{\alpha}{x} \exp\left(-\frac{x}{\beta_{+}}\right) &\:{\rm if}\; x>0,
	\end{array}\right. 
\]
where $\alpha>0$, $\beta_{-} \geq 0$, and $\beta_{+} \geq 0$ (of course, 
$\beta_{-}^{2}+\beta_{+}^{2}>0$). As in the case of Gamma L\'evy processes, 
$\alpha$ controls the overall jump activity, while 
$\beta_{+}$ and $\beta_{-}$ take respectively charge of
the intensity of large positive and negative jumps. In particular, 
the difference between $1/\beta_{+}$ and $1/\beta_{-}$ determines the
frequency of drops relative to rises, while their sum measures the frequency of 
large moves relative to small ones.	
\subsubsection{The simulation procedure}
The above two representations provide straightforward methods 
to simulate a variance Gamma model. 
One way will be to simulate the Gamma L\'evy processes 
$\{X_{+}(t)\}_{0\leq{t}\leq T}$ and $\{X_{-}(t)\}_{0\leq{t}\leq{T}}$ 
of (\ref{DiffGammaDef}) using the series representation method 
of Section \ref{SectSimProcGammaProc}. 
The other approach is to first generate random time change 
$\{U(t)\}_{0\leq{t}\leq{T}}$ of (\ref{TimeChangeDef}), 
and then construct a discrete skeleton from the increments 
$X(i\Delta{t})-X((i-1)\Delta{t})$, $i\geq{1}$. 
The increments of $X$ are simply simulated using normal 
random variables with mean and variances determined by the 
increments of $U$.

\subsubsection{The numerical results}
Notice that, from an algorithmic point of view, 
the estimation for the variance Gamma model
using penalized projection is not different from 
the estimation for the Gamma process. 
We can simply estimate both tails of the 
variance Gamma process separately. 
However, from the point of view of maximum likelihood 
estimation (MLE), the problem is numerically challenging.
Even though the marginal density functions have closed 
form expressions (see \cite{Madan1}), 
there are well-documented issues with MLE 
(see for instance \cite{Prause}). 
The likelihood function is highly flat for a wide 
range of parameters and good starting values as well as 
convergence are critical. Also, the separation of parameters 
and the identification of the variance Gamma process from 
other classes of the generalized hyperbolic L\'evy processes 
is difficult. In fact, difference between subclasses in terms 
of likelihood is small. It is important to mention that these 
issues worsen when dealing with ``high-frequency'' data.

\begin{figure}[htp]
	{\par \centering
	\includegraphics[width=8cm,height=8cm]{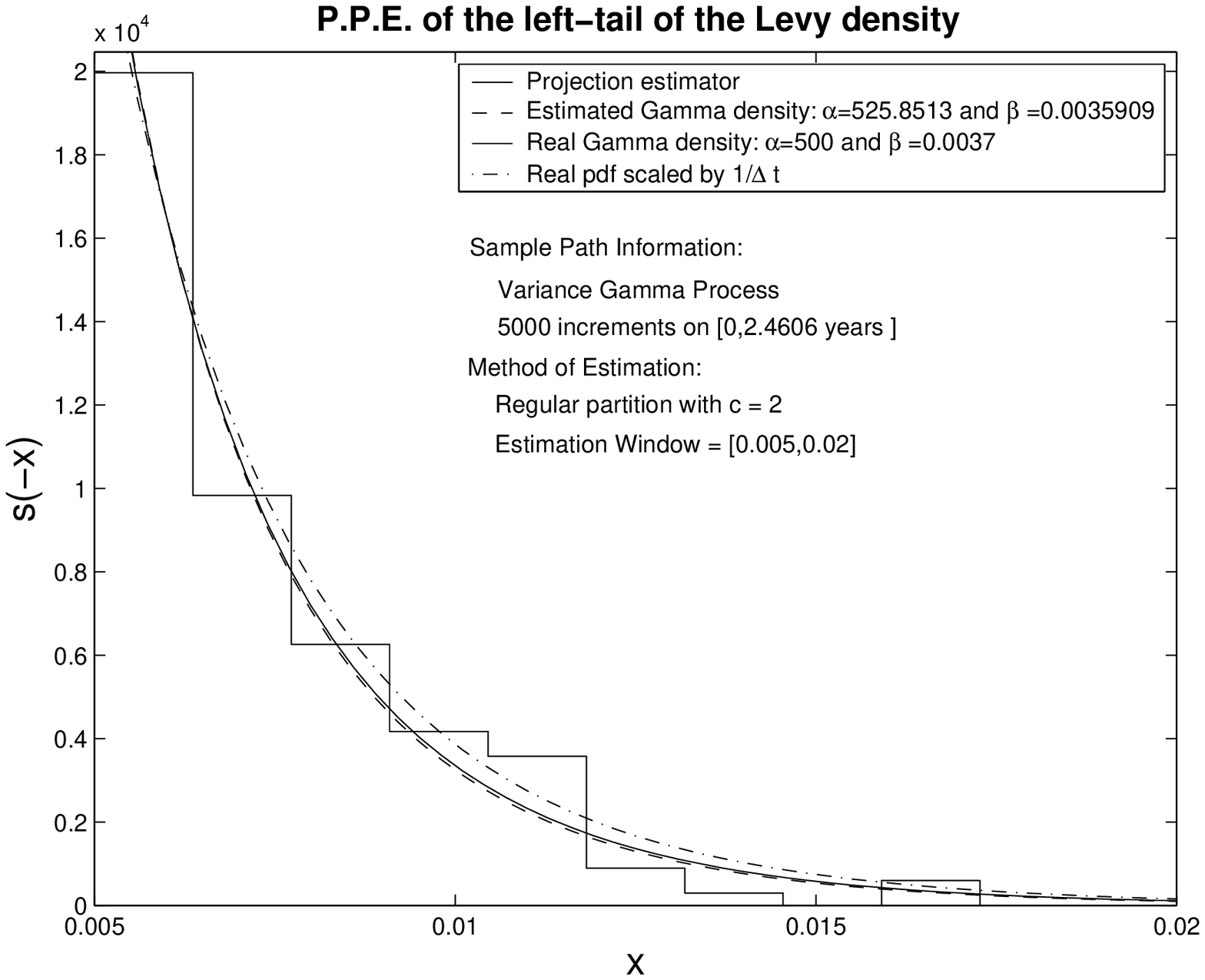} 
	\par}
	\caption{\label{ppe_VG_f1} 
	Penalized projection estimation of the left-tail of 
	the variance Gamma Levy density.}
\end{figure}

\begin{figure}[htp]
	{\par \centering
	\includegraphics[width=8cm,height=8cm]{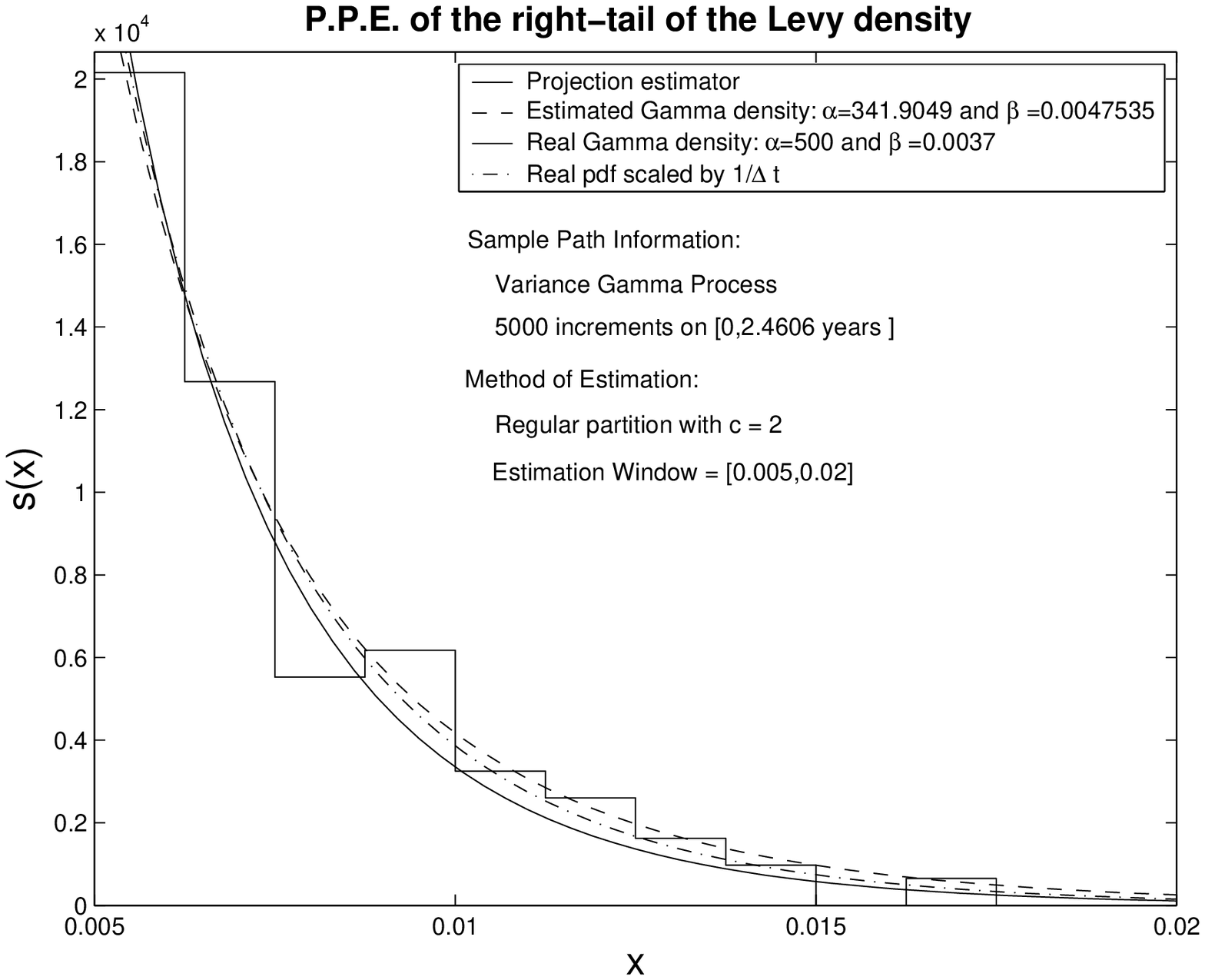} 
	\par}
	\caption{\label{ppe_VG_f2} 
	Penalized projection estimation of the right-tail of 
	the variance Gamma Levy density.}
\end{figure}

\begin{figure}[htp]
	{\par \centering
	\includegraphics[width=8cm,height=8cm]{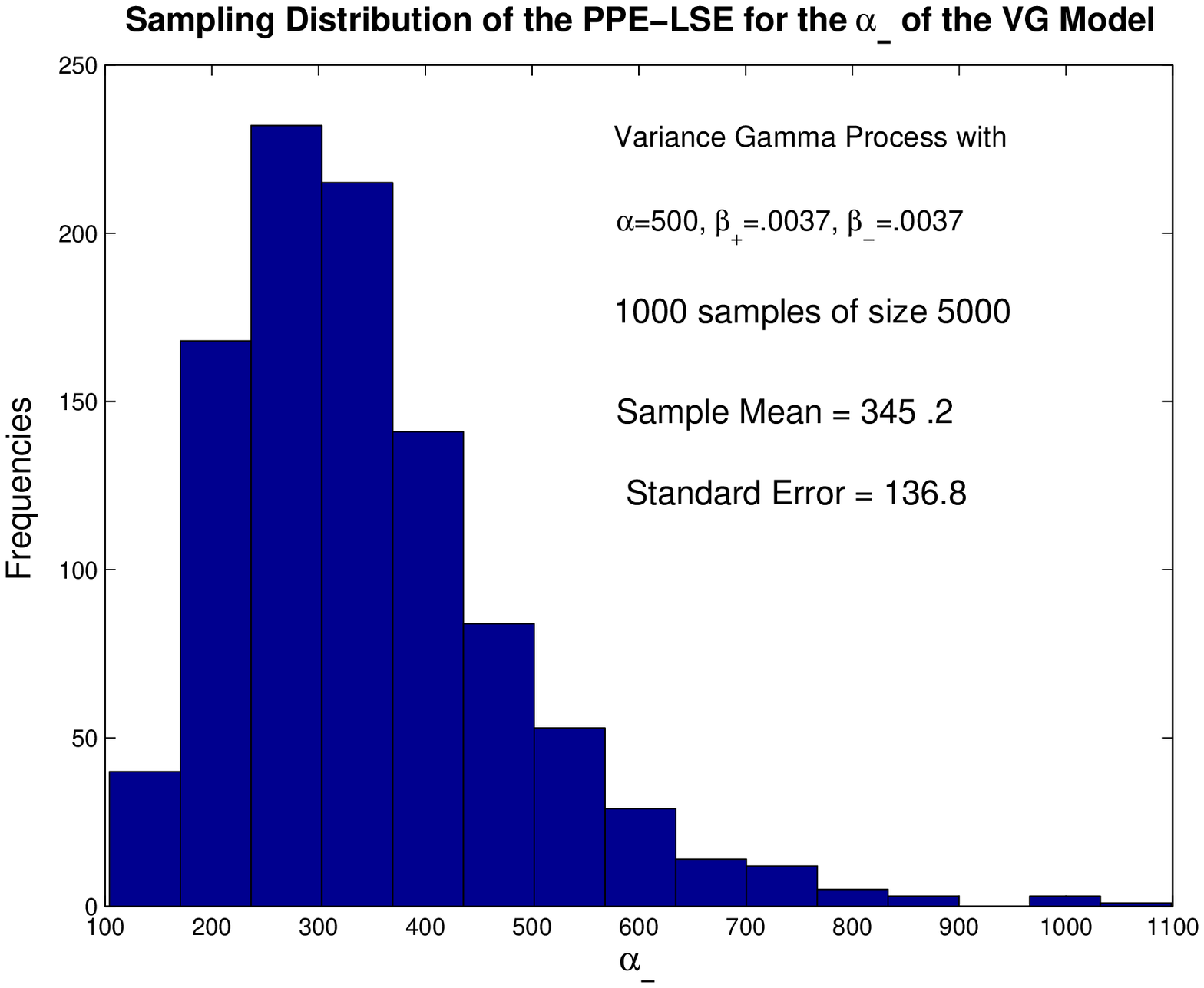} 
	\par}
	\caption{\label{SmplDistr_f2} 
	Sampling Distribution for the Estimates of $\alpha_{-}$
	obtained from the PPE and the LSE method.}
\end{figure}

\begin{figure}[htp]
	{\par \centering
	\includegraphics[width=8cm,height=8cm]{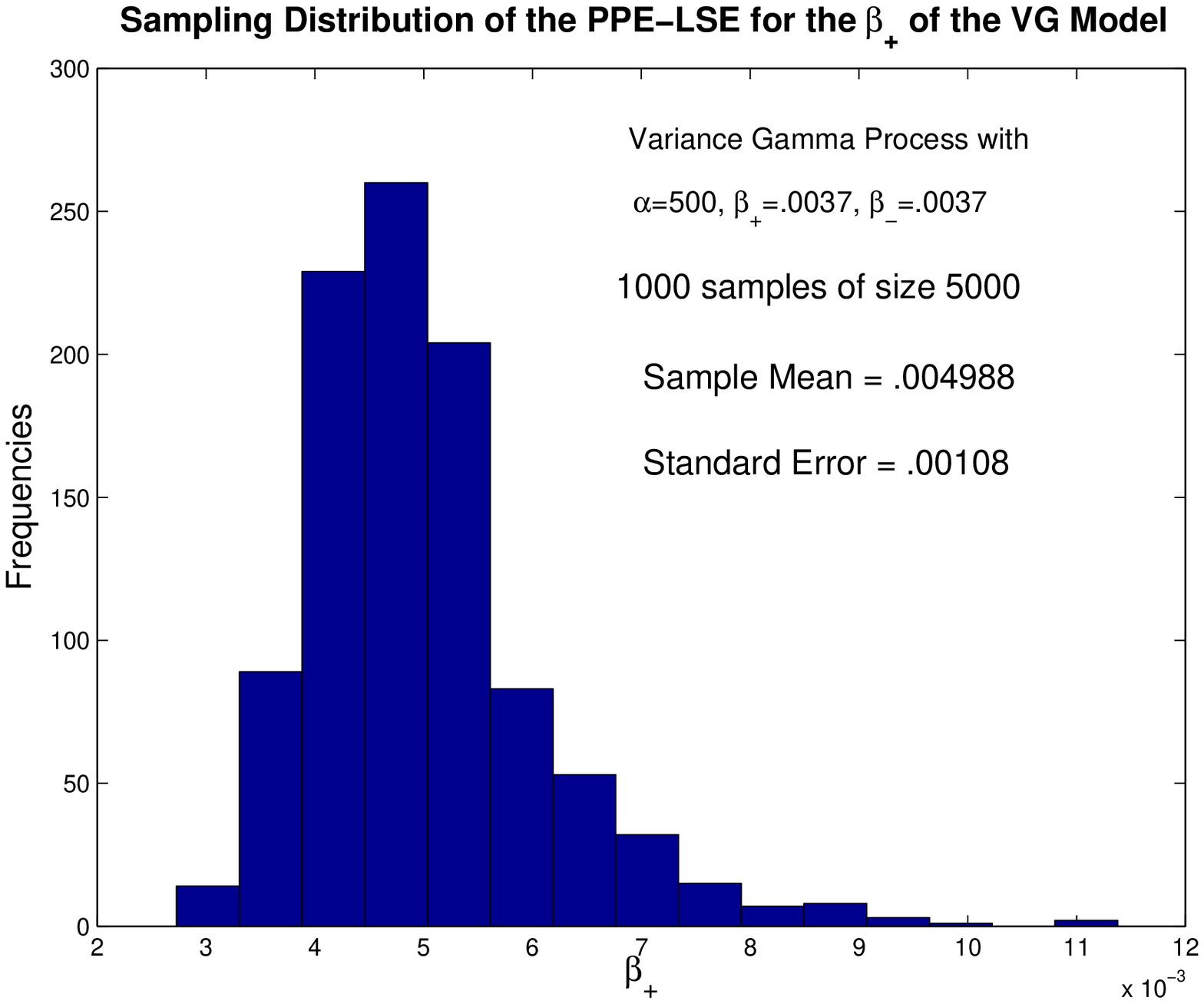} 
	\par}
	\caption{\label{SmplDistr_f3} 
	Sampling Distribution for the Estimates of $\beta_{+}$
	obtained from the PPE and the LSE method.}
\end{figure}

\begin{figure}[htp]
	{\par \centering
	\includegraphics[width=8cm,height=8cm]{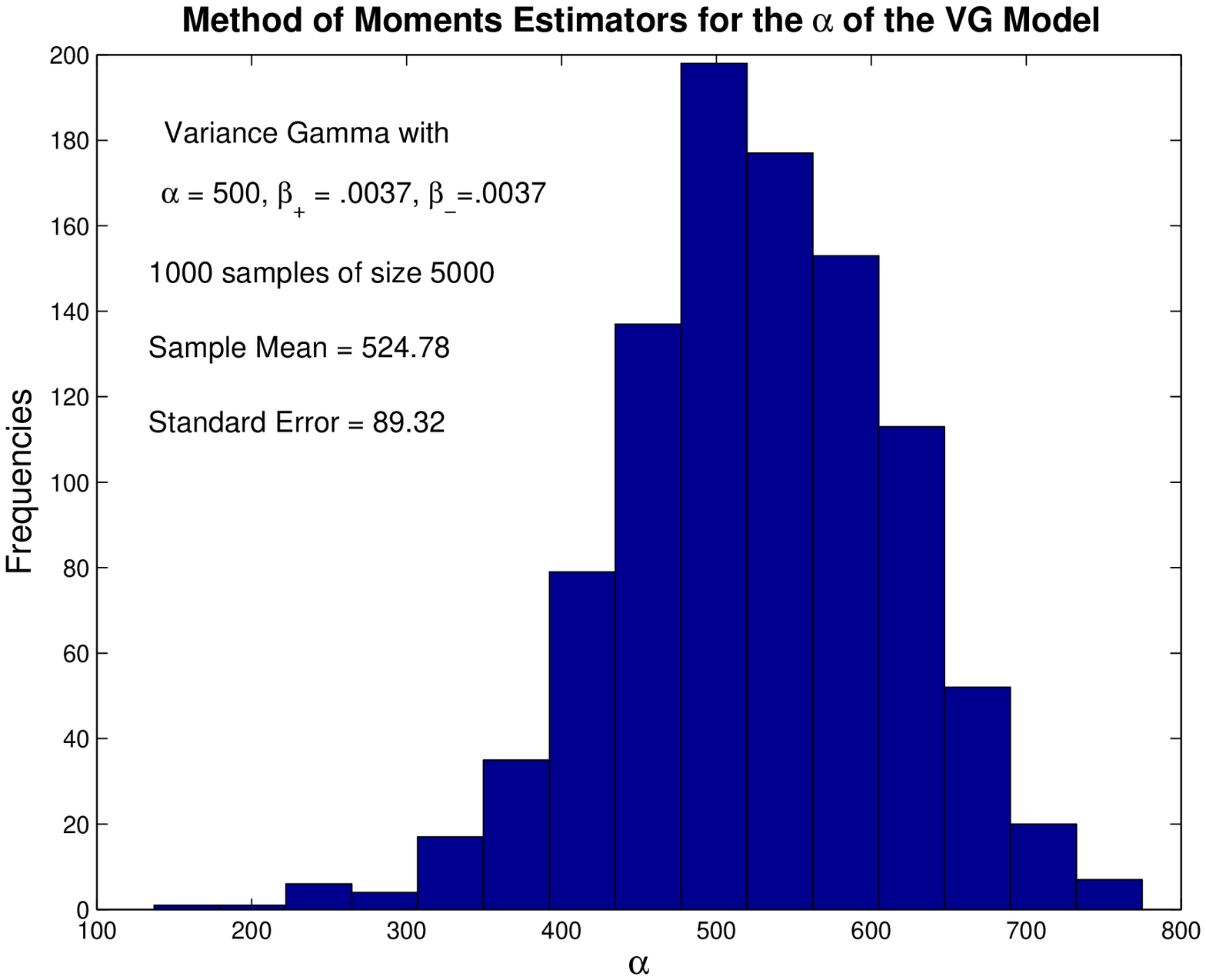} 
	\par}
	\caption{\label{SmplDistr_MME_f1} 
	Sampling Distribution for the Estimator of $\alpha$ 
	obtained by the Method of Moments.}
\end{figure}

\begin{figure}[htp]
	{\par \centering
	\includegraphics[width=8cm,height=8cm]{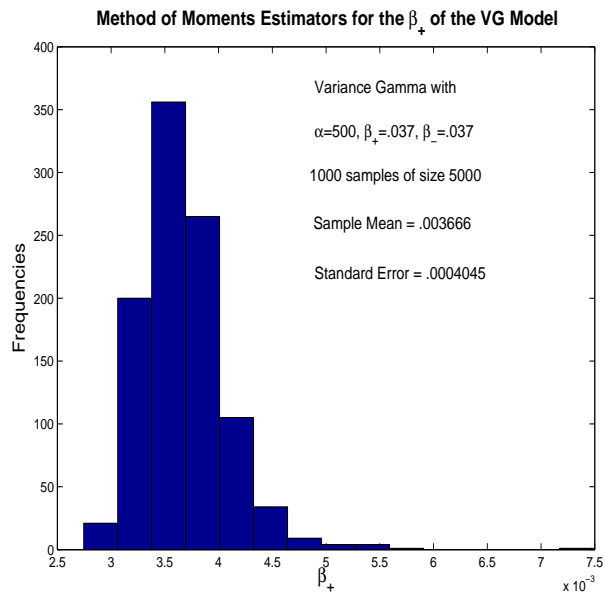} 
	\par}
	\caption{\label{SmplDistr_MME_f2} 
	Sampling Distribution for the Estimator of $\beta_{+}$ 
	obtained by the Method of Moments.}
\end{figure}

Let us consider a numerical example motivated by the empirical 
findings of \cite{Madan1} based on daily returns on the S\&P 
stock index from January 1992 to September 1994 
(see their Table I). 
Using maximum likelihood methods, 
the annualized estimates of the parameters for 
the variance Gamma model were reported to be 
$\hat{\theta}=-0.00056256$, $\hat{\sigma}^{2}=0.01373584$, 
and $\hat{\nu}=0.002$, from where we obtain 
%
$\hat\alpha = 500$, $\hat{\beta}_{+}=0.0037056$, and 
$\hat{\beta}_{-}=0.0037067$.
%
Figures \ref{ppe_VG_f1} and \ref{ppe_VG_f2} 
show respectively the left- and right- tails of the 
L\'evy density and their penalized projection estimators
as well as their corresponding 
best-fit variance Gamma L\'evy densities using a 
least-square method, and 
their marginal probability density functions 
(pdf) scaled by $1/\Delta t$ 
(the reciprocal of the time span between observations).
The estimation was based on $5000$ simulated increments with 
$\Delta{t}$ equal to one-eigth of a day. The figures seem quite 
comforting. To get a better picture, Figures \ref{SmplDistr_f2}
and \ref{SmplDistr_f3} show the 
sampling distributions of the  estimates for 
$\alpha_{-}$ and $\beta_{+}$ obtained from applying 
the least-square method to the penalized proyection estimators. 
The histograms are based on $1000$ samples of size $5000$ 
with $\Delta{t}=1/8$ of a day. This experiment shows clear, 
though not critical,  underestimation of the parameter 
$\alpha$ and overestimation of the parameters $\beta$'s. 
A simple method of moments (based on the first four moments)
yields better results (see Figures \ref{SmplDistr_MME_f1} and 
\ref{SmplDistr_MME_f2}).
Nonparametric methods are not free-lunches and usually 
the gain in robustness is paid by a lost in precision. 


\section{Concluding Remarks}

\begin{itemize}
\item[$\bullet$] 
In the present paper we have developed a new 
methodology for the estimation of the 
L\'evy density of a L\'evy process. 
Our methods are quite flexible in the sense that 
different type of estimating functions can be used; 
for instance, histograms, splines, trigonometric polynomials, 
and wavelets. The estimation is  model free, easily 
implementable, and suitable for ``high-frequency'' data. 
We prove that, based on continuous-time data, our procedures 
enjoy good asymptotic properties. \emph{Oracle inequalities} 
imply that, up to a constant, the procedure will 
achieve the best possible risk among the 
\emph{projection estimators}. Moreover, it is proved that 
penalized projection estimators on splines achieve
the optimal rate of convergence, from the minimax point of view,
on some classes of smooth L\'evy densities. Simulations
show good results in L\'evy models with infinite jump activity 
such as the variance Gamma model. 

\item[$\bullet$] 
Generalization of our procedures and results to
some multivariate L\'evy models can be readily obtained, 
since the results behind our construction 
have multivariate versions. Indeed, 
the L\'evy-It\^o decomposition of the sample paths, 
the concentration inequalities for compensated Poisson integrals,
the inference theory for Locally Asymptotically Normal 
distributions, and the short-term properties of the marginal 
distributions are valid in the multivariate setting.
More precisely, consider a L\'evy process 
${\bf X}=\{ \bfXt \}_{t \geq 0}$ on $\bbr^{d}$ 
with L\'evy measure $\nu$. Assume that, on a 
window of estimation 
$D\in\calB(\bbr^{d}\backslash{\{0\}})$, $\nu$ is absolutely 
continuous with respect to a reference measure $\eta$ and that
$s \equiv d\nu / d\eta$ is bounded with also
$\int_{D} s^{2}({\bf x})\eta(d{\bf x})<\infty$. 
Then, given a finite-dimensional subspace $\calS$ of 
$L^{2}\left((D,\eta)\right)$, the projection estimator of $s$
on $\calS$ is defined as in Section \ref{BscMethod} with 
$\calJ$ being the Poisson measure on $\bbr_{+}\times\bbr^{d}$
associated with the jumps of ${\bf X}$. Similarly, penalized 
projection estimators can be constructed, and the risk bound 
of Theorem \ref{PreOracleIneq}, 
along with the Oracle inequality (\ref{OracleIneq}),
are satisfied. 
The results of Sections \ref{SecAproxPoissnInt} 
and \ref{EstimationmMethod} are valid as well.
However, let us point out that further specifications of our 
methods for some semiparametric models are desirable. 
Important examples of these models include multivariate 
stable processes, and the tempered stable L\'evy processes, 
recently introduced in \cite{Rosinski:2002}.

\item[$\bullet$] We have concentrated here on the 
estimation of the jump part of the L\'evy process. 
It is natural to address the problem of estimating 
the continuous part too. In the one-dimensional 
case, this part is of the form $bt+\sigma W(t)$, where 
$\{ W(t)\}_{t\geq{0}}$ is a standard Brownian motion.
In the multivariate case, it is characterized by a vector 
${\bf b}$ and a symmetric nonnegative-definite matrix $\Sigma$.
There are several approaches to deal with the estimation of 
$\Sigma$, from moment based methods to 
methods based on high-frequency data. A simple approach is
to use the following functional limit:
\[
	\left\{ \frac{1}{\sqrt{h}} {\bf X}(ht)\right\}_{t\geq{0}} 
	\ld \{ {\bf Y}(t)\}_{t\geq{0}},\;\;h\rightarrow{0},
\]
where $\{ {\bf Y}(t) \}_{t \geq 0}$ is a 
centered Gaussian L\'evy process with variance-covariance 
matrix $\Sigma$. This result can be deduced from the proof 
of the uniqueness of the L\'evy-Khintchine representation
as in pp. 40 of \cite{Sato}. Another simple method will be to 
consider empirical versions of the moments:
\[
	\bbe\left[(X_{i}(t)-X_{i}(t))(X_{k}(t)-X_{k}(t))\right] =
	t \left(\Sigma_{i,k} + 
	\int_{\bbr_{0}^{d}} x_{i} x_{k} \nu(d\bfx )\right),
\]
provided that $\int_{\|\bfx\|>1} \|\bfx\|^{2}\nu(d\bfx )<\infty$
(see Section 25 \cite{Sato}).
Here, $X_{i}(t)$ and $x_{i}$ refer to the $i^{th}$ 
component of the vectors $\bfXt$ and $\bfx$, respectively,
while $\Sigma_{i,j}$ is the $(i,j)$ entry of $\Sigma$.
The second term on the left hand side of the above expression 
can be estimated 
using our estimators for $\nu$. 
In the one-dimensional case, another approach is 
to use ``threshold estimators'' of the form:
\[
	\sum_{k=1}^{n}\left(\Delta_{k}X\right)^{2} 
	{\bf 1}\left((\Delta_{k}X)^{2}\leq r(h)\right),
\]
 where $\Delta_{k}X\equiv X(t_{k}^{n})-X(t_{k-1}^{n})$ is the $i^{th}$ 
increment of the process and $r(h)$ is an appropriate cutoff function
(see \cite{Mancini} for details). For a class of semimartingales with 
finite jump activity, 
\cite{Barndorff:2003a} provides other methodology based on the 
\emph{bipower variation} (see also \cite{Barndorff:2003b}). 
In the case of L\'evy processes with finite jump activity, 
\cite{Ait_Sahalia} disentangles the difussion from 
jumps using maximum likelihood and the \emph{Generalized Method of Moments}.
On the other hand, the estimation of the
parameter ${\bf b}$ can be done by different methods. 
For instance, using the empirical version for 
\[
	\bbe\left[\bfX(t)\right] = 
	t \left({\bf b}+\int_{\|\bfx\|>1} \bfx \nu(d\bfx )\right), 
\]
valid if $\int_{\|\bfx\|>1} \|\bfx\|\nu(d\bfx )<\infty$. 
Another approach will be to estimate the ``drift'' 
${\bf b}_{0}\equiv {\bf b}-\int_{\|\bfx\|\leq1} \bfx \nu(d\bfx )$ 
(where the integration is component-wise) 
using the fact that
\[
	\bbp \left[ \lim_{h\rightarrow 0} \frac{1}{h} {\bf X}(h) 
	= {\bf b}_{0} \right]=1.
\]
The above result holds true if 
$\int_{\|\bfx\|\leq1} \|\bfx\|\nu(d\bfx )<\infty$ and $\Sigma=0$
(see \cite{Sato}). Even though our methods are valid for non
necessarily pure-jump L\'evy processes, it is expected that the 
presence of a diffusion component will reduce the efficiency 
(in terms of speed of convergence and accuracy) of our estimators. 
It would be interesting to study in greater detail this 
phenomenon.

\end{itemize}


\section{Main Proofs}


\subsection{Proof of the risk Bound}\label{SectionProofOracle}
We will break the proof of Theorem \ref{PreOracleIneq} into several
preliminary results.
	\noindent\begin{lmma} \label{BasicIneq} 
	For any penalty function ${\rm pen}: {\cal M} \rightarrow [0,\infty)$ and any 
	$m \in {\cal M}$, the penalized projection estimator $\tilde s$ satisfies 
	\begin{equation}
		\label{FundmIneq}
		\| s - \tilde s \|^{2} \leq
		\| s - s_{m}^{\bot} \|^{2} +  2 \chi_{\hat m}^{2}
		+2 \nu_{D} \left( s^{\bot}_{\hat m} - s^{\bot}_{m}\right) + {\rm pen}(m) - {\rm pen}(\hat m),
	\end{equation}
	where $\chi_{m}^{2}\equiv \| s_{m}^{\bot} - \hat{s}_{m}\|^{2}$ and where the functional
	$\nu_{D}: L^{2}\left((D,\eta)\right) \rightarrow \bbr$ is defined by 
	\begin{equation}
		\label{CompensatedIntegral}
		\nu_{D}(f) \equiv \iint\limits_{[0,T] \times D} f(x)\; \frac{{\cal J}(dt,dx)-s(x)\,dt\,\eta(dx)}{T}.
	\end{equation}
	\end{lmma}
	
The general idea to deduce (\ref{PseudoOracleIneq}) is to bound 
the unattainable terms of the right hand side of (\ref{FundmIneq})  
(namely $\chi^{2}_{\hat{m}}$ and  $\nu_{D} \left( s^{\bot}_{\hat m} - s^{\bot}_{m}\right)$) 
 by observable statistics. Then, the form of  ${\rm pen }(\cdot)$ will be determined by this observable statistics so that the right hand side in (\ref{FundmIneq}) does not involve $\hat{m}$. To carry out this plan, we use concentration inequalities for
 $\chi^{2}_{\hat{m}}$ and for the compensated Poisson integrals $\nu_{D}(f)$. The following result gives a concentration inequality for general compensated Poisson integrals. 
	\noindent\begin{prop} 
	\label{ConctrIneq1} \index{Poisson integrals!concentration inequality}
	\index{concentration inequality!Poisson functionals}
	Let $N$ be a Poisson process on a
	measurable space $({\rm V},{\cal V})$ with mean measure $\mu$ and let
	$f:{\rm V}\rightarrow \bbr$ be an essentially bounded measurable function satisfying
	$0<\int_{{\rm V}} f^{2}(v) \mu(dv)$ and 
	$\int_{{\rm V}} \vert f(v)\vert \mu(dv)<\infty$.
	Then, for any $u>0$,
	\begin{equation}\label{ConctrInqForIntgr}
		\bbp \left[ \int_{{\rm V}} f(v) (N(dv)-\mu(dv))\geq 
		 \| f \|_{L^{2}(\mu)}\sqrt{2 u} + \frac{1}{3} \| f \|_{\infty} u \right]\leq e^{-u},
	\end{equation}
	where $ \| f \|_{L^{2}(\mu)}^{2}= \int_{{\rm V}} f^{2}(v) \mu(dv)$.
	In particular, if $f:{\rm V}\rightarrow [0,\infty)$ then, for any $\epsilon >0$ and $u>0$,
	\begin{equation}\label{BoundForMean}
		\bbp \left[ (1+\varepsilon)\left(\int_{{\rm V}} f(v) N(dv)
		+\left( \frac{1}{2\varepsilon}+\frac{5}{6}\right)\|f\|_{\infty}u\right)\geq 
		\int_{{\rm V}} f(v)\mu(dv) \right]\geq 1-e^{-u}.
	\end{equation}
	\end{prop} 
For a proof of the inequality (\ref{ConctrInqForIntgr}), 
see \cite{Reynaud} (Proposition 7) or \cite{Houdre1}
(Corollary 5.1). 
Inequality (\ref{BoundForMean}) is a direct consequence 
of (\ref{ConctrInqForIntgr})
(see Section \ref{proofsChap3} for a proof).

The next result allow us to bound the Poisson functional
$\chi_{m}^{2}$. This results is essentially Proposition 9 of \cite{Reynaud}.
	\noindent\begin{lmma} 
	\label{IneqForChi} \index{concentration inequality!Poisson functionals}
	\index{variance term!concentration inequality}
	Let $N$ be a Poisson process on a
	measurable space $({\rm V},{\cal V})$ with mean measure 
	$\mu(d v)=p(v)\eta(dv)$ and intensity function 
	$p \in L^{2}({\rm V},{\cal V}, \eta)$. Let ${\cal S}$ be a finite dimensional subspace 
	of $L^{2}({\rm V},{\cal V}, \eta)$ with orthonormal basis $\left\{ \tilde{\varphi}_{1},
	\dots,\tilde{\varphi}_{d} \right\}$, and define 
	\begin{align}
		\hat p(v) &\equiv \sum_{i=1}^{d} \left( \int_{{\rm V}} \tilde{\varphi}_{i} (w) N(dw) \right)
		 \tilde{\varphi}_{i} (v)\\
		p^{\bot}(v) &\equiv \sum_{i=1}^{d} \left( \int_{{\rm V}} p(w)\tilde{\varphi}_{i}(w) \eta(dw) \right) 
		\tilde{\varphi}_{i} (v).
	\end{align}
	Then, $\chi^{2}({\cal S})\equiv \| \hat p - p^{\bot} \|_{L^{2}(\eta)}^{2}$ is such that for
	any $u>0$ and $\varepsilon>0$
	\begin{equation}\label{MainInqForChi}
		\bbp \left[ \chi({\cal S}) \geq (1+\varepsilon) \sqrt{\bbe \left[ \chi^{2}({\cal S})\right] } 
		+ \sqrt{2 k {\rm M}_{{\cal S}} u } + k(\varepsilon) {\rm B}_{{\cal S}} u \right]\leq e^{-u},
	\end{equation}
	where we can take $k=6$, $k(\varepsilon)=1.25+32/\varepsilon$, and where
	\begin{align}
		M_{{\cal S}}&\equiv \sup\left\{ \int_{{\rm V}} f^{2}(v)p(v)\eta(dv): f\in {\cal S}, \| f\|_{L^{2}(\eta)} =1 \right\}, \\
		B_{{\cal S}}&\equiv \sup\left\{ \| f\|_{\infty}: f\in {\cal S}, \| f\|_{L^{2}(\eta)} =1 \right\}.
	\end{align}
	\end{lmma} 
Following the same strategy as \cite{Reynaud}, the idea is to deduce a concentration inequality of the form
\[
	\bbp \left[ \| s - \tilde s \|^{2} \leq  C  \left( \| s - s^{\bot}_{m}\|^{2} +
	{\rm pen}(m)\right)+h(\xi)\right] \geq 1- C' e^{-\xi},
\]
for constants $C$ and $C'$, and a function $h(\xi)$ (all independent of $m$). This will prove to be enough in view of the following result (see Section \ref{proofsChap3} for a proof).
	\noindent\begin{lmma} \label{AuxLemma}\index{concentration inequality!bounding the mean from}
	Let $h:[0,\infty)\rightarrow \bbr_{+}$ be an strictly increasing function with continuous derivative and  such that $h(0)=0$ and
	 $\lim_{\xi	\rightarrow \infty} e^{-\xi} h(\xi) =0$. If $Z$ is random variable satisfying 
	\[
		\bbp \left[ Z \geq  h(\xi) \right] \leq K e^{-\xi}, 
	\]
	for every $\xi>0$, then
	\[
		\bbe Z \leq K \int_{0}^{\infty} e^{-u} h(u) du.
	\]
	\end{lmma} 
We are now in position to prove the main result of this section. 
Throughout the proof, we shall have to introduce various constants and 
inequalities that will hold with high probability. In order to clarify the role that the
constants play in these inequalities, we shall make some conventions and give to 
the letters $x$, $y$, $f$, $a$, $b$, $\xi$, ${\cal K}$, c, and $C$, with various sub- or 
superscripts, special meaning. The letters with $x$ are reserved to denote positive 
constants that can be chosen arbitrarily. The letters with $y$ denote arbitrary 
constants greater than $1$.
$f, f_{1}, f_{2}, \dots$ denote quadratic polynomials of a variable $\xi$ whose 
coefficients (denoted by $a's$ and $b's$) are determined by the values of the 
$x's$ and $y's$. The inequalities will be true with probabilities 
greater that $1-{\cal K}e^{-\xi}$, where ${\cal K}$ is determined by the values of the $x's$ and 
the $y's$. Finally, $c's$ and $C's$ are used for constants constrained by the $x's$ and $y's$. 
\emph{It is important to remember that the constants in a given inequality are meant only for  
that inequality}. The pair of equivalent inequalities below will be repeatedly invoked through
 the proof:
	\begin{equation}\label{EasyIneq}
		\begin{array}{l}
		{\rm (i)}\;\;2ab \leq x a^{2} + \frac{1}{x} b^{2},\;\;\; {\rm and}\\
		{\rm (ii)}\;\;(a+b)^{2} \leq \left(1+x\right)a^{2}+
		\left(1+\frac{1}{x}\right)b^{2},\;\;\;({\rm for}\; x>0).
		\end{array}
	\end{equation}

\noindent\textbf{Proof of Theorem \ref{PreOracleIneq}}:
We consider successive improvements of the inequality (\ref{FundmIneq}):

\noindent \textit{Inequality 1: 
For any positive constants $x_{1}$, $x_{2}$,  $x_{3}$,  
and $x_{4}$, there is a positive number ${\cal K}$ and an 
increasing quadratic function $f(\xi)$ (both independent of 
the family of linear models and of $T$) such that,
with probability larger than $1-{\cal K} e^{-\xi}$,
\begin{equation}\label{FundmIneq_v2}\begin{array}{ll}
	\| s - \tilde s \|^{2} &\leq \| s - s_{m}^{\bot} \|^{2} +
	2 \chi_{\hat m}^{2} +  
	2 x_{1} \|s^{\bot}_{\hat m} -s^{\bot}_{m} \|^{2}\\ 
	&\;\;\;+ x_{2}\frac{D_{\hat m}}{T} + 
	x_{3}\frac{D_{m}}{T}+x_{4} \frac{d_{\hat m}}{T} \\
 	&\;\;\;+\;{\rm pen}(m)-{\rm pen}(\hat m) + \frac{f(\xi)}{T}.
	\end{array}
\end{equation}
}
	
	\noindent\textit{Verification:}
	 Let us find an upper bound for  $\nu_{D} \left( s^{\bot}_{m'} - s^{\bot}_{m}\right)$,
$m', m \in {\cal M}$. Since the operator $\nu_{D}$ defined by (\ref{CompensatedIntegral}) 
is just a compensated integral with respect to a Poisson process with mean measure 
$\mu(dtdx) =dt \eta(dx)$, we can apply 
Proposition \ref{ConctrIneq1} to 
		obtain that, for any $x'_{m'}>0$, and with probability larger than $1-e^{-x'_{m'}}$  
		\begin{equation}\label{Eq1}
			\nu_{D} \left( s^{\bot}_{m'} -s^{\bot}_{m}\right)
			\leq \Bigl\|\frac{s^{\bot}_{m'} -s^{\bot}_{m}}{T} \Bigr\|_{L^{2}(\mu)}\sqrt{2 x'_{m'}} +
			\frac{ \| s^{\bot}_{m'} -s^{\bot}_{m} \|_{\infty} x'_{m'}}{3T}.
		\end{equation}
	In that case, the probability that (\ref{Eq1}) holds for
	every $m' \in {\cal M}$ is larger than $1-\sum_{m'\in {\cal M}}e^{-x_{m'}}$ because 
	$P(A \cap B)\geq 1-a-b$, whenever $P(A)\geq 1-a$ and $P(B)\geq1-b$. Clearly,
	\begin{align}
		\Bigl\| \frac{s^{\bot}_{m'} -s^{\bot}_{m}}{T} \Bigl\|^{2}_{L^{2}(\mu)} 
		&=\iint\limits_{[0,T]\times D} \left(\frac{s^{\bot}_{m'}(x) -
		s^{\bot}_{m}(x)}{T}\right)^{2}s(x)dt\eta(dx)
		\nonumber \\
		 &\leq \| s\|_{\infty} \frac{\|s^{\bot}_{m'} -s^{\bot}_{m} \|^{2}}{T}.\nonumber
	\end{align}
Using (\ref{EasyIneq}-i), the first term on the right hand 
side of (\ref{Eq1}) is then bounded as follows:
\begin{equation}\label{Eq2} 
	\Bigl\| 
	\frac{s^{\bot}_{m'} 
	-s^{\bot}_{m}}{T} \Bigl\|_{L^{2}(\mu)}\sqrt{2 x'_{m'}} 
	\leq x_{1} \|s^{\bot}_{m'} 
	-s^{\bot}_{m} \|^{2} + \frac{\|s\|_{\infty}  x_{m'}}{2T x_{1}},
\end{equation}
for any $x_{1}>0$. Using (\ref{DefOfDm}) and (\ref{EasyIneq}-i),
\begin{align}
	\| s^{\bot}_{m'} -s^{\bot}_{m} \|_{\infty}x'_{m'} &
	\leq\left( \| s^{\bot}_{m'} \|_{\infty} + 
	\|s^{\bot}_{m} \|_{\infty} \right)x'_{m'} \nonumber \\
	&\leq \left( \sqrt{D_{m'}}\| s^{\bot}_{m'} \| 
	+  \sqrt{D_{m}}\|s^{\bot}_{m} \|\right)x_{m'}\nonumber \\
	&\leq \sqrt{D_{m'}}\| s \| x'_{m'}
	+  \sqrt{D_{m}}\|s \| x'_{m'} \nonumber \\
	&\leq 3 x_{2} D_{m'} + 3 x_{3} D_{m} +
	\frac{\| s \|^{2}{x'}_{m'}^{2}}{12} \left(\frac{1}{x_{2}}+
		  \frac{1}{x_{3}}\right), \nonumber
\end{align}
for all $x_{2}>0$, $x_{3}>0$. 
It follows that, for any $x_{1}>0$, $x_{2}>0$, and $x_{3}>0$,
\begin{align}
	\nu_{D} \left( s^{\bot}_{m'} -s^{\bot}_{m}\right) &\leq 
	x_{1} \|s^{\bot}_{m'} -s^{\bot}_{m} \|^{2} + 
	x_{2}\frac{D_{m'}}{T} + x_{3}\frac{D_{m}}{T} \nonumber \\
	&\;\;+ \frac{\|s\|_{\infty}  x'_{m'}}{2Tx_{1}} 
	+ \frac{\| s \|^{2}
	{x'}^{2}_{m'}}{36 T \bar{x}}, \nonumber
\end{align}
where we set $\frac{1}{\bar{x}} = \frac{1}{x_{2}}+\frac{1}{x_{3}}$. 
Next, take
\[
	x'_{m'} \equiv x_{4} \sqrt{d_{m'}}\left(\frac{1}{\|s\|}
	\wedge \frac{1}{\|s\|_{\infty}}\right)+\xi.
\]
Then, for any positive $x_{1}$, $x_{2}$, $x_{3}$, and $x_{4}$,  
there is a ${\cal K}$ and a function $f$ such that, 
with probability greater than $1 - {\cal K} e^{-\xi}$,
\begin{equation}\label{BoundForNu}
	\begin{array}{rl}
	\nu_{D} \left( s^{\bot}_{m'} 
	-s^{\bot}_{m}\right) &\leq x_{1} \|s^{\bot}_{m'} 
	-s^{\bot}_{m} \|^{2} + x_{2}\frac{D_{m'}}{T} 
	+ x_{3}\frac{D_{m}}{T}  \vspace{.2 cm}\\
	&\;\;+ \left(\frac{x_{4}^{2}}{18 \bar{x}} +  
	\frac{x_{4}}{2x_{1}}\right) \frac{ d_{m'}}{T} 
	+\frac{f(\xi)}{T},\:\:\:\forall m' \in {\cal M}. 
	\end{array}
\end{equation}
Concretely,  
\begin{equation}\label{FormPol0}
	\begin{array}{c}
	f(\xi) = \frac{\| s\|}{18 \bar{x}} \xi^{2}
	+\frac{\| s\|_{\infty}}{2 x_{1}} \xi,\vspace{.2 cm} \\
	{\cal K} = \Gamma \sum_{n=1}^{\infty} n^{R} 
	\exp\left(-\sqrt{n} x_{4} \left(\frac{1}{\|s\|}\wedge 
	\frac{1}{\|s\|_{\infty}}\right)\right).\vspace{.2 cm}\\
	\end{array}
\end{equation}
Here, we use the assumption of polynomial models 
(Definition\ref{Polynomial}) to come up with the constant 
${\cal K}$. Pluging (\ref{BoundForNu}) in (\ref{FundmIneq}), 
and renaming the coefficient of $d_{m'}/T$,
we can corroborate inequality 1.

\noindent \textit{Inequality 2: For any positive constants 
$y_{1}>1$, $x_{2}$,  $x_{3}$, and $x_{4}$, there are positive 
constants $C_{1}<1$, $C'_{1}>1$, and ${\cal K}$, 
and a strictly increasing quadratic polynomial $f$ 
(all independent of the class of linear models and $T$) 
such that with probability larger than $1-{\cal K} e^{-\xi}$, 
\begin{equation}\label{FundmIneq_v3}
	\begin{array}{ll}
	C_{1}\| s - \tilde s \|^{2} &\leq C'_{1}\| s 
	- s_{m}^{\bot} \|^{2} + y_{1} \chi_{\hat m}^{2}\\ 
	&+ x_{2}\frac{D_{\hat m}}{T} + x_{3}\frac{D_{m}}{T} 
	+ x_{4} \frac{d_{\hat m}}{T} \\
	&+\;{\rm pen}(m)-{\rm pen}(\hat m) + \frac{f(\xi)}{T}.
	\end{array}
\end{equation}
Moreover, if $1<y_{1}<2$, then $C'_{1}=3-y_{1}$ and 
$C_{1}=y_{1}-1$. If $y_{1} \geq 2$, then $C'_{1}=1+4x_{1}$ 
and $C_{1}=1-4x_{1}$, where $x_{1}$ is any positive 
constant related to $f$ according to equation (\ref{FormPol0}).\\
}\\
	\textit{Verification:} Let us combine the term on the left hand side of (\ref{FundmIneq_v2}) 
	with the first three terms on the right hand side. Using the triangle inequality 
	followed by (\ref{EasyIneq}-ii),
	\[
		 \|s^{\bot}_{\hat m}- s^{\bot}_{m} \|^{2}
		\leq 2\|s -s^{\bot}_{m} \|^{2}+
		2\|s^{\bot}_{\hat m} -s \|^{2}.
	\]  
	 Then, since  
	\(
		\chi_{\hat m}^{2} = \|s^{\bot}_{\hat m} -\hat s_{\hat m} \|^{2}
	\), 
	and
	\(
		 \| s^{\bot}_{\hat m} -s \|^{2}
		= \|s -\hat s_{\hat m} \|^{2}-\|s^{\bot}_{\hat m} -\hat s_{\hat m} \|^{2}
	\), it follows that  
	\[
		\nonumber
		\begin{array}{l}
		\| s - s_{m}^{\bot} \|^{2} +  2 \chi_{\hat m}^{2} +  
		2 x_{1} \|s^{\bot}_{\hat m} -s^{\bot}_{m} \|^{2}-\| s - \tilde s \|^{2}\\
		\leq \left(1 + 4 x_{1}\right) \| s - s_{m}^{\bot} \|^{2} 
		+  \left(2 - 4 x_{1}\right)\|s^{\bot}_{\hat m} -\hat s_{\hat m} \|^{2}\\
		\; + \left(4 x_{1} -1 \right)\| s - \tilde s \|^{2},
		\end{array}
	\]
	for every $x_{1}>0$. Then, for any $y_{1}>1$, 
	there are positive constant $C$, $C'_{1}> 1$, and $C_{1}< 1$ such that 
	\begin{equation}\label{Eq4}
		\begin{array}{l}
		\| s - s_{m}^{\bot} \|^{2} +  2 \chi_{\hat m}^{2} + 
		2 C \|s^{\bot}_{\hat m} -s^{\bot}_{m} \|^{2}-\| s - \tilde s \|^{2}\\
		\;\leq \; C'_{1} \| s - s_{m}^{\bot} \|^{2} +  y_{1}  \chi_{\hat m}^{2} -C_{1}\| s - \tilde s \|^{2}.
		\end{array}
	\end{equation}
	Combining (\ref{FundmIneq_v2}) and  (\ref{Eq4}), we obtain (\ref{FundmIneq_v3}).

\noindent \textit{Inequality 3: For any $y_{2}>1$ and positive 
constants $x_{i}$, $i=2,3,4$, there exist positive numbers 
$C_{1}<1$, $C'_{1}>1$, an increasing quadratic polynomial of 
the form $f_{2}(\xi)=a\xi^{2}+b\xi$, and a constant ${\cal K}_{2}>0$
(all independent of the family of linear models and of $T$) 
so that, with probability  greater than $1-{\cal K}_{2}e^{-\xi}$,
\begin{equation}\label{FundmIneq_v4}
	\begin{array}{rl} C_{1}\| s - 
	\tilde s \|^{2} &\leq C'_{1}\| s - s_{m}^{\bot} \|^{2} \\ 
	&\;\;+y_{2}\frac{V_{\hat m}}{T} + x_{2} \frac{D_{\hat m}}{T} 
	+ x_{3} \frac{d_{\hat m}}{T}-{\rm pen}(\hat m) \\
	&\;\;+x_{4} \frac{D_{m}}{T}+{\rm pen}(m) + \frac{f(\xi)}{T}.
	\end{array}
\end{equation}
}\\
\textit{Verification:} We bound $\chi^{2}_{m'}$ using 
Lemma \ref{IneqForChi} with ${\rm V} = \bbr_{+}\times D$ 
and $\mu(d {\bf x}) = s(x) dt \eta(dx)$.  We regard the linear model 
${\cal S}_{m}$ as  a subspace of 
$L^{2}(\bbr_{+} \times D,dt\eta(dx))$ with orthonormal basis 
\(
	\left\{ \frac{\varphi_{1,m}}{\sqrt{T}}, 
	\dots,\frac{\varphi_{d_{m},m}}{\sqrt{T}} \right\}.
\)
Recall that
\[
	\chi_{m}^{2} = \| s_{m}^{\bot} - \hat s_{m} \|^{2} = 
	\sum_{i=1}^{d} \left[\;\iint_{[0,T] \times D} \varphi_{i,m}(x) 
	 \frac{{\cal J}(dt,dx)-s(x)dt\eta(dx)}{T} \;\right]^{2}.
\]
Then, with probability larger than 
$1-\sum_{m'\in{\cal M}} e^{-x'_{m'}}$,
\begin{equation}\label{Eq5}
	\sqrt{T}\chi_{m'} \leq  (1+x_{1}) \sqrt{V_{m'}} 
	+ \sqrt{2 k {\rm M}_{m'} x'_{m'} } 
	+ k(x_{1}) {\rm B}_{m'} x'_{m'},
\end{equation}
for every $m' \in {\cal M}$, where $B_{m'}=\sqrt{D_{m'}/T}$, 
\begin{align}\label{DefForV}
	V_{m'} &\equiv \int_{D}  \left(\sum_{i=1}^{d_{m}} \varphi^{2}_{i,m}(x) \right) s(x)\eta(dx),
		\;\;{\rm and} \\
		M_{m'} &\equiv \sup\left\{ \int_{D} f^{2}(x)s(x)\eta(dx): f\in {\cal S}_{m'}, \|f\| =1 \right\}.\nonumber
	\end{align}
	Since $\int_{D} f^{2}(x)s(x)\eta(dx) \leq \|f\|_{\infty} \|s\|$, $M_{m'}$ is bounded above by 
	$\| s\| \sqrt{D_{m'}}$. In that case, we can use (\ref{EasyIneq}-i) to obtain 
	\[
		\sqrt{2 k {\rm M}_{m'} x'_{m'} } \leq x_{2} \sqrt{  D_{m'}}+ 
		\frac{k \| s \|}{2 x_{2}} x'_{m'},
	\]
	for any $x_{2}>0$. On the other hand, by hypothesis $D_{m'} \leq T$, and (\ref{Eq5})
	implies that 
	\[
		\sqrt{T}\chi_{m'} \leq (1+x_{1}) \sqrt{V_{m'}} + x_{2} \sqrt{  D_{m'}} + 
		 \left( \frac{k \| s \|}{2 x_{2}}+k(x_{1})\right) x'_{m'},
	\]
	where the constants $x'_{m'}$ are chosen as
	\[
		x'_{m'} = \frac{x_{3} \sqrt{  d_{m'}}}
		{ \frac{k \| s \|}{2 x_{2}}+k(x_{1})} + \xi.
	\]
	Then, for any $x_{1}>0$, $x_{2}>0$, $x_{3}>0$, and $\xi>0$, 
	\begin{equation}\label{IneqForChi1}
		\sqrt{T}\chi_{m'} \leq 
		(1+x_{1}) \sqrt{V_{m'}} + x_{2} \sqrt{  D_{m'}} +
		 x_{3} \sqrt{  d_{m'}} + f_{1}(\xi),
	\end{equation}
	with probability larger than $1-{\cal K}_{1} e^{-\xi}$, where
	\begin{equation}
		\label{FormPol}
		\begin{array}{c}
		f_{1}(\xi)= \left( \frac{k \| s \|}{2 x_{2}}+k(x_{1})\right) \xi,\vspace{.2 cm} \\
		{\cal K}_{1}= \Gamma \sum_{n=1}^{\infty} n^{R} 
		\exp\left(- \sqrt{n} x_{3}/\left( \frac{k \| s \|}{2 x_{2}}+k(x_{1})\right)\right).\vspace{.2 cm}\\
		\end{array}
	\end{equation}
	Squaring (\ref{IneqForChi1}) and using (\ref{EasyIneq}-ii) repeatedly, we conclude that, 
	for any $y>1$, $x_{2}>0$, and $x_{3}>0$, there are both
	a constant ${\cal K}_{1}>0$ and a quadratic function of the form 
	$f_{2}(\xi)= a\xi^{2}$
	(independent of $T$, $m'$, and the family of linear models) such that,  with probability 
	greater than $1-{\cal K}_{1} e^{-\xi}$,
	\begin{equation}\label{Eq6}
		\chi^{2}_{m'} \leq 
		y\frac{V_{m'}}{T} + x_{2} \frac{  D_{m'}}{T} +
		x_{3} \frac{  d_{m'}}{T} + \frac{f_{2}(\xi)}{T}, \;\;\;\forall m' \in {\cal M}. 
	\end{equation}
	 Then, (\ref{FundmIneq_v4}) immediately  follows from (\ref{Eq6}) and (\ref{FundmIneq_v3}).

\textit{Proof of  (\ref{PseudoOracleIneq}) for case (c):} \\
By the inequality (\ref{BoundForMean}), we can upper bound 
$V_{m'}$ by $\hat V_{m'}$ on an event of large probability. 
Namely, for every $x'_{m'}>0$ and $x>0$, with probability 
greater than $1-\sum_{m'\in {\cal M}}e^{-x'_{m'}}$
\begin{equation}\label{BoundForVHomog}
	(1+x)\left(\hat V_{m'}+
	\left( \frac{1}{2x}+\frac{5}{6}\right)
	\frac{D_{m'}}{T} x'_{m'}\right)\geq V_{m'},
	\;\;\forall m'\in{\cal M},
\end{equation}
(recall that 
$D_{m}=\|\sum_{i=1}^{d_{m}} \varphi^{2}_{i,m}\|_{\infty}$).  
Since by hypothesis $D_{m'}<T$, and choosing 
\[
	x'_{m'}= x' d_{m'}+\xi,\;\;(x'>0),
\]
it is seen that for any $x>0$ and $x_{4}>0$, there are
a positive constant ${\cal K}_{2}$ and a function $f(\xi)=b \xi$
(independent of $T$ and of the linear models) such that 
with probability greater than $1-{\cal K}_{2}e^{-\xi}$
\begin{equation}\label{BoundForV}
	(1+x)\hat V_{m'}+x_{4}d_{m'}+f(\xi) \geq
	V_{m'},\;\;\forall m'\in {\cal M}.
\end{equation}
Here, we get ${\cal K}_{2}$ from the Polynomial assumption 
on the class of models. Combining (\ref{BoundForV}) and 
(\ref{FundmIneq_v4}), it is clear that for any $y_{2}>1$, 
and positive $x_{i}$, $i=1,2,3$, we can 
choose a pair of positive constants $C_{1}<1$, $C'_{1}>1$, 
an increasing quadratic polynomial of the form 
$f(\xi)=a\xi^{2}+b\xi$, and a constant ${\cal K}>0$
(all independent of the family of linear models and of $T$) so that,  
with probability greater than $1-{\cal K}e^{-\xi}$ 
\begin{equation}\label{FundmIneq_v5}
	\begin{array}{rl}
	C_{1}\| s - \tilde s \|^{2} &\leq
	C'_{1}\| s - s_{m}^{\bot} \|^{2} \\ 
	&\;\;+y_{2}\frac{\hat V_{\hat m}}{T} 
	+ x_{1} \frac{D_{\hat m}}{T} 
	+ x_{2} \frac{d_{\hat m}}{T}-{\rm pen}(\hat m) \\
	&\;\;+x_{3} \frac{D_{m}}{T}+{\rm pen}(m) + \frac{f(\xi)}{T}.
	\end{array}
\end{equation}
Next, we take 
$y_{2}=c$, $x_{1}=c'$, and $x_{2}=c''$ to cancel
$-pen(\hat m)$ in (\ref{FundmIneq_v5}). By Lemma \ref{AuxLemma},
it follows that 
\begin{equation}\label{OracleIneqForm}
	C_{1}\bbe\left[\| s - \tilde s \|^{2}\right] \leq
	C'_{1}\| s - s_{m}^{\bot} \|^{2}+
	\left( 1+\frac{x_{3}}{c'}\right)\bbe\left[ {\rm pen}(m)\right] 
	+\frac{C''_{1}}{T}.
\end{equation}
Since $m$ is arbitrary, we obtain the case (c) of  
(\ref{PseudoOracleIneq}).

	\textit{Proof of  (\ref{PseudoOracleIneq}) for case (a):} \\
	By Remark \ref{DmInOB}, we can bound $V_{m'}$, as given in (\ref{DefForV}), 
	by $D_{m'}\rho$ (assuming that $\rho<\infty$). 
	On the other hand, (\ref{BoundForMean}) implies that 
	\begin{equation}\label{BoundForRho}
		(1+x_{1})\left(\frac{{\cal N}}{T}+
		  \left( \frac{1}{2x_{1}}+\frac{5}{6}\right)\frac{\xi}{T}\right)\geq \rho,
	\end{equation}
	with probability greater than $1-e^{-\xi}$.  Using these bounds for $V_{m'}$ and the assumption 
	that $D_{m'} \leq T$, (\ref{FundmIneq_v4}) reduces to
	\begin{equation}\label{FundmIneq_v6}
		\begin{array}{ll}
		C_{1}\| s - \tilde s \|^{2} &\leq C'_{1}\| s - s_{m}^{\bot} \|^{2} \\ 
		 &+y \frac{D_{\hat m} {\cal N}}{T^{2}} 
		  + x_{1} \frac{d_{\hat m}}{T}-{\rm pen}(\hat m) \\
		 &+x_{2} \frac{D_{m}{\cal N}}{T^{2}}+{\rm pen}(m) + \frac{f(\xi)}{T}, 
		\end{array}
	\end{equation}
	which is valid with probability $1-{\cal K}e^{-\xi}$. In (\ref{FundmIneq_v6}),  
	$y>1$, $x_{1}>0$ and $x_{2}>0$ are arbitrary, while  
	$C_{1}$, $C'_{1}$, the increasing quadratic polynomial of the form $f(\xi)=a\xi^{2}+b\xi$, and
	a constant ${\cal K}>0$ are determined by $y$, $x_{1}$, and $x_{2}$ independently 
	of the family of linear models and of $T$. 
	We point out that we divided and multiplied by $\rho$ the terms
	$D_{\hat m}/T$ and $D_{m}/T$ in (\ref{FundmIneq_v4}), and then applied
	(\ref{BoundForRho}) to get (\ref{FundmIneq_v6}). It is now clear that 
	$y=c$, and $x_{1}=c'$ will produce  the desired cancelation. 

	\textit{Proof of  (\ref{PseudoOracleIneq}) for case (b):} \\
	We first upper bound $D_{\hat m}$
	by $\beta^{-1} V_{\hat m}$ and $d_{\hat m}$ by $(\beta
	\phi)^{-1} V_{\hat m}$ in the inequality (\ref{FundmIneq_v4}):
	\begin{equation}\label{FundmIneq_v7}
		\begin{array}{rl}
		C_{1}\| s - \tilde s \|^{2} &\leq
		C'_{1}\| s - s_{m}^{\bot} \|^{2}
		+\left(y+x_{1}\beta^{-1}+x_{2}(\beta\phi)^{-1} \right)
		 \frac{V_{\hat m}}{T} \\
		 &\;\;-{\rm pen}(\hat m)+x_{3} \beta^{-1} \frac{V_{m}}{T}+{\rm pen}(m) + \frac{f(\xi)}{T}.
		\end{array}
	\end{equation}
	Then, using $d_{m'}\leq (\beta\phi)^{-1} V_{m'}$ in (\ref{BoundForV}) and
	letting $x_{4}(\beta\phi)^{-1}$ vary on $(0,1)$, we verify that for any
	$x'>0$, a positive constant ${\cal K}_{4}$ and a polynomial $f$ can be
	found so that with probability greater than $1-{\cal K}_{4} e^{-\xi}$,
	\begin{equation}\label{BoundForV_2}
		(1+x')\hat V_{m'}+f(\xi) \geq
		V_{m'},\;\;\forall m'\in {\cal M}.
	\end{equation}
	 Putting together 
	(\ref{BoundForV_2}) and (\ref{FundmIneq_v7}), it is clear that for any
	$y>1$ and $x_{1}>0$, we can 
	find a pair of positive constants $C_{1}<1$, $C'_{1}>1$, an increasing quadratic
	polynomial of the form $f(\xi)=a\xi^{2}+b\xi$, and a constant ${\cal K}>0$
	(all independent of the family of linear models and of $T$) so that,
	with probability greater than $1-{\cal K}e^{-\xi}$,  
	\begin{equation}\label{FundmIneq_v8}
		\begin{array}{rl}
		C_{1}\| s - \tilde s \|^{2} &\leq
		C'_{1}\| s - s_{m}^{\bot} \|^{2} + y\frac{\hat V_{\hat m}}{T}-{\rm pen}(\hat m) \\
		 &\;\;+x_{1} \frac{V_{m}}{T}+{\rm pen}(m) + \frac{f(\xi)}{T}.
		\end{array}
	\end{equation}
	In particular, by taking $y=c$, the term $-pen(\hat m)$ cancels out. Lemma \ref{AuxLemma}
	implies that  
	\begin{equation}\label{OracleIneqForm_2}
		C_{1}\bbe\left[\| s - \tilde s \|^{2}\right] \leq
		C'_{1}\| s - s_{m}^{\bot} \|^{2}+\left( 1+x_{1}\right)\bbe\left[ {\rm pen}(m)\right] +
		\frac{C''_{1}}{T}.
	\end{equation}
	Finally, (\ref{PseudoOracleIneq}) (b) follows since $m$ is arbitrary.\hfill \(\Box \)

\noindent\begin{remk} \label{Constants}
Let us analyze more carefully the values that the constants 
$C$ and $C'$ can take in the inequality (\ref{PseudoOracleIneq}). 
For instance, consider the penalty function of part (c). 
As we saw in (\ref{OracleIneqForm}), 
the constants $C$ and $C'$ are determined by $C_{1}$, $C'_{1}$, 
$C''_{1}$, and $x_{3}$. 
The constant $C_{1}$ was proved to be $y_{1}-1$ if $1<y_{1}<2$, 
while it can be made arbitrarily close to one otherwise 
(see the comment immediately after (\ref{FundmIneq_v3})). 
On the other hand, $y_{1}$ itself can be made arbitrarily close 
to the penalization parameter $c$ since 
$c=y_{2}=y_{1}(1+x)y$, 
where $x$ is as in (\ref{BoundForVHomog}) and  
$y$ is in (\ref{Eq6}). 
Then, when $c \geq 2$, 
$C_{1}$ can be made arbitrarily close to one at the cost of 
increasing $C''_{1}$ in (\ref{OracleIneqForm}). 
Similarly, paying the same cost, we are able to
select $C'_{1}$ as close to one as we wish and $x_{3}$ 
arbitrarily small.
Therefore, it is possible to find for any $\varepsilon>0$, a constant
$C'(\varepsilon)$ (increasing in $\varepsilon$) so that 
\begin{equation}\label{PseudoOracleIneq_Better}
	\bbe \| s - \tilde s \|^{2} \leq 
	\left( 1+\varepsilon\right) \inf_{m\in {\cal M}} 
	\left\{ \| s - s^{\bot}_{m}\|^{2} +
	\bbe \left[ {\rm pen}(m)\right]\right\}
	+\frac{C'(\varepsilon)}{T}.
\end{equation}
A more thorough inspection shows that 
\[
	\lim_{\varepsilon \rightarrow 0} C'(\varepsilon) 
	\varepsilon = K,
\]
where $K$ depends only $c$, $c'$, $c''$, $\Gamma$, $R$, $\| s\|$, 
and $\|s \|_{\infty}$. The same reasoning apply to the other two
types of penalty functions when $c \geq 2$. 
In particular, we point out that
$\tilde{C}$ can be made arbitrarily close to 2 in the Oracle inequality
(\ref{OracleIneq}) at the price of having a large $\tilde{C}'$ 
constant.  
\end{remk}


\subsection{Proof of the minimax results}\label{SectionProofMinimax}
\noindent{\bf Proof of Theorem \ref{MinimaxResult1}:}\\
Fix a L\'evy density $s_{0}\in\Theta_{\alpha}\left(L/2;[a,b]\right)$
such that $s_{0}(x)>0$, for all $x\in [a,b]$, 
and a constant $\kappa>0$. 
Consider a bounded function $g:\bbr\rightarrow\bbr_{+}$, 
with compact support $\bbk\subset [-1,1]$, 
that meets (\ref{HolderCond}) with $L/2$ (instead of $L$) 
for all $x_{1}, x_{2}\in\bbr$, $g(0)>0$, increasing for $x<0$, 
and decreasing for $x>0$. 
Moreover, the support and the maximum value of $g$ are chosen 
small enough so that
\[
	s_{0}(x) - 
	\kappa^{-\alpha}g\left(\kappa (x-x_{0})\right)>0,
	\;\; \forall x\in[a,b],
\]
and the support of $g\left(\kappa (x-x_{0})\right)$ is contained 
in $(a,b)$.
Let us consider the parametric model
\[
	s_{\theta}(x):= s_{0}(x) +\theta T^{-\frac{\alpha}{2\alpha+1}} 
	g\left(\kappa T^{\frac{1}{2\alpha+1}}(x-x_{0})\right),\;\; x\in\bbr_{0},
\]
parametrized by $\theta\in(-\kappa^{-\alpha},\kappa^{-\alpha})$.
Notice that the function $s_{\theta}$ is a L\'evy density
for any $T>1$ and $|\theta|<\kappa^{-\alpha}$. 
Now, for $x_{1},x_{2}\in [a,b]$,
\begin{align*}
	| s^{(k)}_{\theta}(x_{1})-s^{(k)}_{\theta}(x_{2})|&\leq
	\left| s^{(k)}_{0}(x_{1})-s^{(k)}_{0}(x_{2})\right|+\\
	|\theta | \kappa^{k} T^{\frac{-\alpha+k}{2\alpha+1}} &
		\left| g^{(k)}\left(\kappa T^{\frac{1}{2\alpha+1}}(x_{1}-x_{0})\right)
		-g^{(k)}\left(\kappa T^{\frac{1}{2\alpha+1}}(x_{2}-x_{0})\right)\right|\\  
		&\leq \frac{L}{2}|x_{1}-x_{2}|^{\beta}+
		\frac{L}{2} |\theta | \kappa^{k+\beta}
		T^{\frac{-\alpha+k+\beta}{2\alpha+1}} |x_{1}-x_{2}|^{\beta}\\
		&\leq L |x_{1}-x_{2}|^{\beta},
	\end{align*}
	implying that $s_{\theta}\in \Theta$ whenever 
$|\theta|<\kappa^{-\alpha}$.

Let $\calM_{0}$ be the space of atomic measures on 
$[0,T]\times[a,b]$ and 
let $\bbp^{(T)}_{\theta}$ be the probability measure on 
$\calM_{0}$ induced by 
those jumps of the L\'evy process
$\left\{ X(t) \right\}_{0\leq{t}\leq{T}}$ 
whose sizes lie on $[a,b]$ and where the L\'evy density of 
the process is $s_{_{\theta}}$. In other words, 
$\bbp^{(T)}_{\theta}$ is the distribution of a Poisson 
process on $[0,T]\times[a,b]$ with mean measure 
$dt s_{_{\theta}}(x)dx$. Using Theorem 1.3 of \cite{Kutoyants},
\begin{align*}
	\frac{d\bbp^{(T)}_{\theta}}{d\bbp^{(T)}_{0}}(\xi)&=
	\exp\left\{\int_{0}^{T}\!\!\!
	\int_{a}^{b}\ln\left\{1+\theta 
	T^{-\frac{\alpha}{2\alpha+1}} 
	s^{-1}_{0}(x)g\left(T^{\frac{1}{2\alpha+1}}
	(x-x_{0})\right)\right\} 
	\xi(dt,dx)\right.\\
	&\;\;\;\;\;\;\;\;\;\;-\left.\theta 
	T^{1-\frac{\alpha}{2\alpha+1}} 
	\int_{a}^{b}g\left(T^{\frac{1}{2\alpha+1}}
	(x-x_{0})\right)dx\right\}.
\end{align*}
The goal is to prove the LAN (\emph{local asymptotic normality}) 
property for the parametric model 
$\left\{\bbp^{(T)}_{\theta}:\theta\in
\left(-\kappa^{-\alpha},\kappa^{-\alpha}\right)\right\}$ 
at $\theta=0$ (see definition 2.1 of \cite{Kutoyants}).
Now, define $R(u)=\ln(1+u)-u+\frac{u^{2}}{2}$. 
The right hand side of the 
above equation can be written as  follows:
\[
	\frac{d\bbp^{(T)}_{\theta}}{d\bbp^{(T)}_{0}}(\xi)=
	\exp\left\{ \theta\Delta_{_{T}}-\frac{\theta^{2}}{2}
	\sigma_{_{T}}^{2}+r_{_{T}}(\theta)\right\},
\]
where
\begin{align*}
	&\Delta_{_{T}}=T^{-\frac{\alpha}{2\alpha+1}}
	\int_{0}^{T}\!\!\!\int_{a}^{b}
	s_{0}^{-1}(x) g\left(T^{\frac{1}{2\alpha+1}}
	(x-x_{0})\right)
	\left[\xi(dt,dx)-s_{0}(x)dtdx\right],\\
	&\sigma_{_{T}}^{2}=T^{1-\frac{2\alpha}{2\alpha+1}}
	\int_{a}^{b} s_{0}^{-1}(x)g^{2}
	\left(T^{\frac{1}{2\alpha+1}}(x-x_{0})\right)dx,\\
	&r_{_{T}}(\theta)=-\frac{\theta^{2}}{2}
	T^{-\frac{2\alpha}{2\alpha+1}}
	\int_{0}^{T}\!\!\!\int_{a}^{b}s_{0}^{-2}(x)
	g^{2}\left( T^{\frac{1}{2\alpha+1}}(x-x_{0})\right)
	\left[\xi(dt,dx)-s_{0}(x)dtdx\right]\\
	&\;\;\;\;\;\;\;\;\;+\int_{0}^{T}\!\!\!\int_{a}^{b} 
	R\left(\theta T^{-\frac{\alpha}{2\alpha+1}} 
	s_{0}^{-1}(x)g\left(T^{\frac{1}{2\alpha+1}}
	(x-x_{0})\right)\right)\xi(dt,dx).	
\end{align*}
We want to prove that there are nomalizing 
constants $\varphi_{T}>0$ such that
\[ 
	\calL_{\bbp^{(T)}_{0}}
	\left(\varphi_{T}\Delta_{_{T}}\right)
	\stackrel {\frak {D}}{\rightarrow}\calN(0,1),\;\;
	\varphi_{T}^{2}\sigma^{2}_{_{T}}\rightarrow 1,\;\; 
	{\rm and}\;\; 
	r_{_{T}}(\theta)\stackrel {\bbp^{(T)}_{0}}{\rightarrow}0
\]
as $T\rightarrow\infty$.
To prove the first limit, we invoke the CLT for Poisson 
integrals by verifying the Liapunov condition 
(see Theorem 1.1 and Remark 1.2 of \cite{Kutoyants}).
For $T>1$, we have that 
		\begin{align*}
			&T^{-\frac{\alpha(2+\delta)}{2\alpha+1}}
			\int_{0}^{T}\int_{a}^{b}s_{0}^{-2-\delta}(x)
			  g^{2+\delta}\left(\kappa T^{\frac{1}{2\alpha+1}}(x-x_{0})\right)\left(s_{0}(x)\right)dxdt=\\
			 &\;\;\;\kappa^{-1}T^{1-\frac{\alpha(2+\delta)}{2\alpha+1}-\frac{1}{2\alpha+1}}
			\int_{\bbk}s_{0}^{-1-\delta}(\kappa^{-1}T^{-\frac{1}{2\alpha+1}}u+x_{0})
			  g^{2+\delta}\left(u\right)du\stackrel{T\rightarrow\infty}{\longrightarrow}0.
		\end{align*}
Notice also that, for large enough $T$, 
\begin{align*}
	{\rm Var}\left(\Delta_{_{T}}\right)
	&=T^{-\frac{2\alpha}{2\alpha+1}}
	\int_{0}^{T}\int_{a}^{b}s_{0}^{-2}(x)
	g^{2}\left(\kappa 
	T^{\frac{1}{2\alpha+1}}(x-x_{0})\right)
	\left(s_{0}(x)\right)dxdt=\\
	&=\kappa^{-1}T^{1-\frac{2\alpha}{2\alpha+1}
	-\frac{1}{2\alpha+1}}
	\int_{\bbk}s_{0}^{-1}
	(\kappa^{-1}T^{-\frac{1}{2\alpha+1}}u+x_{0})
	g^{2}\left(u\right)du\\
	&\stackrel{T\rightarrow\infty}{\longrightarrow}
	\kappa^{-1}s_{0}^{-1}(x_{0})
	\int_{\bbk} g^{2}(u) du.
\end{align*}
Then,  $\calL_{\bbp^{(T)}_{0}}\left(\Delta_{_{T}}\right)
\stackrel {\frak {D}}{\rightarrow}\calN(0,I_{0}^{2})$
with
$I_{0}^{2}:= \kappa^{-1}s_{0}^{-1}(x_{0})\int_{\bbk} g^{2}(u) du$. 
In particular, we also prove that 
$\sigma_{_{T}}^{2}\rightarrow I_{0}^{2}$. 
We now verify that $r_{_{T}}(\theta)$ vanishes in probability. 
Notice that the first term of $r_{_{T}}$ has mean $0$
and variance
\begin{align*}
	&\frac{\theta^{4}}{4}T^{-\frac{4\alpha}{2\alpha+1}}
	\int_{0}^{T}\int_{a}^{b}
	s_{0}^{-4}(x)g^{4}\left(\kappa T^{\frac{1}{2\alpha+1}}
	(x-x_{0})\right)(s_{0}(x))dxdt\\
	&=\frac{\theta^{4}}{4}\kappa^{-1}
	T^{1-\frac{4\alpha}{2\alpha+1}-\frac{1}{2\alpha+1}}
	\int_{\bbk}s_{0}^{-4}(\kappa^{-1}T^{-\frac{1}{2\alpha+1}}u
	+x_{0})g^{4}\left(u\right)du\stackrel
	{T\rightarrow\infty}{\longrightarrow}0.
\end{align*}
Then, the first term converges in probability to $0$. 
Similarly, the second term of $r_{T}(\theta)$ converges 
to $0$ in probability because its mean and variance both
goes to $0$. Indeed, using that $| R(u)| \leq |u|^{3}/3$, 
the absolute value of its expectation satisfies
	\begin{align*}
		&\left|\int_{0}^{T}\int_{a}^{b} R\left(\theta T^{-\frac{\alpha}{2\alpha+1}} s^{-1}_{0}(x)
		g\left(\kappa T^{\frac{1}{2\alpha+1}}(x-x_{0})\right)\right)(s_{0}(x))dx dt\right|\\
		&\leq\frac{|\theta|^{3}}{3} T^{1-\frac{3\alpha}{2\alpha+1}}\int_{a}^{b}
		  s_{0}^{-2}(x)g^{3}\left(\kappa T^{\frac{1}{2\alpha+1}}(x-x_{0})\right)dx
		  \stackrel{T\rightarrow\infty}{\longrightarrow}0.
	\end{align*} 
A similar reasoning applies to the variance. Therefore, 
$\{ \bbp^{(T)}_{\theta}\}_{\theta\in(-\kappa^{-\alpha},
\kappa^{-\alpha})}$ is Locally Asymptotically Normal (LAN)
at $\theta=0$ 
(with the normalizing constants $\varphi_{T}:= I_{0}^{-1}$).
 We are now in position of using the theory
for LAN families (see \cite{Ibragimov} for the general theory 
and \cite{Kutoyants} for the case of Poisson processes). 
In particular, by (2.11) of \cite{Kutoyants},
if for each $T>0$, $\hat{\theta}_{_{T}}$ is an arbitrary  
estimator of $\theta$, 
based on the jumps of the L\'evy process happening on 
or before $T$ and with sizes in $[a,b]$, then
\begin{equation}\label{ParmIneq}
	\liminf_{T\rightarrow\infty} \sup_{|\theta |
	<\kappa^{-\alpha}} 
	\bbe_{\theta} \left[\ell_{0}\left(I_{0}
	\left(\hat{\theta}_{_{T}}-\theta\right)\right)\right]\geq
	B,
\end{equation}
where $B:=\bbe\left[ \ell_{0}(Z)\chi_{\left[|Z|
<I_{0}\kappa^{-\alpha}/2\right]}\right]$
and $Z\sim \mathcal{N}(0,1)$.

Now, for each  $T>0$, let $\hat{s}_{_{T}}$ be an arbitrary estimator, 
based on the jumps of the L\'evy process happening on or 
before $T$ and with sizes in $[a,b]$. 
Clearly, $\hat{s}_{_{T}}$ induces the estimator
	\(
		\hat{\theta}_{_{T}}:= T^{\frac{\alpha}{2\alpha+1}}g^{-1}(0)
		\left(\hat{s}_{_{T}}(x_{0})-s_{0}(x_{0})\right),
	\)
	and since 
	\(
		\theta=T^{\frac{\alpha}{2\alpha+1}}g^{-1}(0)\left(s_{\theta}(x_{0})-s_{0}(x_{0}) \right),
	\)
	we can write
	\[
		g(0)\left(\hat{\theta}_{_{T}}-\theta\right)=T^{\frac{\alpha}{2\alpha+1}}
		\left(\hat{s}_{_{T}}(x_{0})-s_{\theta}(x_{0})\right).
	\]
	If we take $\ell_{0}(u):= \ell\left(g(0)I_{0}^{-1}u\right)$, (\ref{ParmIneq}) becomes:
	\begin{align*}
		B&\leq 
		\liminf_{T\rightarrow\infty} \sup_{|\theta |<\kappa^{-\alpha}} 
		\bbe_{\theta} \left[\ell_{0}\left(I_{0}\left(\hat{\theta}_{_{T}}-\theta\right)\right)\right]\\
		&=\liminf_{T\rightarrow\infty} \sup_{|\theta |<\kappa^{-\alpha}} 
		\bbe_{\theta} \left[\ell\left(
		T^{\frac{\alpha}{2\alpha+1}}
		\left(\hat{s}_{_{T}}(x_{0})-s_{\theta}(x_{0})\right)\right)\right]
	\end{align*}
Since $\{s_{\theta}:\theta\in(-k^{-\alpha},k^{-\alpha})\}\subset\Theta$, 
\begin{equation}\label{MinMaxIneq1c}
	\liminf_{T\rightarrow\infty} \sup_{s\in\Theta}
	\bbe_{s} \left[\ell\left(T^{\frac{\alpha}{2\alpha+1}}
	\left(\hat{s}_{_{T}}(x_{0})-s(x_{0})\right)\right)\right]\geq B,
\end{equation}
where 
\begin{equation}\label{LowerBound}
	B:=2^{-3/2}\pi^{-1/2}\int_{|z|<I_{0}\kappa^{-\alpha}/2}
	\ell(g(0)I_{0}^{-1}z)e^{-z^{2}/2}dz.
\end{equation}
This implies (\ref{MinMaxIneq1b}) because the lower bound $B$ 
does not depend on the family of estimators $\hat{s}_{_{T}}$.
 Indeed, for each $\varepsilon>0$, let $\hat{s}^{(\varepsilon)}_{_{T}}$ 
be such that 
\begin{align*}
	&\sup_{s\in\Theta}\bbe_{s} \left[\ell
	\left(T^{\frac{\alpha}{2\alpha+1}}
	\left(\hat{s}^{(\varepsilon)}_{T}(x_{0})-s(x_{0})\right)\right)
	\right]\\
	&<\inf_{\hat{s}_{_{T}}}\sup_{s\in\Theta}
	\bbe_{s} \left[\ell\left(T^{\frac{\alpha}{2\alpha+1}}
	\left(\hat{s}_{_{T}}(x_{0})-s(x_{0})\right)\right)\right]+
	\varepsilon.
\end{align*}
Taking the $\liminf$ as $T\rightarrow\infty$ on both sides, 
we obtain (\ref{MinMaxIneq1b}) since $\varepsilon$ is arbitrary.
\hfill$\Box$

\noindent{\bf Proof of Corollary \ref{MinimaxResult2}:}\\
We first notice that the proof of Theorem \ref{MinimaxResult1}
can be modified so that  (\ref{MinMaxIneq1c}) holds true even if 
$x_{0}$ is not fixed. That is, 
for any family of estimators $\hat{s}_{_{T}}$ and points 
$x_{_{T}}\in(a,b)$,  
\begin{equation}\label{MinMaxIneq2c}
	\liminf_{T\rightarrow\infty} \sup_{s\in\Theta}
	\bbe_{s} \left[\ell\left(T^{\frac{\alpha}{2\alpha+1}}
	\left(\hat{s}_{_{T}}(x_{_{T}})-s(x_{_{T}})\right)\right)
	\right]\geq C,
\end{equation}
for a constant $C>0$, which is independent of the family of 
estimators and of the points. 
Indeed, we can construct a parametric model of the form 
\[
	s_{_{\theta,T}}(x):= s_{0}(x) +
	\theta T^{-\frac{\alpha}{2\alpha+1}} 
	g_{_{T}}\left(\kappa T^{\frac{1}{2\alpha+1}}
	(x-x_{_{T}})\right),\;\; x\in\bbr_{0},
\]
where $|\theta|<\kappa^{-\alpha}$ and where 
$g_{_{T}}$ is as in the previous 
proof. Moreover,  without loss of generality, 
$0<\inf_{T} g_{_{T}}(0)\leq \sup_{T} g_{_{T}}(0) <\infty$, since 
$s_{0}$ is continuous and strictly positive on $(a,b)$.
Let $\bbp_{\theta}^{(T)}$ be the distribution of a Poisson process on 
$[0,T]\times [a,b]$ with mean measure $dt s_{_{\theta,T}}(x)dx$.
Following the same arguments as above, 
$\{ \bbp^{(T)}_{\theta}:\theta\in(-\kappa^{-\alpha},
\kappa^{-\alpha})\}$ is Locally Asymptotically Normal (LAN)
at $\theta=0$ with the normalizing constants 
\[
	\varphi_{_{T}}:= 
	\kappa^{2}
	\left(\int_{\bbk_{T}}s_{0}^{-1}
	(\kappa^{-1}T^{-\frac{1}{2\alpha+1}}u+x_{_{T}})
	g_{_{T}}^{2}\left(u\right)du\right)^{-2}
\]
Observe that there is an $m>0$ for which 
$\inf_{T} \varphi_{T}\geq m$. By (2.11) of \cite{Kutoyants},
for any $\delta>0$,
\begin{equation}
	\liminf_{T\rightarrow\infty} \sup_{|\theta |
	<\delta\varphi_{_{T}}} 
	\bbe_{\theta} \left[\ell_{0}\left(\varphi_{_{T}}^{-1}
	\left(\hat{\theta}_{_{T}}-\theta\right)\right)\right]\geq
	C,
\end{equation}
where $C:=\bbe\left[ \ell_{0}(Z)\chi_{\left[|Z|
<\delta/2\right]}\right]$ and $Z\sim \mathcal{N}(0,1)$.
Since $\varphi_{_{T}}\geq m$ and $\ell_{0}(|y|)$ is 
increasing in $y$, 
\begin{equation}
	\liminf_{T\rightarrow\infty} \sup_{|\theta |
	<\delta\varphi_{_{T}}} 
	\bbe_{\theta} \left[\ell_{0}\left(m^{-1}
	\left(\hat{\theta}_{_{T}}-\theta\right)\right)\right]\geq
	C,
\end{equation}
Now, take  
\[
	\hat{\theta}_{_{T}}:= T^{\frac{\alpha}{2\alpha+1}}
	g_{_{T}}^{-1}(0)
	\left(\hat{s}_{_{T}}(x_{_{T}})-s_{0}(x_{_{T}})\right).
\]

Since 
\(
	\theta=T^{\frac{\alpha}{2\alpha+1}}g_{_{T}}^{-1}(0)
	\left(s_{_{\theta,T}}(x_{_{T}})-s_{0}(x_{_{T}}) \right)
\),
\begin{align*}
	\sup_{s\in\Theta}
	\bbe_{s} \left[\ell\left(T^{\frac{\alpha}{2\alpha+1}}
	\left(\hat{s}_{_{T}}(x_{_{T}})-s(x_{_{T}})\right)\right)
	\right]
	&\geq 
	\sup_{|\theta |<\delta\varphi_{_{T}}} 
	\bbe_{\theta} \left[\ell
	\left(g_{_{T}}(0)\left(\hat{\theta}_{_{T}}
	-\theta\right)\right)\right]\\
	&\geq \sup_{|\theta |<\delta\varphi_{_{T}}} 
	\bbe_{\theta} \left[\ell
	\left(\tilde{m}\left(\hat{\theta}_{_{T}}
	-\theta\right)\right)\right],
\end{align*}
where $\tilde{m}=\inf_{_{T}} g_{_{T}}(0)$. 
Taking $\liminf$ as $T\rightarrow\infty$, 
(\ref{MinMaxIneq2c}) is obtained with 
\begin{equation}\label{LowerBound2}
	C=2^{-3/2}\pi^{-1/2}\int_{|z|<\delta/2}
	\ell(\tilde{m}\,m\, z)e^{-z^{2}/2}dz.
\end{equation}
Finally, (\ref{MinMaxIneq2a}) can be deduced as follows.
For each $\varepsilon>0$, 
 let $\hat{s}^{(\varepsilon)}_{_{T}}\in\Theta$ and 
$x_{_{T}}^{(\varepsilon)}\in(a,b)$ be such that 
\begin{align*}
&\sup_{s\in\Theta}\bbe_{s} 
\left[\ell\left(T^{\frac{\alpha}{2\alpha+1}}
\left(\hat{s}^{(\varepsilon)}_{T}
\left(x_{_{T}}^{(\varepsilon)}\right)-
s\left(x_{_{T}}^{(\varepsilon)}\right)\right)\right)\right]\leq\\
&\inf_{x\in(a,b)}\inf_{\hat{s}_{_{T}}}\sup_{s\in\Theta}
\bbe_{s} \left[\ell\left(T^{\frac{\alpha}{2\alpha+1}}
\left(\hat{s}_{_{T}}(x)-s(x)\right)\right)\right]+\varepsilon.
\end{align*}
Next, take the $\liminf$ as $T\rightarrow\infty$ on both sides
above, and  apply (\ref{MinMaxIneq2c}). Finally,  let $\varepsilon\rightarrow{0}$.
\hfill$\Box$

\noindent{\bf Proof of Corollary \ref{MinimaxResult3}:}\\
	Fix a measurable estimator $\hat{s}_{_{T}}$  and a $s\in\Theta$. By Fubini's Theorem, 
	\[
		\bbe_{s} \left[\int_{a}^{b}\left(\hat{s}_{_{T}}(x)-s(x)\right)^{2}dx\right]
		=\int_{a}^{b}\bbe_{s} \left[\left(\hat{s}_{_{T}}(x)-s(x)\right)^{2}\right] dx.
	\]
	Now, for each $\varepsilon>0$, there exists an $x_{0}^{(\varepsilon)}\in(a,b)$ satisfying 
	\[
		\frac{1}{b-a}\int_{a}^{b}\bbe_{s} \left[\left(\hat{s}_{_{T}}(x)-s(x)\right)^{2}\right] dx \geq
		\bbe_{s} \left[\left(\hat{s}_{_{T}}\left(x_{0}^{(\varepsilon)}\right)-
		s\left(x_{0}^{(\varepsilon)}\right)\right)^{2}\right] -\varepsilon.
	\]
Then,
\begin{align*}
	\frac{1}{b-a}\sup_{s\in\Theta}
	\bbe_{s} \left[\int_{a}^{b}\left(\hat{s}_{_{T}}(x)-s(x)\right)^{2}dx\right] &\geq
	\sup_{s\in\Theta}\bbe_{s} \left[\left(\hat{s}_{_{T}}(x_{0}^{(\varepsilon)})-
	s(x_{0}^{(\varepsilon)})\right)^{2}\right] -\varepsilon\\
	&\geq \inf_{x\in(a,b)}
	\sup_{s\in\Theta}\bbe_{s} \left[\left(\hat{s}_{_{T}}\left(x\right)-
		s\left(x\right)\right)^{2}\right] -\varepsilon.
	\end{align*}
	Letting $\varepsilon\rightarrow{0}$, (\ref{MinMaxIneq3a}) becomes a consequence of 
	(\ref{MinMaxIneq2a}) with $\ell(u)=u^{2}$.
	 \hfill$\Box$


\subsection{Some additional proofs}\label{proofsChap3}

\noindent{\bf Proof of Corollary \ref{ConvergenceBesov}:} 
The idea is to estimate the bias and the penalized term 
in (\ref{PseudoOracleIneq}). Clearly, the dimension 
$d_{m}$ of $\calS_{m}^{k}$ is $m(k+1)$. 
Also, $D_{m}$ is bounded by $(k+1)^{2} m /(b-a)$ 
(see (7) in \cite{Birge}), and 
\[
	\bbe\left[\hat{V}_{m}\right] = \int_{a}^{b} \left(\sum_{i} 
	\varphi^{2}_{i,m}(x)\right) s(x) dx
	\leq (k+1)m\| s \|_{\infty},
\]
since the $\varphi_{i,m}$'s are orthonormal. On the other hand, 
by Chapter 2 (10.1) in \cite{Devore}, 
if $s\in\mathcal{B}^{\alpha}_{\infty}\left(\bbl^{p}([a,b])\right)$,
there is a polynomial $q\in\calS_{m}^{k}$ such that
\[
	\|s - q\|_{\bbl^{p}}\leq c_{[\alpha]}
	| s |_{\mathcal{B}^{\alpha}_{\infty}\left(\bbl^{p}\right)}
	(b-a)^{\alpha} m^{-\alpha}.
\]
Thus,
\[
	\|s - s_{m}^{\bot} \|\leq 
	c_{[\alpha]}(b-a)^{\frac{1}{2}-\frac{1}{p}+\alpha} 
	| s |_{\mathcal{B}^{\alpha}_{\infty}\left(\bbl^{p}\right)}
	m^{-\alpha}.
\]
By (\ref{PseudoOracleIneq})), there is a constant $M$ 
(depending on $C$, $c$, $c'$, $c''$, $\alpha$, $k$, $b-a$, $p$, 
$| s |_{\mathcal{B}^{\alpha}_{\infty}\left(\bbl^{p}\right)}$, and 
$\| s \|_{\infty}$), for which
\[
	\bbe \left[ \| s - \tilde s_{_{T}} \|^{2}\right]
	\leq M
	\inf_{m\in\calM_{T}}
	\left\{m^{-2\alpha}+\frac{m}{T}\right\} + \frac{C'}{T}.
\]
It is not hard to see that, for large enough $T$, the infimum on 
the above right hand side is 
$O_{\alpha}(T^{-2\alpha/(2\alpha+1)})$ (where $O_{\alpha}$ means
that the term depends only on $\alpha$).
Since $M$ is monotone in
$| s |_{\mathcal{B}^{\alpha}_{\infty}\left(\bbl^{p}\right)}$ and 
$\| s \|_{\infty}$, (\ref{LongRunRiskppe}) is verified.
\hfill $\Box$


\noindent{\bf Verification of Remark \ref{DmInOB}:} 
Suppose that $D_{m}$ is finite, and thus each $f \in S$, 
with $\| f \|=1$ is bounded. 
It follows using Lagrange multipliers that, for each $x\in D$,
\[
	D(x) \equiv \sup\left\{ \big| 
	\sum_{i=1}^{d_{m}} c_{i} \varphi_{i}(x)\big|^{2}:  
	\sum_{i=1}^{d_{m}} 		
	c_{i}^{2} =1\right\}=\sum_{i=1}^{d_{m}} \varphi_{i}^{2}(x).
\]
Since $D_{m} \geq D(x)$ for every $x\in D$, we obtain 
$D_{m}\geq\|\sum_{i=1}^{d_{m}}\varphi^{2}_{i}\|_{\infty}$. 
On the other hand, for every $\varepsilon>0$, 
there are $b_{1},\dots,b_{n}$ satisfying
$\sum_{i=1}^{d_{m}} b_{i}^{2} =1$ and an $x \in D$ such that 
\[
	D_{m}-\varepsilon <  
	\big| \sum_{i=1}^{d_{m}} b_{i} \varphi_{i}(x) 
	\big |^{2} \leq D(x)
	= \sum_{i=1}^{d_{m}} \varphi_{i}^{2}(x)
	\leq  \big\|\sum_{i=1}^{d_{m}}\varphi^{2}_{i}\big\|_{\infty}.
\]
Letting $\varepsilon \rightarrow 0$, it follows that $D_{m} =  	
\big\|\sum_{i=1}^{d_{m}}\varphi^{2}_{i}\big\|_{\infty}$. \hfill$\Box$

\noindent{\bf Proof of Lemma \ref{BasicIneq}:} 
Clearly, $\gamma_{D}$ as defined by (\ref{Contrast}) can be 
written as 
\[
	\bbe\gamma_{D}(f) 
	= \| f \| ^{2} - 2 \langle f, s_{D}\rangle - 2 \nu_{D}(f)
		 =  \| f - s_{D} \| ^{2} - \| s_{D} \|^{2} - 2 \nu_{D}(f).
	\]
	By the very definition of $\tilde{s}$ as the penalized projection estimator and by Remark \ref{Remark_1},
	\[
		\gamma_{D}(\tilde{s}) + {\rm pen}(\hat{m}) \leq 	\gamma_{D}(\hat{s}_{m}) + {\rm pen}(m) 
		\leq 	\gamma ( s^{\bot}_{m}) + {\rm pen}(m),
	\] 
	for any $m\in{\cal M}$. Using the previous two equations:
	\begin{align}
		\|\tilde{s}  - s_{D} \| ^{2} &= \gamma_{D}(\tilde{s})+ \| s_{D} \|^{2} + 2 \nu_{D}(\tilde{s}) 
		\nonumber \\
		 &\leq \gamma ( s^{\bot}_{m}) + \| s_{D} \|^{2} + 2 \nu_{D}(\tilde{s}) +{\rm pen}(m) - 
		 {\rm pen}(\hat{m}) \nonumber \\	
		 &= \|s_{m}^{\bot} - s_{D}\|^{2}+  2 \nu_{D}(\tilde{s} -s_{m}^{\bot}) +{\rm pen}(m) - 
		 {\rm pen}(\hat{m}). \nonumber
	\end{align}
	Finally, notice that $\nu_{D}(\tilde{s} -s_{m}^{\bot})=\nu_{D}(\tilde{s} -s_{\hat{m}}^{\bot}) +  
	\nu_{D}(s_{\hat{m}}^{\bot} -s_{m}^{\bot}) $ and $\nu_{D}(\hat{s}_{m} -s^{\bot}_{m}) = \chi_{m}^{2}$. 
	\hfill $\Box$
	
\noindent{\bf Verification of  inequality (\ref{BoundForMean}):}
	Notice just that for any $a, b, \varepsilon>0$:
	\begin{equation}\label{Eq1PBoundForMean}
		 a-\sqrt{2ab}-\frac{1}{3}b \geq \frac{a}{1+\varepsilon}
		 -\left( \frac{1}{2\varepsilon}+\frac{5}{6}\right)b.
	\end{equation}
	Evaluating the integral in (\ref{ConctrInqForIntgr}) for $-f$, we can write 
	\[
		\bbp \left[ \int_{{\rm X}} f(x) N(dx) \geq \int_{{\rm X}} f(x) \mu(dx) -
		 \| f \|_{\mu}\sqrt{2 u} - \frac{1}{3} \| f \|_{\infty} u \right]\geq 1- e^{-u}.
	\]
Using that $ \| f \|_{\mu}^{2} \leq \| f \|_{\infty} 
\int_{{\rm X}} |f(x)| \mu(dx)$ and
	 (\ref{Eq1PBoundForMean}),
	\[
		\bbp \left[ \int_{{\rm X}} f(x) N(dx) \geq \frac{1}{1+\varepsilon}\int_{{\rm X}} f(x) \mu(dx) 
		-\left(\frac{1}{2\varepsilon}+\frac{5}{6}\right) \|f\|_{\infty}u \right]\geq 1- e^{-u},
	\]
which is precisely inequality (\ref{BoundForMean}).\hfill$\Box$
	
\noindent{\bf Proof of Lemma \ref{AuxLemma}:}\\
	Let $Z^{+}$ be the positive part of $Z$. First, 
	\[
		\bbe\left[Z\right]\leq\bbe\left[Z^{+}\right] = \int_{0}^{\infty} \bbp[ Z >x] dx.
	\]
	Since $h$ is continuous and strictly increasing, $\bbp[ Z >x]  \leq K \exp(-h^{-1}(x))$, where $h^{-1}$ is the inverse of $h$. Then, changing variables to $u=h^{-1}(x)$,
	\[
		\int_{0}^{\infty} \bbp[ Z >x] dx \leq K \int_{0}^{\infty} e^{-h^{-1}(x)} dx
		= K \int_{0}^{\infty} e^{u} h'(u)du.
	\]
	Finally, an integration by parts yields  $\int_{0}^{\infty} e^{u} h'(u)du=\int_{0}^{\infty} h(u) e^{-u} du$.
	\hfill $\Box$

\noindent{\bf Proof of Proposition \ref{ApprProjAsymt}:}\\
	From the orthonormality  property,
	\begin{align}
		\bbe\left[\|\hat{s}^{n}_{m} - s^{\bot}_{m}\|^{2}\right]&= 	
	       \sum_{i=1}^{d_{m}} \bbe\left[\left(\hat\beta^{n}_{i,m}- \beta_{i,m}\right)^{2}\right] \nonumber \\
	      & =  \sum_{i=1}^{d_{m}} \left\{{\rm Var}\left(\hat\beta^{n}_{i,m}\right) 
	       +\left(\bbe\left[\hat\beta^{n}_{i,m}\right] - \beta_{i,m}\right)^{2}\right\}.\nonumber
	\end{align}
By remark \ref{AsympMeanVarAppxInt}, 
	\begin{align*}
	\lim_{n\rightarrow \infty} 
	\bbe \left[ I_{n}(\varphi_{i,m}) \right]&=
	T \int_{\bbr_{0}} \varphi(x) s(x)\eta(dx),\\
	 \lim_{n\rightarrow \infty} {\rm Var} 
	\left[ I_{n}(\varphi_{i,m}) \right]&=
	T \int_{\bbr_{0}} \varphi_{i.m}^{2}(x) s(x)\eta(dx).
	\end{align*}
	Then, (\ref{AsymRisk}) is true from (\ref{Chi}) and (\ref{RiskPEonOP}). 
	The second statement in the proof is straightforward since
	\[
		\bbe\left[\|\hat{s}^{n}_{m} - s\|^{2}\right]= 
		\bbe\left[\|\hat{s}^{n}_{m} - s^{\bot}_{m}\|^{2}\right]
		+\|{s}^{\bot}_{m} - s\|^{2}.
	\]
	
	\hfill $\Box$


\section{Figures}

\newpage
\noindent
{\bf Acknowledgments:}  It is a pleasure to thank 
P.~Reynaud-Bouret for very helpful discussions.
\bibliographystyle{plain}
\bibliography{Bibliography}

\end{document}